\documentclass[11pt]{article}
\usepackage{amsmath}
\usepackage{amsfonts}
\usepackage{latexsym}
\usepackage{amssymb}
\usepackage{graphicx}
\usepackage{amsthm}
\usepackage{chemarrow}
\usepackage{fullpage}
\usepackage{framed}

\title{On the Mathematics of the Law of Mass Action}
\author{\small Leonard Adleman, Manoj Gopalkrishnan, Ming-Deh Huang, \\ \small Pablo Moisset, Dustin Reishus}
\begin{document}
\maketitle
\theoremstyle{plain}
\newtheorem{theorem}{\bf{Theorem}}[section]
\newtheorem{lemma}[theorem]{\bf{Lemma}}
\newtheorem{corollary}[theorem]{\bf{Corollary}}
\newtheorem{proofl}[theorem]{\bf{Proof}}
\newtheorem{open}{\bf{Open Problem}}

\theoremstyle{definition}
\newtheorem{definition}{\bf{Definition}}[section]

\theoremstyle{remark}
\newtheorem{example}{\bf{Example}}[section]
\newtheorem{algorithm}{\bf{Algorithm}}[section]
\newtheorem{notation}[definition]{\bf{Notation}}
\newtheorem{fact}{\bf{Fact}}
\newtheorem{note}[definition]{\bf{Note}}
\newcommand{\CE}{{\cal E}}
\newcommand{\MC}{\mathbb{C}}
\newcommand{\MR}{\mathbb{R}}
\newcommand{\MZ}{\mathbb{Z}}
\newcommand{\MN}{\mathbb{N}}
\newcommand{\iton}{i=1,2,\ldots,n}
\newcommand{\jtom}{j=1,2,\ldots, m}
\newcommand{\MP}{\mathbb{R}_{>0}}
\newcommand{\bs}{\boldsymbol}
\newcommand{\op}{\operatorname}
\begin{abstract}
In 1864, Waage and Guldberg formulated the ``law of mass action.'' Since that time, chemists, chemical engineers, physicists and mathematicians have amassed a great deal of knowledge on the topic. In our view, sufficient understanding has been acquired to warrant a formal mathematical consolidation.  A major goal of this consolidation is to solidify the mathematical foundations of mass action chemistry --- to provide precise definitions, elucidate what can now be proved, and indicate what is only conjectured. In addition, we believe that the law of mass action is of intrinsic mathematical interest and should be made available in a form that might transcend its application to chemistry alone. We present the law of mass action in the context of a dynamical theory of sets of binomials over the complex numbers.
\end{abstract}
\section{Introduction}\label{es:sec:intro}
The study of mass action kinetics dates back at least to 1864, when Waage and Guldberg~\cite{Waage} formulated the ``law of mass action.'' Since that time, a great deal of knowledge on the topic has been amassed in the form of empirical facts, physical theories and mathematical theorems by chemists, chemical engineers, physicists and mathematicians. In recent years, Horn and Jackson~\cite{horn72general}, and Feinberg~\cite{feinberg95existence} have made significant mathematical contributions, and these have guided our work.

It is our view that a critical mass of knowledge has been obtained, sufficient to warrant a formal mathematical consolidation.  A major goal of this consolidation is to solidify the mathematical foundations of this aspect of chemistry --- to provide precise definitions, elucidate what can now be proved, and indicate what is only conjectured.  In addition, we believe that the law of mass action is of intrinsic mathematical interest and should be made available in a form that might transcend their application to chemistry alone.

To make the law of mass action available for consideration by researchers in areas other than chemistry, we present mass action kinetics in a new form, which we call event-systems. Our formulation begins with the observation that systems of chemical reactions can be represented by sets of binomials. This gives us an opportunity to extend the law of mass action to arbitrary sets of binomials. Once this extension is made, there is no reason to restrict ourselves to binomials with real coefficients. Hence, we are led to a dynamical theory of sets of binomials over the complex numbers. Possible mathematical applications of this theory include:
\begin{enumerate}
\item Binomials are objects of intrinsic mathematical interest~\cite{eisenbud1994binomial}. For example, they occur in the study of toric varieties, and hence in string theory. With each set of binomials over the complex numbers, we associate a corresponding system of differential equations. Ideally, this dynamical viewpoint will help advance the theory of binomials, and enhance our understanding of their associated algebraic sets.

\item When we extend the study of the law of mass action to sets of binomials over the complex numbers, we can consider reactions that involve complex rates, complex concentrations, and move through complex time. Extending to the complex numbers gives us direct access to the powerful theorems of complex analysis. Though this clearly transcends conventional chemistry, it may have applications in pure mathematics.

    For example, in ongoing work, we seek to exploit an analogy between number theory and chemistry, where atoms are to molecules as primes are to numbers. We associate a distinct species with each natural number. Then each multiplication rule $m\times n = mn$ is encoded by a reaction where the species corresponding to the number $m$ reacts with the species corresponding to the number $n$ to form the species corresponding to the number $mn$. With an appropriate choice of specific rates of reactions the resulting event-system has the property that the sum of equilibrium concentrations of all species at complex temperature $s$ is the value of the Riemann zeta function at $s$. We hope to pursue this approach to study questions related to the distribution of the primes.
\item Systems of linear differential equations are well understood. In contrast, systems of ordinary non-linear differential equations can be notoriously intractable. Differential equations that arise from event-systems lie somewhere in between --- more structured than arbitrary non-linear differential equations, but more challenging than linear differential equations. As such, they appear to be an important new class for consideration in the theory of ordinary differential equations.
\end{enumerate}

In addition to their use in mathematics, event-systems provide a vehicle by which ideas in algebraic geometry may be made readily available to the study of mass action kinetics. As such, they may help solidify the foundations of this aspect of chemistry. We expand on this in Section~\ref{conclusion}.

Part of our motivation for this research comes from the emerging field of nanotechnology. To quote from \cite{adleman99toward}, ``Self-assembly is the ubiquitous process by which objects autonomously assemble into complexes. Nature provides many examples: Atoms react to form molecules. Molecules react to form crystals and supramolecules. Cells sometimes coalesce to form organisms. Even heavenly bodies self-assemble into astronomical systems. It has been suggested that self-assembly will ultimately become an important technology, enabling the fabrication of great quantities of small objects such as  computer circuits\dots\@ Despite its importance, self-assembly is poorly understood.'' Hopefully, the theory of event-systems is a step towards understanding this important process.

The paper is organized as follows:

In Section~\ref{sec:definitions}, we present the basic mathematical notations and definitions for the study of event-systems.

In Section~\ref{sec:finite}, and all of the sections that follow, we restrict to finite event-systems. Theorem~\ref{thm:consclass} demonstrates that the stoichiometric coefficients give rise to flow-invariant affine subspaces --- ``conservation classes.''

In Section~\ref{sec:physical}, and all of the sections that follow, we restrict to ``physical event-systems.'' Though we have defined event-systems over the complex numbers, in this paper we focus on consolidating results from the mass action kinetics of reversible chemical reactions. Physical event-systems capture the idea that the specific rates of chemical reactions are always positive real numbers. The main result of this section is Theorem~\ref{thm:stays_positive}, which demonstrates that for physical event-systems, if initially all concentrations are non-negative, then they stay non-negative for all future real times so long as the solution exists. Further, the concentration of every species whose initial concentration is positive, stays positive.

In Section~\ref{sec:natural}, and all the sections that follow, we restrict to ``natural event-systems.'' Natural event-systems capture the concept of detailed balance from chemistry. In Theorem~\ref{NatL}, we give four equivalent characterizations of natural event-systems; in particular, we show that natural event-systems are precisely those physical event-systems that have no ``energy cycles.'' In Theorem~\ref{LyapunovExists}, following Horn and Jackson~\cite{horn72general}, we show that natural event-systems have associated Lyapunov functions. This theorem is reminiscent of the second law of thermodynamics. The main result of this section is Theorem~\ref{AsymptoticStability}, which establishes that for natural event-systems, given non-negative initial conditions:
\begin{enumerate}
\item Solutions exist for all forward real times.
\item Solutions are uniformly bounded in forward real time.
\item All positive equilibria satisfy detailed balance.
\item Every conservation class containing a positive point also contains exactly one positive equilibrium point.
\item Every positive equilibrium point is asymptotically stable relative to its conservation class.
\end{enumerate}
For systems of reversible reactions that satisfy detailed balance, must concentrations approach equilibrium? We believe this to be the case, but are unable to prove it. In 1972, an incorrect proof was offered~\cite[Lemma 4C]{horn72general}. This proof was retracted in 1974~\cite{horn74dynamics}. To the best of our knowledge, this question in mass action kinetics remains unresolved~\cite[p.\hspace{3pt}10]{sontag01structure}. We pose it formally in Open~Problem~\ref{MainTheorem}, and consider it the fundamental open question in the field.

In Section~\ref{sec:atomic}, we introduce the notion of ``atomic event-systems.'' As the name suggests, this is an attempt to capture mathematically the atomic hypothesis that all species are composed of atoms. The main theorem of this section is Theorem~\ref{AtomicMainTheorem}, which establishes that for natural, atomic event-systems, solutions with positive initial conditions asymptotically approach positive equilibria. Hence, Open~Problem~\ref{MainTheorem} is resolved in the affirmative for this restricted class of event-systems.

\section{Basic Definitions and Notation}\label{sec:definitions}
Before formally defining event-systems, we give a very brief, informal introduction to chemical reactions. All reactions are assumed to take place at constant temperature in a well-stirred vessel of constant volume.

Consider \[A+2B\autorightleftharpoons{\small$\sigma$}{\small$\tau$}C.\]
This chemical equation concerns the reacting species $A,B$ and $C$. In the forward direction, one mole of $A$ combines with two moles of $B$ to form one mole of $C$. The symbol ``$\sigma$'' represents a real number greater than zero. It denotes, in appropriate units, the rate of the forward reaction when the reaction vessel contains one mole of $A$ and one mole of $B$. It is called the specific rate of the forward reaction. In the reverse direction, one mole of $C$ decomposes to form one mole of $A$ and two moles of $B$. The symbol ``$\tau$'' represents the specific rate of the reverse reaction. Chemists typically determine specific rates empirically. Though irreversible reactions (those with $\sigma=0$ or $\tau=0$) have been studied, they will not be considered in this paper.

Inspired by the law of mass action, we introduce a multiplicative notation for chemical reactions, as an alternative to the chemical equation notation. In our notation, each chemical reaction is represented by a binomial. Consider the following examples. On the left are chemical equations. On the right
are the corresponding binomials.
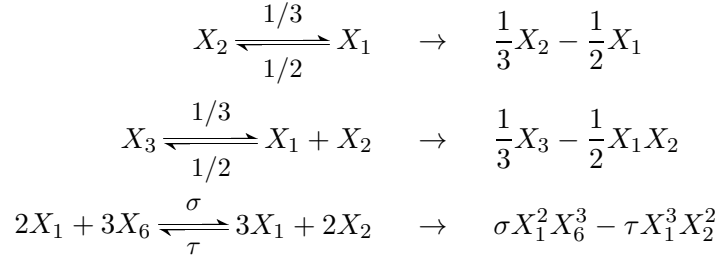
\begin{eqnarray*}
X_2 \autorightleftharpoons{\small 1/3}{\small 1/2} X_1 \ \
&\rightarrow&\ \
 \frac{1}{3}X_2 - \frac{1}{2}X_1\\
X_3\autorightleftharpoons{\small 1/3}{\small 1/2}X_1+X_2\ \
&\rightarrow&\ \
\frac{1}{3}X_3 - \frac{1}{2}X_1 X_2  \\
2X_1+ 3 X_6\autorightleftharpoons{\small $\sigma$}{\small $\tau$}3 X_1+2X_2 \ \ &\rightarrow&\ \ \sigma X^2_1 X^3_6 - \tau X^3_1
X^2_2
\end{eqnarray*}

Our notation leads us to view every set of binomials over an
arbitrary field $\mathbb{F}$ as a formal system of reversible
reactions with specific rates in $\mathbb{F}\setminus\{0\}$. For our present purposes, we will restrict our attention to binomials over the complex numbers. With this in mind, we now define our notion of event-system.

\begin{notation}
Let $\MC_{\infty}=\bigcup_{n=1}^{\infty}\MC[X_1, X_2, \cdots, X_n]$. A monic monomial of $\MC_{\infty}$ is a product of the form $\prod_{i=1}^{\infty}X_i^{e_i}$ where the $e_i$ are non-negative integers all but finitely many of which are zero. We will write ${\mathbb M}_\infty$ to denote the set of all monic
monomials of $\MC_{\infty}$. More generally, if
$S\subset\{X_1,X_2,\cdots\}$, we let $\MC[S]$ be the ring of polynomials
with indeterminants in $S$ and we let
$\mathbb{M}_S=\mathbb{M}_\infty\cap\MC[S]$ (i.e. the monic monomials in
$\MC[S]$).

If $n\in\MZ_{>0}$, $p\in\MC[X_1,X_2,\cdots,X_n]$, and
${\bs a}=\langle a_1,a_2,\cdots,a_n\rangle\in\MC^n$ then, as is usual,
we will let $p({\bs a})$ denote the value of $p$ on argument ${\bs a}$.

Given two monic monomials $M=\prod_{i=1}^{\infty}X_i^{e_i}$ and
$N=\prod_{i=1}^{\infty}X_i^{f_i}$ from $\mathbb{M}_\infty$, we will say
$M$ {\em precedes} $N$ (and we will write $M\prec N$) iff $M \neq N$ and
for the least $i$ such that $e_i \not = f_i$, $e_i<f_i$.
\end{notation}

It follows that $1$ is a monic monomial of $\ \MC_{\infty}$ and that
each element of $\MC_{\infty}$ is a $\MC$-linear combination of finitely
many monic monomials. We will be particularly concerned with the set of
binomials ${\mathbb B}_{\infty}=\{\sigma M+\tau N\mid\sigma,\tau\in{\mathbb
C} \setminus \{0\}$ and $M, N$ are distinct monic monomials of $\MC_\infty\}$.

\begin{definition}[Event-system]\label{def:event-system}
An event-system $\CE$ is a nonempty subset of ${\mathbb B}_{\infty}$.
\end{definition}

If $\CE$ is an event-system, its elements will be called
``$\CE$-events'' or just ``events.'' Note that if $\sigma M + \tau
N$ is an event then $M\neq N$.

Our map from chemical equations to events is as follows. A
chemical equation
\begin{eqnarray*}
\sum_i a_i X_i \autorightleftharpoons{\small $\sigma$}{\small $\tau$}
\sum_j b_j X_j \mbox{\ \ goes to:}
\end{eqnarray*}
\begin{eqnarray*}
&1.&\ \ \sigma \prod_i X^{a_i}_i - \tau \prod_j X^{b_j}_j \mbox{\ \ if\ \ }
\prod_i X^{a_i}_i \prec \prod_j X^{b_j}_j\\
\mbox{or\ }&2.&\ \ \tau \prod_j X^{b_j}_j - \sigma \prod_i X^{a_i}_i
\mbox{\ \ if\ \ } \prod_j X^{b_j}_j \prec \prod_i X^{a_i}_i
\end{eqnarray*}

For example:
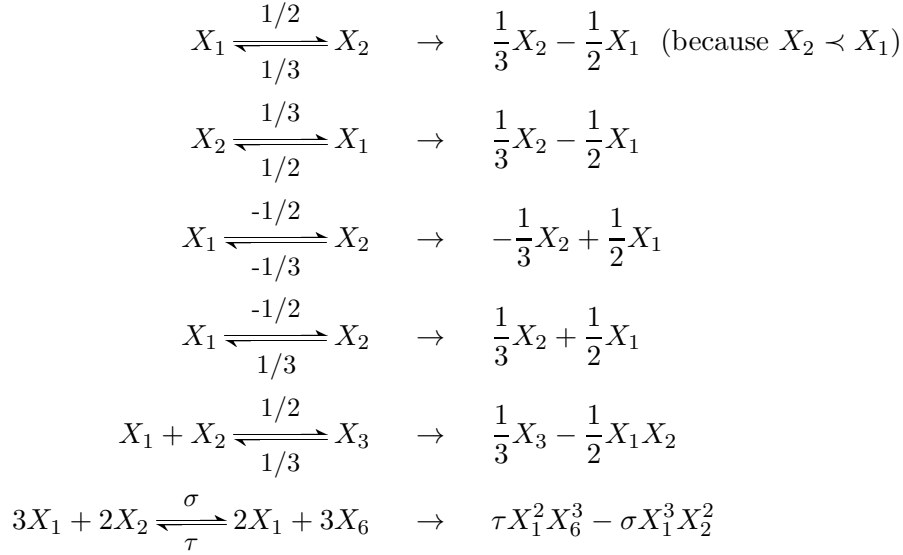
\begin{eqnarray*}
X_1 \autorightleftharpoons{\small 1/2}{\small 1/3} X_2 \ \ &\rightarrow&\ \
 \frac{1}{3}X_2 - \frac{1}{2}X_1  \mbox{\ \ (because $X_2\prec X_1$)}\\
X_2 \autorightleftharpoons{\small 1/3}{\small 1/2} X_1 \ \ &\rightarrow&\ \
 \frac{1}{3}X_2 - \frac{1}{2}X_1  \\
X_1\autorightleftharpoons{\small -1/2}{\small -1/3} X_2 \ \ &\rightarrow&\ \
-\frac{1}{3}X_2 +\frac{1}{2}X_1  \\
X_1 \autorightleftharpoons{\small -1/2}{\small 1/3} X_2 \ \ &\rightarrow&\ \
 \frac{1}{3}X_2 +\frac{1}{2}X_1 \\
X_1+X_2\autorightleftharpoons{\small 1/2}{\small 1/3}X_3\ \ &\rightarrow&\ \
\frac{1}{3}X_3 - \frac{1}{2}X_1 X_2  \\
3 X_1+2X_2\autorightleftharpoons{\small $\sigma$}{\small $\tau$}2 X_1+ 3 X_6
\ \ &\rightarrow&\ \
\tau X^2_1 X^3_6 - \sigma X^3_1 X^2_2
\end{eqnarray*}

Note that our order of monomials is arbitrary. Any linear order
would do. The order is necessary to achieve a one-to-one map from chemical reactions to events.

Our definition of event-systems allows for an infinite number of reactions, and an infinite number of reacting species. Indeed, polymerization reactions are commonplace in nature and, in principle, they are capable of creating arbitrarily long polymers (for example, DNA molecules).

The next definition introduces the notion of systems of reactions for which the number of reacting species is finite.

\begin{definition}[Finite-dimensional event-system]
An event-system $\CE$ is {\em finite-dimensional} iff there
exists an $n\in \mathbb{Z}_{>0}$ such that
$\CE\subset\MC[X_1,X_2,\cdots,X_n]$.
\end{definition}

\begin{definition}[Dimension of event-systems]
Let $\CE$ be a finite-dimensional event-system. Then the least $n$ such that
$\CE\subset\MC[X_1,X_2,\cdots,X_n]$ is the {\em dimension} of $\CE$.
\end{definition}

\begin{definition}[Physical event, Physical event-system] A binomial
$e\in\mathbb{B_{\infty}}$ is a {\em physical event} iff there exist
$\sigma,\tau\in\MR_{>0}$ and $M$, $N\in{\mathbb M}_\infty$ such that $M\prec N$
and $e=\sigma M -
\tau N$. An event-system $\CE$ is {\em physical} iff each
$e\in\CE$ is physical.
\end{definition}

Chemical reaction systems typically have positive real forward and
backward rates. Physical event-systems generalize this notion.

\begin{definition}$\label{Point}$
Let $n\in\MZ_{>0}$. Let $\bs{\alpha}=\langle \alpha_1,\alpha_2,\ldots,\alpha_n \rangle\in\mathbb{C}^n$.
\begin{enumerate}
\item $\bs{\alpha}$ is a {\em non-negative point}
iff for $\iton$, $\alpha_i\in\MR_{\geq 0}$.

\item $\bs{\alpha}$ is a {\em positive point}
iff for $\iton$, $\alpha_i\in\MR_{>0}$.

\item $\bs{\alpha}$ is a {\em $z$-point} iff
there exists an $i$ such that $\alpha_i=0$.
\end{enumerate}
\end{definition}

In chemistry, a system is said to have achieved detailed balance when it is at a point where the net flux of each reaction is zero. Given the corresponding event-system, points of detailed balance corresponds to points where each event evaluates to zero, and vice versa. We call such points ``strong equilibrium points.''

\begin{definition}[Strong equilibrium point]$\label{Strong
Equilibrium Point}$ Let $\CE$ be a finite-dimensional event-system
of dimension $n$. $\bs{\alpha}\in\MC^n$ is a {\em strong
$\CE$-equilibrium point} iff for all $e\in\CE$, $e(\bs{\alpha})=0$.
\end{definition}

In the language of algebraic geometry, when $\CE$ is a finite-dimensional event-system, its corresponding algebraic set is precisely the set of its strong $\CE$-equilibrium points.

It is widely believed that all ``real'' chemical reactions achieve detailed balance. We now introduce natural event-systems, a restriction of finite-dimensional, physical event-systems to those that can achieve detailed balance.

\begin{definition}[Natural event-system]$\label{Natural event-systems}$ A finite-dimensional event-system $\CE$ is {\em natural}
iff it is physical and there exists a positive strong
$\CE$-equilibrium point.
\end{definition}

Our next goal is to introduce atomic event-systems: finite-dimensional event-systems obeying the atomic hypothesis that all species are composed of atoms. Towards this goal, we will define a graph for each finite-dimensional event-system. The vertices of this graph are the monomials from $\mathbb{M_\infty}$ and the edges are determined by the events. If a weight $r$ is assigned to an edge, then $r$ represents the energy released when a reaction corresponding to that edge takes place. For the purpose of defining atomic event-systems, the reader may ignore the weights; they are included here for use elsewhere in the paper (Definition~\ref{EnergyCycle}).

Though graphs corresponding to systems of chemical reactions have been defined elsewhere (e.g. ~\cite{feinberg95existence}, \cite[p.\hspace{3pt}10]{sontag01structure}), it is important to note that these definitions do not coincide with ours.

\begin{definition}[Event-graph$\label{EventGraph}$]
Let $\CE$ be a finite-dimensional event-system. The event-graph
$G_\CE = \langle V, E, w\rangle$ is a weighted, directed
multigraph such that:
\begin{enumerate}
\item $V=\mathbb{M_\infty}$
\item For all $M_1$, $M_2\in\mathbb{M_\infty}$, for all $r\in\MC$, \\$\langle M_1, M_2\rangle\in E$
and $r\in w\left(\langle M_1, M_2\rangle\right)$ iff \\there exist
$e\in\CE$ and $\sigma, \tau\in\MC$ and $M,N,T\in\mathbb{M_\infty}$
such that $e=\sigma M + \tau N$ and $M \prec N$ and either \\(a)
$M_1 = TM$ and $M_2=TN$ and
$r=\ln\left(-\frac{\sigma}{\tau}\right)$ or
\\(b) $M_1 = T N$ and $M_2=TM$ and
$r=-\ln\left(-\frac{\sigma}{\tau}\right)$
\end{enumerate}
\end{definition}

Notice that two distinct weights $r_1$ and $r_2$ could be assigned to a single edge. For example, let $\CE=\{X_1X_2-2X_1^2,X_2-5X_1\}$. Consider the edge in $G_\CE$ from the monomial $X_1^2$ to the monomial $X_1X_2$. Weight $\ln 2$ is assigned to this edge due to the event $X_1X_2-2X_1^2$, with $T=1$. Weight $\ln 5$ is also assigned to this edge due to the event $X_2-5X_1$, with $T=X_1$.

\begin{definition} Let $\CE$ be a finite-dimensional event-system.
For all $M \in \mathbb{M}_\infty$, the \emph{connected component}
of $M$, denoted $C_\CE(M)$, is the set of all $N \in
\mathbb{M}_\infty$ such that there is a path in $G_\CE$ from $M$
to $N$.
\end{definition}

It follows from the definition of ``path'' that every monomial belongs to its connected component.

\begin{definition}[Atomic event-system]\label{atomic} Let $\CE$ be a finite-dimensional event-system of dimension $n$.
Let $S = \{X_1,X_2,\cdots,X_n\}$. Let $A_\CE=\big\{X_i\in S\mid
C_\CE(X_i)=\{X_i\}\big\}$. $\CE$ is \emph{atomic} iff for all $M \in
\mathbb{M}_S$, $C(M)$ contains a unique monomial in
$\mathbb{M}_{A_\CE}$.
\end{definition}

If $\CE$ is atomic then the members of $A_\CE$ will be called {\em the atoms of $\CE$}. It follows from the definition that in atomic event-systems, atoms are not decomposable, non-atoms are uniquely decomposable into atoms and events preserve atoms.

Since the set $\mathbb{M}_{\{X_1,X_2\dots,X_n\}}$ is infinite, it is not possible to decide whether $\CE$ is atomic by exhaustively checking the connected component of every monomial in $\mathbb{M}_{\{X_1,X_2\dots,X_n\}}$. The following is sometimes helpful in deciding whether a finite-dimensional event-system is atomic (proof not provided).

Let $\CE$ be an event-system of dimension $n$ with no event of the form $\sigma + \tau N$. Let $B_\CE=\{X_i \mid$~For~all~$\sigma,\tau\in\MC\setminus\{0\}$~and~$N\in\mathbb{M}_\infty: \sigma X_i + \tau N\notin\CE$\}. Then $\CE$ is atomic iff there exist $M_1\in C_\CE(X_1)\cap\mathbb{M}_{B_\CE}, M_2\in C_\CE(X_2)\cap\mathbb{M}_{B_\CE},\dots,M_n\in C_\CE(X_n)\cap\mathbb{M}_{B_\CE}$ such that:
\begin{align}\label{Eq:AtomCons}
\text{For all }\sigma \prod_{i=1}^n X_i^{a_i} - \tau \prod_{i=1}^n X_i^{b_i} \in \CE, \quad\prod_{i=1}^n M_i^{a_i} = \prod_{i=1}^nM_i^{b_i}.
\end{align}

We have shown (proof not provided) that if $\CE$ and $B_\CE$ are as above, and there exist $M_1\in C_\CE(X_1)\cap\mathbb{M}_{B_\CE}, M_2\in C_\CE(X_2)\cap\mathbb{M}_{B_\CE},\dots,M_n\in C_\CE(X_n)\cap\mathbb{M}_{B_\CE}$ and there exists $\sigma \prod_{i=1}^n X_i^{a_i} - \tau \prod_{i=1}^n X_i^{b_i} \in \CE$ such that $\prod_{i=1}^n M_i^{a_i} \neq \prod_{i=1}^nM_i^{b_i}$, then $\CE$ is not atomic. Hence, to check whether an event-system with no event of the form $\sigma + \tau N$ is atomic, it suffices to examine an arbitrary choice of $M_1\in C_\CE(X_1)\cap\mathbb{M}_{B_\CE}, M_2\in C_\CE(X_2)\cap\mathbb{M}_{B_\CE},\dots,M_n\in C_\CE(X_n)\cap\mathbb{M}_{B_\CE}$, if one exists, and check whether (\ref{Eq:AtomCons}) above holds.

\begin{example}
Let $\CE=\{ X_2^2-X_1^2\}$. Then $B_\CE=\{X_1,X_2\}$. Let $M_1 = X_1$ and $M_2 = X_2$. Trivially, $M_1,M_2\in\mathbb{M}_{B_\CE}$, $M_1\in C_\CE(X_1)$ and $M_2\in C_\CE(X_2)$. Consider the event $X_2^2-X_1^2$. Since $M_2^2=X_2^2 \neq X_1^2 = M_1^2$, $\CE$ is not atomic. Note that the event $ X_2^2-X_1^2$ does not preserve atoms.
\end{example}

\begin{example}
Let $\CE=\{X_4^2-X_2,X_5^2-X_3,X_2X_3-X_1\}$. Then $B_\CE=\{X_4,X_5\}$. Let $M_1 = X_4^2X_5^2, M_2 = X_4^2, M_3 = X_5^2, M_4 = X_4, M_5 = X_5$. Clearly these are all in $\mathbb{M}_{B_\CE}$. $X_5^2 - X_3\in\CE$ implies $M_3 \in C_\CE(X_3)$. $X_4^2-X_2\in\CE$ implies $M_2\in C_\CE(X_2)$. Since $(X_1, X_2X_3, X_2X_5^2, X_4^2X_5^2)$ is a path in $G_\CE$, we have $M_1\in C_\CE(X_1)$. For the event $X_4^2-X_2$, we have $M_4^2 = X_4^2 = M_2$. For the event $X_5^2-X_3$, we have $M_5^2 = X_5^2 = M_3$. For the event $X_2X_3-X_1$, we have $M_2M_3 = X_4^2X_5^2 = M_1$. Therefore, $\CE$ is atomic.
\end{example}

Note that it is possible to have an atomic event-system where $A_\CE$ is the empty set. For example:
\begin{example}
Let $\CE=\{1-X_1\}$. In this case, $S=\{X_1\}$ and $\mathbb{M}_S$ is the set $\{1,X_1,X_1^2,X_1^3,\dots\}$. It is clear that $\mathbb{M}_S$ forms a single connected component $C$ in $G_\CE$. Hence, $X_1$ is not in $A_\CE$, and $A_\CE=\emptyset$. $1$ is the only monomial in $\mathbb{M}_{A_\CE}$. Since $1$ is in $C$, $\CE$ is atomic.
\end{example}

\section{Finite Event-systems}\label{sec:finite}
The study of infinite event-systems is embryonic and appears
to be quite challenging. In the rest of this paper only finite
event-systems (i.e., where the set $\CE$ is finite) will be
considered. It is clear that all finite event-systems are
finite-dimensional.

\begin{definition}[Stoichiometric matrix]\label{gamma}
Let $\CE=\{e_1, e_2, \cdots, e_m\}$ be an event-system of
dimension $n$. Let $i\leq n$ and $j\leq m$ be positive integers.
Let $e_j=\sigma M + \tau N$, where $M\prec N$. Then $\gamma_{j,i}$
is the number of times $X_i$ divides $N$ minus the number of times
$X_i$ divides $M$. The {\em stoichiometric matrix} $\Gamma_{\CE}$
of ${\CE}$ is the $m\times n$ matrix of integers
$\Gamma_{\CE}=(\gamma_{j,i})_{m\times n}$.
\end{definition}

\begin{example} Let $e_1=0.5 X_2^5 - 500 X_1X_2^{3}X_7$. Let $\CE=\{e_1\}$. Then
$\gamma_{1,1}=1$, $\gamma_{1,2}=-2$, $\gamma_{1,7}=0$ and for all
other $i$, $\gamma_{1,i}=0$, hence $\Gamma_\CE=\left( \begin{array}{ccccccc}
1 & -2 & 0 & 0 & 0 & 0 & 1\end{array} \right)$.
\end{example}

\begin{definition}\label{def:LoMA}
Let $\CE=\{e_1,\cdots,e_m\}$ be a finite event-system of dimension
$n$. Then:
\begin{enumerate}
\item ${\bs P}_{\CE}$ is the column vector $\langle P_1, P_2, \ldots,
P_n\rangle^T=\Gamma_{\CE}^T\langle e_1,e_2,\ldots,e_m\rangle^T$.
\item Let $\bs{\alpha}\in \MC^n$. Then $\bs{\alpha}$ is an $\CE$-equilibrium
point iff for $\iton:P_{i}(\bs{\alpha})=0$.
\end{enumerate}
\end{definition}

The $P_{i}$'s arise from the
Law of Mass Action in chemistry. For a system of chemical reactions, the $P_i$'s are the right-hand sides of the differential equations that
describe the concentration kinetics. Definition~\ref{def:LoMA} extends the Law of Mass Action to arbitrary event-systems, and hence, arbitrary sets of binomials.

It follows from the definition that for finite event-systems, all
strong equilibrium points are equilibrium points, but the converse need
not be true.

\begin{example}\label{eg:eqmpt1}
Let $e_1 = X_2 - X_1$ and $e_2 = X_2 - 2 X_1$. Let $\CE = \{ e_1, e_2 \}$. Then
$\Gamma_\CE =\left( \begin{array}{cc}
1 & -1\\
1 & -1\end{array} \right)$ and $\bs{P}_\CE = \left( \begin{array}{c}
P_1\\
P_2\end{array} \right) = \left( \begin{array}{c}
2X_2 - 3X_1\\
3X_1 - 2X_2\end{array} \right)$. Therefore $(2,3)$ is an $\CE$-equilibrium point. Since $e_1(2,3)=1$, $(2,3)$ is not a strong $\CE$-equilibrium point.
\end{example}

\begin{example}\label{eg:eqmpt2}
Let $e_1 = 6 - X_1 X_2$ and $e_2 = 2 X_2^2 - 9 X_1$. Let $\CE = \{ e_1, e_2\}$. Then $\Gamma_\CE = \left( \begin{array}{cc}
1 & 1\\
1 & -2\end{array} \right)$ and $\bs{P}_\CE = \left(\begin{array}{c}
P_1\\
P_2\end{array}\right)=\left(\begin{array}{c}
6 - X_1X_2 + 2X_2^2 - \phantom{1}9X_1\\
6 - X_1X_2 - 4X_2^2 + 18 X_1 \end{array}\right)$. The point $(2,3)$ is a strong equilibrium point because $e_1(2,3)=0$ and $e_2(2,3)=0$. Since $P_1(2,3)=e_1(2,3)+e_2(2,3)=0$ and $P_2(2,3)=e_1(2,3) -2e_2(2,3)=0$, the point $(2,3)$ is also an equilibrium point.
\end{example}

The event-system in Example~\ref{eg:eqmpt1} is not natural, whereas the one in Example~\ref{eg:eqmpt2} is. In Theorem~\ref{alleqstrong}, it is shown that if $\CE$ is a finite, natural event-system then all positive $\CE$-equilibrium points are strong $\CE$-equilibrium points.

\begin{definition}[Event-process]
\label{Process} Let $\CE$ be a finite event-system of dimension
$n$. Let $\langle P_1,P_2,\ldots,P_n\rangle^T=\bs{P}_\CE$. Let $\Omega\subseteq\MC$ be a non-empty simply-connected open set.
Let ${\bs f}=\langle f_1,f_2,\cdots,f_n\rangle$ where for
$\iton$, $f_i:\MC\rightarrow\MC$ is defined on $\Omega$. Then
${\bs f}$ is an $\CE$-process on $\Omega$ iff for $\iton$:
\begin{enumerate}
\item $f_i'$ exists on $\Omega$.
\item $f_i'=P_{i}\circ f$ on $\Omega$.
\end{enumerate}
\end{definition}

Note that $\CE$-processes evolve through complex time, and hence generalize the idea of the time-evolution of concentrations in a system of
chemical reactions.

Definition~\ref{Process} immediately
implies that if $\bs{f}=\langle f_1, f_2, \dots,f_n\rangle$ is an $\CE$-process on $\Omega$, then for $\iton$, $f_i$ is holomorphic on $\Omega$. In particular, for each $i$ and all $\alpha\in\Omega$, there is a power series around $\alpha$ that agrees with $f_i$ on a disk of non-zero radius.

Systems of chemical reactions sometimes obey certain conservation laws. For example, they may conserve mass, or the total number of each kind of atom. Event-systems also sometimes obey conservation laws.

\begin{definition}[Conservation law, Linear conservation law]
\label{Conservation Law} Let $\CE$ be a finite event-system of
dimension $n$. A function $g:\MC^n \rightarrow \MC$ is a {\em
conservation law of $\CE$} iff $g$ is holomorphic on $\MC^n$, $g(\langle
0, 0, \cdots, 0\rangle)=0$ and $\nabla g\cdot {\bs P}_{\CE}$ is
identically zero on $\MC^n$. If $g$ is a conservation law of $\CE$ and
$g$ is linear (i.e. $\forall c\in\MC, \forall
{\bs\alpha},{\bs\beta}\in\MC^n,
g(c{\bs\alpha}+{\bs\beta})=cg({\bs\alpha})+g({\bs\beta})$), then $g$ is
a {\em linear conservation law of $\CE$}.
\end{definition}

The event-system described in Example~\ref{eg:eqmpt1} has a linear conservation law $g(X_1,X_2) = X_1 + X_2$. The next theorem shows that conservation laws of $\CE$ are dynamical invariants of $\CE$-processes.

\begin{theorem} For all finite event-systems $\CE$, for all
conservation laws $g$ of $\CE$, for all simply-connected open sets
$\Omega\subseteq \MC$, for all $\CE$-processes ${\bs f}$ on $\Omega$,
there exists $k\in\MC$ such that $g\circ {\bs f} - k$ is
identically zero on $\Omega$.\label{Dynamical Invariant}
\end{theorem}
\begin{proof}
Let $n$ be the dimension of $\CE$. Let $\langle P_1, P_2,\ldots,P_n\rangle^T=\bs{P}_\CE$. For all $t\in \Omega$, by
Definition~\ref{Process}, for $\iton$, $f_i(t)$ and $f'_i(t)$ are
defined. Further, by Definition~\ref{Conservation Law}, $g$ is
holomorphic on $\MC^n$. Hence, $g\circ {\bs f}$ is holomorphic on
$\Omega$. Therefore, by the chain rule, $(g \circ {\bs f})'(t) =
(\nabla g|_{{\bs f}(t)}) \cdot \langle f_1'(t),f_2'(t),
\ldots, f_n'(t)\rangle$. By Definition~\ref{Process}, for all $t\in
\Omega$, $\langle f_1'(t),f_2'(t), \ldots, f_n'(t)\rangle=\langle
P_1({\bs f}(t)),P_2({\bs f}(t)), \ldots, P_n({\bs f}(t))\rangle$. From
these, it follows that $(g \circ {\bs f})'(t) = (\nabla g\cdot {\bs
P}_{\CE})({\bs f}(t))$. But by Definition~\ref{Conservation Law},
$\nabla g\cdot {\bs P}_{\CE}$ is identically zero. Hence, for all $t\in
\Omega$, $(g\circ {\bs f})'(t)=0$. In addition,  $\Omega$ is a
simply-connected open set. Therefore,  by~\cite[Theorem~11]{AH79}, there exists $k\in\MC$ such that $g\circ {\bs
f} -k$ is identically zero on $\Omega$.
\end{proof}

The next theorem shows a way to derive linear conservation laws of an event-system from its stoichiometric matrix.

\begin{theorem}
Let $\CE$ be a finite event-system of dimension $n$. For all
$\bs{v}\in\ker \Gamma_{\CE}$, $\bs{v}\cdot\langle X_1,\cdots,X_n\rangle$ is a
linear conservation law of $\CE$. \label{linear conservation law}
\end{theorem}
\begin{proof}
Let $\Gamma=\Gamma_{\CE}$, then $\ker\Gamma$ is orthogonal to the image
of $\Gamma^T$. By the definition of ${\bs P}={\bs P}_{\CE}$, for all
$\bs{w}\in\MC^n$, ${\bs P}(\bs{w})$ lies in the image of $\Gamma^T$. Hence, for all
$\bs{v}\in\ker\Gamma$, for all $\bs{w}\in\MC^n$, $\bs{v}\cdot {\bs P}(\bs{w})=0$. But $\bs{v}$ is the gradient of $\bs{v}\cdot\langle X_1,\cdots,X_n\rangle$. It now follows from
Definition~\ref{Conservation Law} that $\bs{v}\cdot\langle
X_1,\cdots,X_n\rangle$ is a linear conservation law of $\CE$.
\end{proof}

\begin{definition}[Primitive conservation law]
\label{Primitive} Let $\CE$ be a finite event-system of dimension
$n$. For all $\bs{v}\in \ker\Gamma_{\CE}$, the linear conservation law
$\bs{v}\cdot\langle X_1, X_2,\cdots,X_n \rangle$ is a {\em primitive
conservation law}.
\end{definition}

We can show (manuscript under preparation) that in physical event-systems all linear conservation laws are primitive and, in natural event-systems, all conservation laws arise from the primitive ones.

\begin{definition}[Conservation class, Positive conservation class]
$\label{Conservation Class}$ Let $\CE$ be a finite event-system of
dimension $n$. A coset of $(\ker\Gamma_{\CE})^\perp$ is a {\em
conservation class of $\CE$}. If a conservation class of $\CE$
contains a positive point, then the class is a {\em positive
conservation class of $\CE$}.
\end{definition}

Equivalently, $\bs{\alpha},\bs{\beta}\in\MC^n$ are in the same conservation class if and only if they agree on all primitive conservation laws. Note that if $H$ is a conservation class of $\CE$ then it is closed in $\MC^n$. The following theorem shows that the name ``conservation class'' is appropriate.

\begin{theorem}\label{thm:consclass}
Let $\CE$ be a finite event-system. Let $\Omega\subset\MC$ be a
simply-connected open set containing $0$. Let ${\bs f}$ be an
$\CE$-process on $\Omega$. Let $H$ be a conservation class of $\CE$
containing $\bs{f}(0)$. Then for all $t\in \Omega$, ${\bs f}(t)\in H$.
\end{theorem}
\begin{proof}
Let $\CE$, $\Omega$, ${\bs f}$, $H$ and $t$ be as in the statement
of this theorem. For all $\bs{v}\in \ker\Gamma_{\CE}$, the primitive
conservation law $\bs{v}\cdot\langle X_1, X_2,\cdots,X_n \rangle$ is a
dynamical invariant of ${\bs f}$, from Theorem~\ref{linear
conservation law} and Theorem~\ref{Dynamical Invariant}. Hence,
\[\bs{v}\cdot\langle f_1(0), f_2(0),\cdots,f_n(0) \rangle=\bs{v}\cdot\langle
f_1(t), f_2(t),\cdots,f_n(t) \rangle\] That is, \[\bs{v}\cdot\langle
f_1(0) - f_1(t), f_2(0)-f_2(t),\cdots,f_n(0)-f_n(t)\rangle=0\]
Hence, ${\bs f}(t)-{\bs f}(0)$ is in
$(\ker\Gamma_\CE)^\perp$. By Definition~\ref{Conservation
Class}, ${\bs f}(t)\in H$.
\end{proof} 
\section{Finite Physical Event-systems}\label{sec:physical}
In this section, we investigate finite, physical event-systems --- a generalization of systems of chemical reactions.

It is widely believed that systems of chemical reactions that begin with positive (respectively, non-negative) concentrations will have positive (respectively, non-negative) concentrations at all future times.  This property has been addressed mathematically in numerous papers~\cite[p.~6]{gatermann02family},\cite[Remark~3.4]{feinberg95existence},
\cite[Theorem~3.2]{bernstein99nonnegative},
\cite[Lemma~2.1]{sontag01structure}. The notion of ``system of chemical reactions'' varies between papers.  Several papers have provided no proof, incomplete proofs or inadequate proofs that this property holds for their systems. Sontag~\cite[Lemma~2.1]{sontag01structure} provides a lovely proof of this property for the systems he considers --- zero deficiency reaction networks with one linkage class. We shall prove in Theorem~\ref{thm:stays_positive} that the property holds for finite, physical event-systems. Finite, physical event-systems have a large intersection with the systems considered by Sontag, but each includes a large class of systems that the other does not. We remark that our methods of proof differ from Sontag's, but it is possible that Sontag's proof might be adaptable to our setting.

Lemma~\ref{lem:solution} and Lemma~\ref{invariant} are proved here because they apply to finite, physical event-systems. However, they are only invoked in subsequent sections. Lemma~\ref{lem:solution} relates $\CE$-processes to solutions of ordinary differential equations over the reals. Lemma~\ref{invariant} establishes that if an $\CE$-process defined on the positive reals starts at a real, non-negative point, then its $\omega$-limit set is invariant and contains only real, non-negative points.

The next lemma shows that if two $\CE$-processes evaluate to the same real point on a real argument then they must agree and be real-valued on an open interval containing that argument. The proof exploits the fact that $\CE$-processes are analytic, by considering their power series expansions.

\begin{lemma}
\label{lem:uni_real_local}
Let $\CE$ be a finite, physical event-system of dimension $n$, let
$\Omega,\Omega'\subseteq\MC$ be open and simply-connected, let ${\bs
f}=\langle f_1, f_2, \ldots, f_n\rangle$ be an $\CE$-process on
$\Omega$ and let ${\bs g}=\langle g_1, g_2, \ldots, g_n\rangle$ be an
$\CE$-process on $\Omega'$. If $t_0 \in \Omega\cap\Omega'\cap\MR$ and
${\bs f}(t_0)\in\MR^n$ and ${\bs f}(t_0)={\bs g}(t_0)$, then
there exists an open interval $I\subseteq\MR$ such that $t_0\in I$ and
for all $t\in I$:

\begin{enumerate}
\item \label{uni}   ${\bs f}(t)={\bs g}(t)$.
\item \label{taylor}For $\iton:$ if $\sum_{j=0}^{\infty}c_j(z-t_0)^j$ is the Taylor series expansion of $f_i$ at $t_0$ then for all $j\in\MZ_{\geq 0}$, $c_j\in\MR$.
\item \label{reals} ${\bs f}(t)\in\MR^n$.
\end{enumerate}
\end{lemma}

\begin{proof}
Let $k\in\MZ_{\geq 0}$. By Definition~\ref{Process}, ${\bs f}$ and ${\bs
g}$ are vectors of functions analytic at $t_0$. For $\iton$, let
$f^{(k)}_i$ be the $k^{th}$ derivative of $f_i$ and let ${\bs
f}^{(k)}=\langle f^{(k)}_1,  f^{(k)}_2, \ldots, f^{(k)}_n \rangle$.
Define $g^{(k)}_i$ and ${\bs g}^{(k)}$ similarly. To prove~\ref{uni}, it
is enough to show that for $\iton$, $f_i$ and $g_i$ have the same Taylor
series around $t_0$. Let ${\bs V}_0=\langle X_1, X_2,\ldots,
X_n\rangle$. Let ${\bs V}_k=\op{Jac}({\bs V}_{k-1}) {\bs P}_\CE$ (recall that
if ${\bs H}=\langle h_1(X_1, X_2, \ldots, X_m), h_2(X_1, X_2, \ldots,
X_m), \ldots, h_n(X_1, X_2,\ldots,X_m)\rangle$ is a vector of functions
in $m$ variables then $\op{Jac}({\bs H})$ is the $n\times m$ matrix
$(\frac{\partial h_i}{\partial x_j})$, where $i=1,2,\ldots, n$ and
$j=1,2, \ldots, m$). Let $\langle V_{k,1}, V_{k,2}, \ldots, V_{k,n}
\rangle = {\bs V}_k$. We claim that  ${\bs f}^{(k)}={\bs V}_k \circ {\bs
f}$ on $\Omega$ and ${\bs g}^{(k)}={\bs V}_k \circ {\bs g}$ on $\Omega'$
and for  $\iton$, $V_{k,i} \in \MR[X_1, X_2, \ldots, X_n]$. We prove the
claim by induction on $k$. If $k=0$, the proof is immediate. If $k\geq
1$, on $\Omega$:

\begin{align*}
{\bs f}^{(k)}&=({\bs f}^{(k-1)})'\\
&=({\bs V}_{k-1}\circ {\bs f})'&\text{(Inductive hypothesis)}\\
&=(\op{Jac}({\bs V}_{k-1})\circ {\bs f}) {\bs f}'&\text{(Chain-rule of derivation)}\\
&=(\op{Jac}({\bs V}_{k-1})\circ {\bs f}) ({\bs P}_\CE\circ {\bs f})&\text{(${\bs f}$ is an $\CE$-process)}\\
&=(\op{Jac}({\bs V}_{k-1}) {\bs P}_\CE) \circ {\bs f}\\
&={\bs V}_{k} \circ {\bs f}
\end{align*}

By a similar argument, we conclude that ${\bs g}^{(k)}={\bs V}_{k} \circ
{\bs g}$ on $\Omega'$.  By the inductive hypothesis, ${\bs V}_{k-1}$ is
a vector of polynomials in $\MR[X_1,X_2,\ldots,X_n]$. It follows that
$\op{Jac}({\bs V}_{k-1})$ is an $n\times n$ matrix of polynomials in
$\MR[X_1,X_2,\ldots,X_n]$. Since $\CE$ is physical, ${\bs P}_\CE$ is a
vector of polynomials in $\MR[X_1,X_2,\ldots,X_n]$. Therefore, ${\bs
V}_k=\op{Jac}({\bs V}_{k-1}) {\bs P}_\CE$ is a vector of polynomials in
$\MR[X_1,X_2,\ldots,X_n]$. This establishes the claim.

We have proved that ${\bs f}^{(k)}={\bs V}_{k} \circ {\bs f}$ on
$\Omega$ and ${\bs g}^{(k)}={\bs V}_{k} \circ {\bs g}$ on $\Omega'$.
Since, by assumption, $t_0\in\Omega\cap\Omega'$ and ${\bs f}(t_0)={\bs
g}(t_0)$, it follows that ${\bs f}^{(k)}(t_0)={\bs g}^{(k)}(t_0)$.
Therefore, for $\iton$, $f_i$ and $g_i$ have the same Taylor series
around $t_0$. For $\iton$, let $a_i$ be the radius of convergence of the
Taylor series of $f_i$ around $t_0$. Let $r_{\bs
f}=\min_{i\in\{1,2,\ldots,n\}}a_i$. Define $r_{\bs g}$ similarly. Let
$D\subseteq\Omega\cap\Omega'$ be some non-empty open disk centered at
$t_0$ with radius $r\leq\min(r_{\bs f}, r_{\bs g})$. Since $\Omega$ and
$\Omega'$ are open sets and $t_0\in\Omega\cap\Omega'$, such a disk must
exist. Letting $I=(t_0-r, t_0+r)$ completes the proof of~\ref{uni}.

By assumption, ${\bs f}(t_0)\in\MR^n$, and we have proved that ${\bs
f}^{(k)}={\bs V}_{k} \circ {\bs f}$ and ${\bs V}_k$ is a vector of
polynomials in $\MR[X_1,X_2,\ldots,X_n]$. It follows that ${\bs
f}^{(k)}(t_0)\in\MR^n$. Therefore, for $\iton$, all coefficients in the
Taylor series of $f_i$ around $t_0$ are real. It follows that $f_i$ is
real valued on $I$, completing the proof of~\ref{reals}.
\end{proof}

The next lemma is a kind of uniqueness result. It shows that if two $\CE$-processes evaluate to the same real point at $0$ then they must agree and be real-valued on every open interval containing $0$ where both are defined. The proof uses continuity to extend the result of Lemma~\ref{lem:uni_real_local}.

\begin{lemma}
\label{lem:uni_real_global}
Let $\CE$ be a finite, physical event-system of dimension $n$, let
$\Omega,\Omega'\subseteq\MC$ be open and simply-connected, let ${\bs
f}=\langle f_1, f_2, \ldots, f_n\rangle$ be an $\CE$-process on
$\Omega$ and let ${\bs g}=\langle g_1, g_2, \ldots, g_n\rangle$ be an
$\CE$-process on $\Omega'$. If $0 \in \Omega\cap\Omega'$ and ${\bs
f}(0)\in\MR^n$ and ${\bs f}(0)={\bs g}(0)$, then for all open intervals
$I\subseteq\Omega\cap\Omega'\cap\MR$ such that $0\in I$, for all $t\in
I$, ${\bs f}(t)={\bs g}(t)$ and ${\bs f}(t)\in\MR^n$.
\end{lemma}
\begin{proof}
Assume there exists an open interval
$I\subseteq\Omega\cap\Omega'\cap\MR$ such that $0\in I$ and $B=\{t\in
I\mid {\bs f}(t)\not={\bs g}(t)$ or ${\bs
f}(t)\not\in\MR^n\}\not=\varnothing$. Let $B_P=B\cap\MR_{\geq 0}$ and let
$B_N = B\cap\MR_{<0}$. Note that $B=B_P\cup B_N$, hence,
$B_P\not=\varnothing$ or $B_N\not=\varnothing$. Suppose $B_P\not=\varnothing$
and let $t_P=\inf(B_P)$. By Lemma~\ref{lem:uni_real_local}, there
exists an $\varepsilon\in\MR_{>0}$ such that $(-\varepsilon,
\varepsilon)\cap B=\varnothing$. Hence, $t_P\geq\varepsilon>0$.  By
definition of $t_P$, for all $t\in[0, t_P)$, ${\bs f}(t)={\bs g}(t)$ and
${\bs f}(t)\in\MR^n$. Since ${\bs f}$ and ${\bs g}$ are analytic at
$t_P$, they are continuous at $t_P$. Therefore, ${\bs f}(t_P)={\bs
g}(t_P)$ and ${\bs f}(t_P)\in\MR^n$. By Lemma~\ref{lem:uni_real_local},
there exists an $\varepsilon'\in\MR_{>0}$ such that for all $t\in
(t_P-\varepsilon', t_P+\varepsilon')$, ${\bs f}(t)={\bs g}(t)$ and ${\bs
f}(t)\in\MR^n$, contradicting $t_P$ being the infimum of $B_P$.
Therefore, $B_P=\varnothing$. Using a similar agument, we can prove that
$B_N=\varnothing$. Therefore, $B=\varnothing$, and for all $t\in I$, ${\bs
f}(t_P)={\bs g}(t_P)$ and ${\bs f}(t_P)\in\MR^n$.
\end{proof}

The next lemma is a convenient technical result that lets us ignore the choice of origin for the time variable.

\begin{lemma}\label{time invariance}
Let $\CE$ be a finite, physical event-system of dimension $n$, let
$\Omega, \widetilde{\Omega} \subseteq \MC$ be open and simply connected,
let $\bs{f}=\langle f_1, f_2, \dots, f_n \rangle$ be an
$\CE$-process on $\Omega$ and let $\bs{\tilde{f}}=\langle
\tilde{f}_1, \tilde{f}_2, \ldots, \tilde{f}_n \rangle$ be an
$\CE$-process on $\widetilde{\Omega}$. Let $u \in \Omega$ and
$\tilde{u} \in \widetilde{\Omega}$ and $\bs{\alpha} \in \MR^n$. Let $I \subseteq \MR$ be an open interval. If
\begin{enumerate}
\item $\bs{f}(u) = \bs{\tilde{f}}(\tilde{u}) = \bs{\alpha}$ and
\item $0 \in I$ and
\item for all $s \in I$, $u+s\in \Omega$ and $\tilde{u}+s \in \widetilde{\Omega}$
\end{enumerate}
then for all $t \in I$, $\bs{f}(u+t) = \bs{\tilde{f}}(\tilde{u}+t)$.
\end{lemma}

\begin{proof}
Suppose $\bs{f}(u) = \bs{\tilde{f}}(\tilde{u}) = \bs{\alpha} \in \MR^n$. Let $\Omega_u = \{z \in \MC
\mid u+z \in \Omega\}$ and $\widetilde{\Omega}_{\tilde{u}} = \{z \in \MC \mid
\tilde{u}+z \in \widetilde{\Omega}\}$. Let $\bs{h}=\langle h_1, h_2,
\ldots, h_n \rangle$ where for $\iton$, $h_i:\Omega_u \rightarrow
\MC$ is such that for all $z\in\Omega_u$, $h_i(z) = f_i(u+z)$ and
let $\bs{\tilde{h}}=\langle \tilde{h}_1, \tilde{h}_2, \ldots,
\tilde{h}_n \rangle$ where for $\iton$,
$\tilde{h}_i:\widetilde{\Omega}_{\tilde{u}} \rightarrow \MC$ is such
that for all $z \in \widetilde{\Omega}_{\tilde{u}}$, $\tilde{h}_i(z) =
\tilde{f}_i(\tilde{u}+z)$. Since $u+z$ is differentiable on
$\Omega_u$ and for $\iton$, $f_i$ is differentiable on $\Omega$, it
follows that for $\iton$, $h_i$ is differentiable on $\Omega_u$.
Further, for $\iton$, for all $z \in \Omega_u$, $h_i'(z) = f_i'(u+z)
= \bs{P}_\CE(f_i(u+z)) = \bs{P}_\CE(h_i(z))$, so $\bs{h}$ is an
$\CE$-process on $\Omega_u$. Similarly,
$\bs{\tilde{h}}$ is an $\CE$-process on
$\widetilde{\Omega}_{\tilde{u}}$. Note that $0 \in \Omega_u \cap
\widetilde{\Omega}_{\tilde{u}}$ because $u \in \Omega$ and $\tilde{u}
\in \widetilde{\Omega}$ and that $\bs{h}(0) = \bs{\tilde{h}}(0) = \bs{\alpha}$ because
$\bs{f}(u) = \bs{\tilde{f}}(\tilde{u}) = \bs{\alpha}$. By
Lemma~\ref{lem:uni_real_global}, for all open intervals $I \subseteq \Omega_u \cap
\widetilde{\Omega}_{\tilde{u}} \cap \MR$ such that $0 \in I$, for all $t \in I$, $\bs{h}(t) =
\bs{\tilde{h}}(t)$, so $\bs{f}(u+t) = \bs{\tilde{f}}(\tilde{u}+t)$.
\end{proof}

Because event-systems are defined over the complex numbers, we have access to results from complex analysis. However, there is a considerable body of results regarding ordinary differential equations over the reals. Definition~\ref{def:realprocess} and Lemma~\ref{lem:solution} establish a relationship between $\CE$-processes and solutions to systems of ordinary differential equations over the reals.

\begin{definition}[Real event-process]\label{def:realprocess}
Let $\CE$ be a finite, physical event-system of dimension $n$. Let $\langle P_1,P_2,\dots,P_n\rangle^T = \bs{P}_\CE$. Let
$I\subseteq\MR$ be an interval. Let $\bs{h}=\langle h_1, h_2,
\ldots, h_n\rangle$ where for $\iton$, $h_i:\MR\rightarrow\MR$ is
defined on $I$. Then $\bs{h}$ is a real-$\CE$-process on $I$ iff for
$\iton$:
\begin{enumerate}
\item $h_i'$ exists on $I$.
\item $h_i'=P_{i}\circ \bs{h}$ on $I$.
\end{enumerate}
\end{definition}

\begin{lemma}[All real-$\CE$-processes are restrictions of $\CE$-processes]\label{lem:solution}
Let $\CE$ be a finite, physical event-system of dimension $n$. Let
$I\subseteq\MR$ be an interval. Let $\boldsymbol{h}=\langle h_1,
h_2, \ldots, h_n\rangle$ be a real-$\CE$-process on $I$. Then there
exist an open, simply-connected $\Omega\subseteq\MC$ and an
$\CE$-process $\boldsymbol{f}$ on $\Omega$ such that:
\begin{enumerate}
\item $I\subset\Omega$
\item For all $t\in I: \boldsymbol{f}(t)=\boldsymbol{h}(t)$.
\end{enumerate}
\end{lemma}
\begin{proof}
Let $\boldsymbol{P}=\langle
P_1,P_2,\ldots,P_n\rangle=\boldsymbol{P}_\CE$. For $\iton$, $P_i$ is
    a polynomial and therefore analytic on $\MC^n$. By Cauchy's
    existence theorem for ordinary differential equations with analytic
    right-hand sides~\cite{edm87}, for all $a\in I$, there exist a
non-empty open disk $D_a\subseteq\MC$ centered at $a$ and functions
$f_{a,1}, f_{a,2}, \ldots, f_{a,n}$ analytic on $D_a$ such that for
$\iton:$

\begin{enumerate}
\item $f_{a,i}(a)=h_i(a)$
\item $f_{a,i}'$ exists on $D_a$ and for all $t\in D_a:
f_{a,i}'(t)=P_i(f_{a,1}(t),f_{a,2}(t), \ldots, f_{a,n}(t))$. That
is, $\boldsymbol{f}_a=\langle f_{a,1},f_{a,2}, \ldots,
f_{a,n}\rangle$ is an $\CE$-process on $D_a$.
\end{enumerate}

Claim: For all $a\in I$, there exists $\delta_a\in\MR_{>0}$ such
that for all $t\in I\cap(a-\delta_a,
a+\delta_a):\bs{f}_a(t)=\bs{h}(t)$. To see this, by
Lemma~\ref{lem:uni_real_local}, for all $a\in I$ there exists
$\beta_a\in\MR_{>0}$ such that for all $t\in
(a-\beta_a,a+\beta_a)\cap D_a$, $\bs{f}_a(t)\in\MR^n$. Let
$I_a=(a-\beta_a,a+\beta_a)\cap D_a$. Note that $\bs{f}_a|_{I_a}$ is
a real-$\CE$-process on $I_a$ . By the theorem of uniqueness of
solutions to differential equations with ${\cal{C}}^1$ right-hand
sides \cite{hirsch04textbook}, there exists $\gamma_a\in\MR_{>0}$ such
that for all $t\in(a-\gamma_a, a+\gamma_a)\cap I_a\cap I$,
$\bs{f}_a(t)=\bs{h}(t)$. Clearly, we can choose
$\delta_a\in\MR_{>0}$ such that $(a-\delta_a,a+\delta_a)\subseteq$
$(a-\gamma_a, a+\gamma_a)\cap I_a$. This establishes the claim.

For all $a\in I$, let $\delta_a\in\MR_{>0}$ be such that for all
$t\in I\cap(a-\delta_a, a+\delta_a):\bs{f}_a(t) = \bs{h}(t)$. Let
$\widehat{D}_a$ be an open disk centered at $a$ of radius $\delta_a$.

Claim: For all $a_1, a_2\in I$, for all $t\in \widehat{D}_{a_1}\cap
\widehat{D}_{a_2}: \bs{f}_{a_1}(t) = \bs{f}_{a_2}(t)$. To see this,
suppose $\widehat{D}_{a_1}\cap \widehat{D}_{a_2}\neq\varnothing$. Let $J =
\widehat{D}_{a_1}\cap \widehat{D}_{a_2} \cap \MR$. Since $\widehat{D}_{a_1}$ and
$\widehat{D}_{a_2}$ are open disks centered on the real line, $J$ is a
non-empty open real interval. For all $t\in J$, by the claim above,
$\bs{f}_{a_1}(t)=\bs{h}(t)$ and $\bs{f}_{a_2}(t)=\bs{h}(t)$. Hence,
$\bs{f}_{a_1}(t)=\bs{f}_{a_2}(t)$. Since $J$ is a non-empty
interval, $J$ contains an accumulation point. Since $\bs{f}_{a_1}$
and $\bs{f}_{a_2}$ are analytic on $\widehat{D}_{a_1} \cap
\widehat{D}_{a_2}$ and $\widehat{D}_{a_1} \cap \widehat{D}_{a_2}$ is simply
connected, for all $t \in \widehat{D}_{a_1}\cap
\widehat{D}_{a_2}:\bs{f}_{a_1}(t)=\bs{f}_{a_2}(t)$. This establishes the
claim.

Let $\Omega=\bigcup_{a\in I}\widehat{D}_a$. Clearly, $I\subset\Omega$.
$\Omega$ is a union of open discs, and is therefore open.

For all $t\in\Omega$, there exists $a\in I$ such that
$t\in\widehat{D}_a$. Since $\widehat{D}_a$ is a disk, $t$ and $a$ are
path-connected in $\Omega$. Since $I$ is path-connected, and
$I\subseteq\Omega$, it follows that $\Omega$ is path-connected.

To see that $\Omega$ is simply-connected, consider the function
$R:[0,1]\times\Omega\rightarrow\Omega$ given by $(u,z)\mapsto
\op{Re}(z)+i \op{Im}(z)(1-u)$. Observe that $R$ is continuous on
$[0,1]\times\Omega$, and for all $z\in\Omega$:
$R(0,z)=z$, $R(1,\Omega)\subset\Omega$, and for all $u\in[0,1]$, for
all $z\in\Omega\cap\MR: R(u,z)\in\Omega$. Therefore, $R$ is a
deformation retraction. Note that $R(0,\Omega)=\Omega$ and
$R(1,\Omega)\subseteq\MR$, and $\Omega$ is path-connected together
imply that $R(1,\Omega)$ is a real interval. Hence, $R(1,\Omega)$ is
simply-connected. Since $R$ was a deformation retraction, $\Omega$
is simply-connected.

Let $\bs{f}:\Omega\rightarrow\MC^n$ be the unique function such that
for all $a\in I$, for all $t\in \widehat{D}_a:\bs{f}(t)=\bs{f}_a(t)$. By
the claim above and from the definition of $\Omega$, $\bs{f}$ is
well-defined.

Observe that for all $t\in I$,
\begin{align*}
\bs{h}(t)&=\bs{f}_t(t) && \text{(Definition of $\bs{f}_t$)} \\
&=\bs{f}(t)  && \text{($I\subset\Omega$ and definition of $\bs{f}$).}
\end{align*}

Claim: $\bs{f}$ is an $\CE$-process on $\Omega$. From the
definitions of $\Omega$ and $\bs{f}$, for all $t\in\Omega$, there
exists $a\in I$ such that $t\in \widehat{D}_a$ and for all
$s\in\widehat{D}_a$, $\bs{f}(s)=\bs{f}_a(s)$. Since $\bs{f}_a$ is an
$\CE$-process on $\widehat{D}_a$, the claim follows.
\end{proof}

In Theorem~\ref{thm:stays_positive}, we prove that if $\CE$ is a finite, physical event-system, then $\CE$-processes that begin at positive (respectively non-negative) points remain positive (respectively non-negative) through all forward real time where they are defined. In fact, Theorem~\ref{thm:stays_positive} establishes more detail about $\CE$-processes. In particular, if at some time a species' concentration is positive, then it will be positive at subsequent times.

\begin{theorem}
\label{thm:stays_positive}
Let $\CE$ be a finite, physical event-system of dimension $n$, let
$\Omega\subseteq\MC$ be open and simply-connected, and let $\bs{f}=\langle f_1, f_2,\ldots, f_n\rangle$ be an $\CE$-process on $\Omega$.
If $I\subseteq\Omega\cap\MR_{\geq 0}$ is connected and $0\in I$ and
${\bs f}(0)$ is a non-negative point then for $k=1,2,\ldots,n$ either:
\begin{enumerate}
\item For all $t\in I$, $f_k(t)=0$, or
\item For all $t\in I\cap\MR_{>0}$, $f_k(t)\in\MR_{>0}$.
\end{enumerate}
\end{theorem}

The proof of Theorem~\ref{thm:stays_positive} is highly technical, and relies on a detailed examination of the vector of polynomials $\bs{P}_\CE$. This allows us to show~(Lemma~\ref{lem:RNN}) that if $\bs{f}=\langle f_1,f_2,\dots,f_n\rangle$ is an $\CE$-process that at real time $t_0$ is non-negative, then each $f_i$ is ``right non-negative.'' That is, the Taylor series expansion of $f_i$ around $t_0$ has real coefficients and the first non-zero coefficient, if any, is positive. Further, (Lemma~\ref{lem:bdy_O1}) if $f_i(t_0)=0$ and its Taylor series expansion has a non-zero coefficient, then there exists $k$ such that $f_k(t_0)=0$ and the first derivative of $f_k$ with respect to time is positive at $t_0$.

\begin{definition} Let $n\in\MZ_{>0}$ and let $k\in\{1,2,\ldots,n\}$. A polynomial $f\in\MR[X_1,X_2,\ldots,X_n]$ is {\em non-nullifying with
respect to} $k$ iff there exist $m\in\MN$, $c_1, c_2, \cdots, c_m  \in \MR_{>0}, M_1, M_2, \ldots, M_m \in
\mathbb{M}_{\{X_1,X_2,\ldots,X_n\}}$ and $h\in\MR[X_1,X_2,\ldots,X_n]$
such that $f = \sum_{i=1}^m c_i
M_i + X_k h$.
\end{definition}

Observe that for all $k$, the polynomial $0$ is non-nullifying with respect to $k$.

\begin{lemma}\label{lem:esnonnull}
Let $\CE$ be a finite, physical event-system of dimension $n$. Let
$\langle P_1, P_2, \ldots, P_n\rangle=\bs{P}_{\CE}$. Then, for all
$i\in\{1,2,\ldots,n\}$, $P_i$ is non-nullifying with respect to $i$.
\end{lemma}
\begin{proof}
Let $m=|\CE|$. Let $(\gamma_{j,i})_{m\times n}=\Gamma_\CE.$ Since $\CE$ is physical, there exist $\sigma_1, \sigma_2,\ldots,\sigma_m$, $\tau_1, \tau_2, \ldots,\tau_m\in\MR_{>0}$ and $M_1,M_2,\ldots,M_m, N_1,N_2,\ldots,N_m\in\mathbb{M}_\infty$ such that for $j=1,2,\ldots,m: M_j\prec N_j$ and $\{\sigma_1M_1-\tau_1N_1, \sigma_2M_2,\tau_2N_2,\ldots,\sigma_mM_m-\tau_mN_m\}=\CE$. Let $i\in\{1,2,\ldots,n\}$.

From the definition of $\bs{P}_\CE$, $P_i=\sum_{j=1}^m\gamma_{j,i} (\sigma_jM_j-\tau_jN_j)$. It is sufficient to prove that for $j=1,2,\ldots,m: \gamma_{j,i} (\sigma_jM_j-\tau_jN_j)$ is non-nullifying with respect to $i$. Let $j\in\{1,2,\ldots,m\}$. If $\gamma_{j,i}=0$ then $\gamma_{j,i} (\sigma_jM_j-\tau_jN_j)=0$ which is non-nullifying with respect to $i$. If $\gamma_{j,i}>0$ then, from the definition of $\Gamma_{\CE}$, $X_i \mid N_j$ and
\[\gamma_{j,i} (\sigma_jM_j-\tau_jN_j)=\gamma_{j,i} \sigma_jM_j+X_i\left(-\gamma_{j,i}\tau_j\frac{N_j}{X_i}\right)\] which is non-nullifying with respect to $i$ since $\gamma_{j,i} \sigma_j>0$. Similarly, if $\gamma_{j,i}<0$ then $X_i\mid M_j$ and \[\gamma_{j,i} (\sigma_jM_j-\tau_jN_j)=-\gamma_{j,i} \tau_jN_j+X_i\gamma_{j,i}\sigma_j\frac{M_j}{X_i}\] which is non-nullifying with respect to $i$ since $-\gamma_{j,i} \tau_j>0$. Hence, $P_i$ is non-nullifying with respect to $i$.
\end{proof}

\begin{definition}
Let $t_0\in\MC$, let $f:{\MC}\rightarrow{\MC}$ be
analytic at $t_0$ and let $f(t)=\sum_{k=0}^\infty c_k (t-t_0)^k$
be the Taylor series expansion of $f$ around
$t_0$. Then $O(f,t_0)$ is the least $k$ such that $c_k\not=0$. If
for all $k$,
$c_k=0$, then $O(f,t_0)=\infty$.
\end{definition}

\begin{definition}[Right non-negative]
Let $t_0\in\MR$, let $f:{\MC}\rightarrow{\MC}$ be
analytic at $t_0$ and let $f(t)=\sum_{k=0}^\infty c_k (t-t_0)^k$
be the Taylor series expansion of $f$ around
$t_0$. Then $f$ is {\em RNN at $t_0$} iff both:
\begin{enumerate}
\item For all $k\in \mathbb{N}$, $c_k\in\MR$ and
\item Either $O(f,t_0)=\infty$ or $c_{O(f,t_0)}\in\MR_{>0}$.
\end{enumerate}
\end{definition}

\begin{lemma}
\label{lem:cdot} Let $t_0\in\MC$. Let $f,g:{\MC}\rightarrow{\MC}$ be functions analytic at $t_0$. Then:
\begin{enumerate}
\item\label{lem:cdot:orders}$O(f\cdot g,t_0)=O(f,t_0)+O(g,t_0)$.
\item\label{lem:cdot:rnn}If $t_0\in\MR$ and $f,g$ are RNN at $t_0$ then $f\cdot g$ is RNN at $t_0$.
\end{enumerate}
\end{lemma}
The proof is obvious.

\begin{lemma}\label{lem:RNN}
Let $\CE$ be a finite, physical event-system of dimension $n$, let $\Omega\subseteq\MC$ be open and simply-connected and let $\bs{f}=\langle f_1,f_2,\ldots,f_n\rangle$ be
an $\CE$-process on $\Omega$. For all $t_0\in\Omega\cap\MR$, if
$\bs{f}(t_0)\in\MR^n_{\geq 0}$ then for $\iton: f_i$ is RNN at $t_0$.
\end{lemma}
\begin{proof}
Suppose $t_0\in\Omega\cap\MR$ and $\bs{f}(t_0)\in\MR^n_{\geq 0}$. Let $\bs{P}=\langle P_1,P_2,\ldots,P_n\rangle=\bs{P}_\CE$. Let $C=\{i\mid f_i$ is not RNN at $t_0\}$.

For the sake of contradiction, suppose $C\neq\varnothing$. Let $m=\min_{i\in C} O(f_i, t_0)$. Let $k\in C$ be such that $O(f_{k}, t_0)=m$. Let $f_k(t)=\sum_{i=0}^\infty a_i (t-t_0)^i$
be the Taylor series expansion of $f_k$ around
$t_0$. Since $\CE$ is physical and $t_0\in\MR$ and $\bs{f}(t_0)\in\MR^n_{\geq 0}$, it follows from Lemma~\ref{lem:uni_real_local}.\ref{taylor} that for all $i\in\mathbb{N}$, $a_i\in\MR$. Further:
\begin{align}\label{lem:RNN:A1}
a_0=a_1=\ldots=a_{m-1}=0&&(O(f_{k}, t_0)=m.)
\\\label{lem:RNN:A2}a_m\in\MR_{<0}&&(f_k\text{ is not RNN at }t_0.)
\end{align}
Since $\bs{f}(t_0)\in\MR^n_{\geq 0}$ and $a_m\in\MR_{<0}$ and $a_0=f_k(t_0)$, it follows that $m>0$.

Consider $f_k'=P_k\circ\bs{f}$. By differentiation, the Taylor series expansion of $f_k'$ at $t_0$ is:
\begin{align}\label{eqn:RNN:1}
f_k'(t)=\sum_{i=0}^\infty (i+1)a_{i+1} (t-t_0)^i.
\end{align}
From Lemma~\ref{lem:esnonnull}, $P_k$ is non-nullifying. Hence, there exist $l\in\MN$, $b_1,
b_2, \ldots, b_l\in\MR_{>0}$, $M_1, M_2, \ldots,
M_l\in{\mathbb M}_{\{X_1, X_2, \ldots, X_n\}}$ and $h\in\MR[X_1,X_2,\ldots,X_n]$ such that $P_k=\sum_{j=1}^l b_j M_j + X_k\cdot h$. Then for all $t\in\Omega$:
\begin{align}\label{lem:RNN:B}
f_k'(t)&=P_k\circ\bs{f}(t)=\sum_{j=1}^l b_j M_j\circ\bs{f}(t) + f_k(t)\cdot (h\circ\bs{f}(t))
\end{align}

Since $h$ is a polynomial, $h\circ\bs{f}$ is analytic at $t_0$. Therefore, $f_k\cdot(h\circ\bs{f})$ is analytic at $t_0$. Let $\sum_{i=0}^\infty c_i(t-t_0)^i$ be the Taylor series expansion of $f_k\cdot(h\circ\bs{f})$ at $t_0$. Similarly, for $j=1,2,\ldots,l$, $b_jM_j\circ\bs{f}$ is analytic at $t_0$. Let $ \sum_{i=0}^\infty d_{j,i}(t-t_0)^i$ be the Taylor series expansion of $b_jM_j\circ\bs{f}$ at $t_0$.
From \eqref{eqn:RNN:1},\eqref{lem:RNN:B}, equating Taylor series coefficients, for $i=0,1,\ldots,m-1$:
\begin{align}\label{eqn:RNN:F}
(i+1)a_{i+1}=c_i+\sum_{j=1}^l d_{j,i}
\end{align}

From Lemma~\ref{lem:cdot}.\ref{lem:cdot:orders},
\begin{align*}
O(f_k\cdot(h\circ\bs{f}),t_0)=O(f_k,t_0)+O(h\circ\bs{f},t_0)\geq O(f_k,t_0)=m
\end{align*}
Hence,
\begin{align}\label{lem:RNN:C}
c_0=c_1=\ldots=c_{m-1}=0.
\end{align}

From \eqref{lem:RNN:A1}, \eqref{eqn:RNN:F}, \eqref{lem:RNN:C}, for $i=0,1,\ldots,m-2$:
\begin{align}\label{lem:RNN:D}
\sum_{j=1}^l d_{j,i}=0
\end{align}
Since $m>0$, from \eqref{lem:RNN:A2}, \eqref{eqn:RNN:F}, \eqref{lem:RNN:C}:
\begin{align}\label{lem:RNN:E}
\sum_{j=1}^l d_{j,m-1}=ma_m\in\MR_{<0}
\end{align}

Let $i_0=\min_{j=1,2,\ldots,l}\{O(b_jM_j\circ\bs{f},t_0)\}$. From \eqref{lem:RNN:E}, it follows that $i_0\leq m-1$.\
\\
\\Case 1: For $j=1,2,\ldots,l: d_{j,i_0}\in\MR_{\geq 0}$. From the definition of $i_0$ it follows that $\sum_{j=1}^l d_{j,i_0} \in\MR_{>0}$. If $i_0<m-1$, this contradicts \eqref{lem:RNN:D}. If $i_0=m-1$, this contradicts \eqref{lem:RNN:E}.\
\\
\\Case 2: There exists $j_0\in\{1,2,\ldots,l\}$ such that $d_{j_0,i_0}\in\MR_{<0}$. From the definition of $i_0$, $O(b_{j_0}M_{j_0},t_0)=i_0\leq m-1$. Therefore, for each $i$ such that $X_i \mid M_{j_0}$, $O(f_i,t_0)\leq m-1$. From the definitions of $C$ and $m$, this implies that for each $i$ such that $X_i \mid M_{j_0}$, $f_i$ is RNN at $t_0$. Since $b_{j_0}\in\MR_{>0}$, it follows that $b_{j_0}M_{j_0}\circ\bs{f}$ is a product of RNN functions. Hence, by Lemma~\ref{lem:cdot}.\ref{lem:cdot:rnn}, $b_{j_0}M_{j_0}\circ\bs{f}$ is RNN at $t_0$ and $d_{j_0,i_0}\in\MR_{>0}$, a contradiction.

Hence, for $\iton$, $f_i$ is RNN at $t_0$.
\end{proof}

\begin{lemma}
\label{lem:rnn_sign}
Let $t_0\in{\MR}$ and let $f$ be a function RNN at $t_0$. There exists an
$\varepsilon\in{\MP}$ such that either for all $t\in(t_0,t_0+\varepsilon)$, $f(t)\in\MR_{>0}$
or for all $t\in(t_0,t_0+\varepsilon)$, $f(t)=0$.
\end{lemma}
\begin{proof}
Let $m=O(f,t_0)$. If $m=\infty$, $f$ is identically zero and the lemma follows immediately. Otherwise, let $f^{(m)}$ denote the $m^\text{th}$ derivative of $f$. Since $f$ is RNN at $t_0$ and has order $m$, $f^{(m)}(t_0)\in\MR_{>0}$. Since $f$ is analytic at $t_0$, $f^{(m)}$ is analytic at $t_0$, and hence continuous at $t_0$. By continuity, there exists $\varepsilon\in\MR_{>0}$ such that for all $\tau\in[t_0,t_0+\varepsilon]: f^{(m)}(\tau)\in\MR_{>0}$. From Taylor's theorem, for all $t\in (t_0,t_0+\varepsilon)$, there exists $\tau\in [t_0,t_0+\varepsilon]$ such that:
\begin{align*}
f(t)=\frac{(t-t_0)^m}{m!}f^{(m)}(\tau)
\end{align*}
Therefore, $f(t)\in\MR_{>0}$.
\end{proof}

Note that Lemma~\ref{lem:RNN} and Lemma~\ref{lem:rnn_sign} together already imply that if $\CE$ is a finite, physical event-system, then $\CE$-processes that begin at non-negative points remain non-negative through all forward real time where they are defined. This result is weaker than Theorem~\ref{thm:stays_positive}.

\begin{lemma}
\label{lem:bdy_O1}Let $\CE$ be a finite, physical event-system of
dimension $n$, let $\Omega \subseteq \MC$ be open and
simply-connected, let $\bs{f}=\langle f_1,f_2,\ldots,f_n\rangle$ be
an $\CE$-process on $\Omega$. Let $t_0\in\Omega$. If
$\bs{f}(t_0)$ is non-negative and there exists
$j\in\{1,2,\ldots,n\}$ such that $0 < O(f_j, t_0) < \infty$ then
there exists $k\in\{1,2,\ldots,n\}$ such that $O(f_{k},t_0)=1$.
\end{lemma}
\begin{proof}
Suppose $\bs{f}(t_0)\in\MR^n_{\geq 0}$. Let $C=\{i\mid 0<O(f_i,t_0)<\infty\}$. Suppose $C\neq\varnothing$. Let $m=\min_{i\in C} O(f_i, t_0)$. There exists $k\in C$ such that $O(f_{k}, t_0)=m$.

Let $\bs{P}=\langle P_1,P_2,\ldots,P_n\rangle=\bs{P}_\CE$. From Lemma~\ref{lem:esnonnull}, $P_k$ is non-nullifying with respect to $k$. Hence, there exist $l\in\MN$, $b_1,
b_2, \ldots, b_l\in\MR_{>0}$, $M_1, M_2, \ldots,
M_l\in\mathbb{M}_{\{X_1, X_2, \ldots, X_n\}}$ and $h\in\MR[X_1,X_2,\ldots,X_n]$ such that $P_k=\sum_{j=1}^l b_j M_j + X_k\cdot h$.

For all $t\in\Omega$: $f_k'(t)=P_k\circ\bs{f}(t)=\sum_{j=1}^l b_j M_j\circ\bs{f}(t) + f_k(t)\cdot (h\circ\bs{f}(t))$. From Lemma~\ref{lem:cdot}.\ref{lem:cdot:orders}, $O(f_k\cdot(h\circ\bs{f}),t_0)=O(f_k,t_0)+O(h\circ\bs{f},t_0)\geq O(f_k,t_0)=m$. It follows that:
\begin{align}\label{eqn:ord:C}
m-1=O(f_k',t_0)=O(\sum_{j=1}^l b_j M_j\circ\bs{f},t_0)
\end{align}

From Lemma~(\ref{lem:cdot}.\ref{lem:cdot:rnn}) and Lemma~\eqref{lem:RNN}, for $j=1,2,\ldots,l:$ $b_jM_j\circ\bs{f}$ is RNN at $t_0$. It follows that $O(\sum_{j=1}^l b_j M_j\circ\bs{f},t_0)=\min_{j=1,2,\ldots,l}O(b_jM_j\circ\bs{f},t_0)$. From Equation~\eqref{eqn:ord:C}, $m-1=\min_{j=1,2,\ldots,l}O(b_jM_j\circ\bs{f},t_0)$.
Hence, there exists $j_0$ such that $O(b_{j_0}M_{j_0}\circ\bs{f},t_0)=m-1$. From Lemma~(\ref{lem:cdot}.\ref{lem:cdot:orders}), for all $i$ such that $X_i\mid M_{j_0}$, $O(f_i,t_0)\leq m-1$. From the definition of $m$, for all $i$ such that $X_i\mid M_{j_0}$, $O(f_i,t_0)=0$. It follows that $m-1=O(b_{j_0}M_{j_0}\circ\bs{f},t_0)=0$. Hence, $m=1$.
\end{proof}

We are now ready to prove Theorem~\ref{thm:stays_positive}.

\begin{proof}[Proof of Theorem~\ref{thm:stays_positive}]
Suppose $I\subseteq\Omega\cap\MR_{\geq 0}$ is connected and $0\in I$ and
${\bs f}(0)$ is a non-negative point. If $I\cap\MR_{>0}=\varnothing$, the theorem is immediate. Suppose $I\cap\MR_{>0}\neq\varnothing$.

It is clear that for all $k$, $O(f_k,0)=\infty$ iff for all $t\in I$, $f_k(t)=0$. Let $C=\{i\mid O(f_i,0)\neq \infty\}$. From Lemma~\eqref{lem:RNN} and Lemma~\eqref{lem:rnn_sign}, for all $k\in C$, there exists $\varepsilon_k\in I\cap\MR_{>0}$ such that for all $t\in (0,\varepsilon_k): f_k(t)\in\MR_{>0}$.

Suppose for the sake of contradiction that there exist $i\in C$ and $t\in I\cap\MR_{>0}$ such that $f_i(t)\notin\MR_{>0}$. From Lemma~\eqref{lem:uni_real_global}, $f_i(t)\in\MR$. Since $f_i(\varepsilon_i/2)\in\MR_{>0}$ and $f_i(t)\in\MR_{\leq 0}$, by continuity there exists $t'\in I\cap\MR_{>0}$ such that $f_i(t')=0$.

Let $t_0=\inf\{t\in I\cap\MR_{>0}\mid$ There exists $i\in C$ with $f_i(t)=0\}$.
It follows that:
\begin{enumerate}
\item $t_0\in\MR_{>0}$ because $t_0\geq\min_{i\in C}\{\varepsilon_i\}$.
\item $\bs{f}(t_0)\in\MR^n_{\geq 0}$, from the definition of $t_0$.
\item There exists $i_1\in C$ such that $O(f_{i_1},t_0)=1$. This follows because there exist $i_0\in C$ and $T\subseteq I\cap\MR_{>0}$ such that $t_0=\inf(T)$ and for all $t\in T$: $f_{i_0}(t)=0$. By continuity, $f_{i_0}(t_0)=0$. Hence, $O(f_{i_0},t_0)>0$. Since $i_0\in C$, $O(f_{i_0},0)\neq \infty$. By connectedness of $I$, $O(f_{i_0},t_0)\neq\infty$. Therefore, $0<O(f_{i_0},t_0)<\infty$. Since $\bs{f}(t_0)\in\MR^n_{\geq 0}$, by Lemma~\eqref{lem:bdy_O1}, there exists $i_1\in\{1,2,\ldots,n\}$ such that $O(f_{i_1},t_0)=1$. Assume $i_1\notin C$. Then $O(f_{i_1},0)=\infty$. By connectedness of $I$, $O(f_{i_1},t_0)=\infty$, contradicting that $O(f_{i_1},t_0)=1$. Hence, $i_1\in C$.
\end{enumerate}

Hence, $f_{i_1}(t_0)=0$. Since $\bs{f}(t_0)\in\MR^n_{\geq 0}$, by Lemma~\eqref{lem:RNN} $f_{i_1}'(t_0)\in\MR_{>0}$.

From the definition of $t_0$, for all $t\in (0,t_0)$, $f_{i_1}(t)\in\MR_{>0}$. Since $t_0\in\MR_{>0}$,
\[
f_{i_1}'(t_0)=\lim_{h\to 0^+}\frac{f_{i_1}(t_0)-f_{i_1}(t_0-h)}{h}=\lim_{h\to 0^+}\frac{-f_{i_1}(t_0-h)}{h}\in\MR_{\leq0},
\]
a contradiction. The theorem follows.
\end{proof}

There is a notion in chemistry that, for systems of chemical reactions, concentrations evolve through time to reach equilibrium. In later sections of this paper, we will investigate this notion. In the remainder of this section of the paper, we will prepare for that investigation.

\begin{definition}\label{omega limit}
Let $\CE$ be a finite event-system of dimension $n$, let
$\Omega \subseteq \MC$ be open, simply connected and such that
$\MR_{\geq 0} \subseteq \Omega$, let $\bs{f}$ be an $\CE$-process on
$\Omega$, and let $\bs{q} \in \MC^n$. Then $\bs{q}$ is an
\emph{$\omega$-limit point of $\bs{f}$} iff for all
$\varepsilon\in\MR_{>0}$ there exists a sequence of non-negative
reals $\{t_i\}_{i\in\MZ_{>0}}$ such that $t_i\rightarrow\infty$ as
$i\rightarrow\infty$ and for all $i\in\MZ_{>0}$, $\|\bs{f}(t_i) -
\bs{q}\|_2 < \varepsilon$.
\end{definition}

Sometimes, an $\omega$-limit is defined by the existence of a single sequence of times such that the value approaches the limit. The above definition is easily seen to be equivalent.

\begin{definition}\label{invariant set}
Let $\CE$ be a finite event-system of dimension $n$ and let
$S \subseteq \MC^n$. $S$ is an \emph{invariant set of $\CE$} iff for
all $\bs{q} \in S$, for all open, simply-connected $\Omega \subseteq
\MC$, for all $\CE$-processes $\bs{f}$ on $\Omega$, if $0 \in
\Omega$ and $\bs{f}(0) = \bs{q}$ then for all $t \in \MR_{\geq 0}$
such that $[0,t] \subseteq \Omega$, $\bs{f}(t) \in S$.
\end{definition}

\begin{lemma}\label{invariant}
Let $\CE$ be a finite, physical event-system of dimension $n$, let
$\Omega \subseteq \MC$ be open and simply connected, and let
$\bs{f}$ be an $\CE$-process on $\Omega$. If $\MR_{\geq 0} \subseteq
\Omega$ and $\bs{f}(0)$ is a non-negative point, then the set of all
$\omega$-limit points of $\bs{f}$ is an invariant set of $\CE$ and
is contained in $\MR^n_{\geq 0}$.
\end{lemma}
\begin{proof}
Let $S$ be the set of all $\omega$-limit points of $\bs{f}$. By
Lemma~\ref{thm:stays_positive}, for all $t\in \MR_{\geq 0}$, $\bs{f}(t) \in
\MR^n_{\geq 0}$, hence $S \subseteq \MR^n_{\geq 0}$.

Let $\bs{q} \in S$, let $\widetilde{\Omega} \subseteq \MC$ be open,
simply-connected, and such that $0 \in \widetilde{\Omega}$, and let
$\bs{h}$ be an $\CE$-process on $\widetilde{\Omega}$ such that
$\bs{h}(0) = \bs{q}$. Suppose $u \in \MR_{\geq 0}$ and $[0,u]
\subseteq \widetilde{\Omega}$. Since $\CE$ is finite and
physical, $\bs{P}_{\CE}|_{\MR^n}$ can be viewed as a map
$\bs{F}:\MR^n\rightarrow\MR^n$ of class ${\cal{C}}^1$. By
Lemma~\ref{lem:uni_real_global}, for all $t\in [0,u]$, $\bs{h}(t) \in \MR^n$, so
$\bs{h}|_{[0,u]}$ can be viewed as a map $\bs{X}:[0,u] \rightarrow
\MR^n$ such that $\bs{X}'=\bs{F}(\bs{X})$. By~\cite[p.\hspace{3pt}147]{hirsch04textbook}, there exists a neighborhood $U\subset\MR^n$ of $\bs{q}$ and a
constant $K$ such that for all $\bs{\alpha}\in U$, there exists a
unique real-$\CE$-process $\bs{\rho_\alpha}$ defined on $[0,u]$ with
$\bs{\rho_\alpha}(0) = \bs{\alpha}$ and
$\|\bs{\rho_\alpha}(u)-\bs{h}(u)\|_2 \leq K \|\bs{\alpha}-\bs{q}\|_2
\exp(K u)$. Observe that necessarily $K \in \MR_{\geq 0}$. By
Lemma~\ref{lem:solution} for all $\bs{\alpha} \in U$ there exists an
open, simply-connected $\Omega_{\bs{\alpha}} \subseteq \MC$ and an
$\CE$-process $\bs{\varrho_\alpha}$ on $\Omega_{\bs{\alpha}}$ such
that $[0,u] \subseteq \Omega_{\bs{\alpha}}$ and for all $t \in
[0,u]$, $\bs{\varrho_\alpha}(t)=\bs{\rho_\alpha}(t)$. Therefore,
$\|\bs{\varrho_\alpha}(u)-\bs{h}(u)\|_2 \leq K \|\bs{\alpha} -
\bs{q}\|_2 \exp(K u)$.

Let $\varepsilon \in \MR_{>0}$ and let $\delta_1, \delta_2 \in
\MR_{>0}$ be such that $K \delta_1 \exp(K u) \leq \varepsilon$ and
the open ball centered at $\bs{q}$ of radius $\delta_2$ is contained
in $U$. Let $\delta = \min(\delta_1, \delta_2)$. Since $\bs{q}$ is
an $\omega$-limit point of $\bs{f}$, there exists a sequence of
non-negative reals $\{t_i\}_{i\in \MZ_{>0}}$ such that $t_i
\rightarrow \infty$ as $i \rightarrow \infty$ and for all
$i\in\MZ_{>0}$, $\|\bs{f}(t_i)-\bs{q}\|_2 < \delta$. Then for all $i
\in \MZ_{>0}$, $\bs{f}(t_i) \in U$, so by Lemma~\ref{time
invariance} for all $t \in [0,u]$,
$\bs{f}(t_i+t)=\bs{\varrho}_{\bs{f}(t_i)}(t)$. Then
\begin{align*}
\|\bs{f}(t_i+u)-\bs{h}(u)\|_2 &= \|\bs{\varrho}_{\bs{f}(t_i)}(u)-\bs{h}(u)\|_2 \\
&\leq K \|\bs{f}(t_i)-\bs{q}\|_2 \exp(K u) \\
&\leq K \delta \exp(K u) \\
&\leq \varepsilon
\end{align*}

Thus $\bs{h}(u)$ is an $\omega$-limit point of $\bs{f}$, so $S$ is
an invariant set of $\CE$.
\end{proof} 
\section{Finite Natural Event-systems}\label{sec:natural}
In this section, we focus on finite, natural event-systems --- a subclass of finite, physical event-systems which has much in common with systems of chemical reactions that obey detailed balance.

In chemical reactions, the total bond energy of the
reactants minus the total bond energy of the products is a measure
of the heat released.
For example, in the reaction,
$\sigma X_2 - \tau X_1$, $\ln
\left(\frac{\sigma}{\tau}\right)$ is taken to be the quantity of heat
released. If there are
multiple reaction paths that take the same reactants to the same
products, then the quantity of heat released along each path
must be the same.

The finite, physical event-system $\CE=\{2 X_2 -X_1, X_2-X_1\}$ does not
behave like a chemical reaction system since, when $X_2$ is converted to
$X_1$ by the first reaction, $\ln \left(2\right)$ units of heat are
released; however, when $X_2$ is converted to $X_1$ by the second
reaction, $\ln \left(1\right)=0$ units of heat are released. When an
event-system admits a pair of paths from the same reactants to the
same products but with different quantities of heat released, we say
that the system has an ``energy cycle.''

\begin{definition}[Energy cycle$\label{EnergyCycle}$] Let $\CE$
be a finite, physical event-system. $\CE$ has an energy cycle iff
$G_\CE$ has a cycle of non-zero weight.
\end{definition}

\begin{example}
For the physical event-system $\CE_1=\{2 X_2 -X_1, X_2-X_1\}$, the
event $X_2-X_1$ induces an edge $\langle X_2, X_1\rangle$ in the
event graph with weight $\ln\left(\frac{1}{1}\right)=0$. The event
$2X_2-X_1$ induces an edge $\langle X_1, X_2\rangle$ with weight
$-\ln\left(\frac{2}{1}\right)=-\ln\left(2\right)$. The weight of the
cycle from $X_2$ to $X_1$ and back to $X_2$ using these two edges,
is $-\ln\left(2\right)\neq 0$. Hence, $\CE_1$ has an energy cycle by
Definition~\ref{EnergyCycle}.
\end{example}

\begin{example} For the physical event-system
$\CE_2=\{X_2-X_1, 2X_3X_4 - X_2X_3, X_4X_5-X_1X_5\}$, the cycle
$\langle X_3X_4X_5,X_2X_3X_5,X_1X_3X_5,X_3X_4X_5\rangle$ is induced
by the sequence of events $2X_3X_4 - X_2X_3, X_2-X_1, X_4X_5-X_1X_5$
and has corresponding weight
$\ln\frac{2}{1}+\ln\frac{1}{1}+\ln\frac{1}{1}=\ln\left(2\right)\neq0$.
Hence, $\CE_2$ has an energy cycle.
\end{example}

The following theorem gives multiple characterizations of natural event-systems.

\begin{theorem}
$\label{NatL}$ Let $\CE$ be a finite, physical event-system of
dimension $n$. The following are equivalent:
\begin{enumerate}
\item $\CE$ is natural.
\item $\CE$ has a strong equilibrium point that is not a z-point.
$($i.e. there exists $\bs{\alpha}\in\MC^n$ such that for all $i=1$ to
$n$, $\alpha_i\neq0$ and for all $e\in\CE$,
$e\left(\bs{\alpha}\right)=0.)$
\item $\CE$ has no energy cycles.
\item If $\CE=\{\sigma_1M_1-\tau_1N_1, \sigma_2M_2 - \tau_2N_2, \ldots,
\sigma_mM_m - \tau_mN_m\}$ and for all $j=1$ to $m$, $M_j\prec N_j$
and $\sigma_j, \tau_j>0$ then there exists $\bs{\alpha}\in\MR^n$ such
that $\Gamma_\CE \bs{\alpha} = \left\langle
\ln\left(\frac{\sigma_1}{\tau_1}\right), \ldots,
\ln\left(\frac{\sigma_m}{\tau_m}\right)\right\rangle^T$.
\end{enumerate}
\end{theorem}

To prove Theorem~\ref{NatL}, we will use the following lemma.

\begin{lemma}
\label{StrongEqCorr} Let $\CE=\{\sigma_1M_1-\tau_1N_1,
\sigma_2M_2-\tau_2N_2, \ldots, \sigma_mM_m-\tau_mN_m\}$ be a finite,
physical event-system of dimension $n$ such that for all $j=1$ to
$m$, $\sigma_j, \tau_j>0$ and $M_j\prec N_j$. Then for all
$\bs{\alpha}=\langle\alpha_1, \alpha_2, \ldots, \alpha_n\rangle^T\in\MR^n$,
$\Gamma_\CE\cdot\bs{\alpha}=\left\langle\ln\left(\frac{\sigma_1}{\tau_1}\right),
\ln\left(\frac{\sigma_2}{\tau_2}\right), \ldots,
\ln\left(\frac{\sigma_m}{\tau_m}\right)\right\rangle^T$ iff $\langle
\mathrm{e}^{\alpha_1}, \cdots, \mathrm{e}^{\alpha_n}\rangle$ is a
positive strong $\CE$-equilibrium point.
\end{lemma}
\begin{proof}
Let $\CE=\{\sigma_1M_1-\tau_1N_1, \sigma_2M_2 - \tau_2N_2, \ldots,
\sigma_mM_m - \tau_mN_m\}$ and for all $j=1$ to $m$, $M_j\prec N_j$
and $\sigma_j, \tau_j>0$. Let $\Gamma=\Gamma_\CE$. For all
$\bs{\alpha}=\langle\alpha_1, \ldots, \alpha_n\rangle\in\MR^n$,
\begin{align*}
&\Gamma\bs{\alpha}=\left\langle\ln\left(\frac{\sigma_1}{\tau_1}\right),
\ln\left(\frac{\sigma_2}{\tau_2}\right), \ldots,
\ln\left(\frac{\sigma_m}{\tau_m}\right)\right\rangle^T\\
\Leftrightarrow\,&\sum_{i=1}^n\gamma_{j,i}\alpha_i=\ln\left(\sigma_j/\tau_j\right),
\,\forall\jtom\\
\Leftrightarrow\,&\prod_{i=1}^n\left(\mathrm{e}^{\alpha_i}\right)^{\gamma_{j,i}}
=\sigma_j/\tau_j,\,\forall\jtom&\text{(Exponentiation.)}\\
\Leftrightarrow& N_j\left(\langle \mathrm{e}^{\alpha_1}, \ldots,
\mathrm{e}^{\alpha_n}\rangle\right)/M_j\left(\langle
\mathrm{e}^{\alpha_1}, \ldots,
\mathrm{e}^{\alpha_n}\rangle\right)=\sigma_j/\tau_j,\, \forall\jtom&\text{(Definition of $\Gamma$.)}\\
\Leftrightarrow\, &\sigma_jM_j\left(\langle
\mathrm{e}^{\alpha_1}, \ldots, \mathrm{e}^{\alpha_n}\rangle\right) -
\tau_jN_j\left(\langle \mathrm{e}^{\alpha_1}, \ldots,
\mathrm{e}^{\alpha_n}\rangle\right)=0,\, \forall\jtom\\
\Leftrightarrow\, &\langle \mathrm{e}^{\alpha_1}, \ldots,
\mathrm{e}^{\alpha_n}\rangle\text{ is a positive strong $\CE$-equilibrium
point.}
\end{align*}
\end{proof}

\begin{proof}[Proof of Theorem~\ref{NatL}]
$\boldsymbol{\left(4\right)\Rightarrow\left(1\right)}:$ Follows from
Lemma~\ref{StrongEqCorr}.
\\$\boldsymbol{\left(1\right)\Rightarrow\left(2\right)}:$ Follows immediately from
definitions.\
\\$\boldsymbol{\left(2\right)\Rightarrow\left(3\right)}:$
\\Consider an arbitrary cycle $\mathcal{C}$ in $G_\CE$ given by the sequence of
$k$ edges \\$\{\langle v_0, v_1\rangle,\langle v_1, v_2\rangle,
\ldots, \langle v_{k-1}, v_k=v_0\rangle\}$ with corresponding
weights $r_1, r_2, \ldots, r_k$. By Definition~ \ref{EventGraph},
for $i=1,2, \ldots, k$, there exist $T_i\in \mathbb{M}_\infty$ and
$e_i\in\CE$ with $e_i=\sigma_i M_i - \tau_i N_i$ where $\sigma_i,
\tau_i>0$ and $M_i, N_i\in
\mathbb{M}_\infty$ and $M_i \prec N_i$  such that either \ \\
1) $v_{i-1}=T_i M_i$ and $v_i=T_i N_i$ and
$r_i=\ln\frac{\sigma_i}{\tau_i}\in w\left(\langle v_{i-1},
v_i\rangle\right)$ or \ \\ 2) $v_{i-1}=T_i N_i$ and $v_i=T_i M_i$
and $r_i=-\ln\frac{\sigma_i}{\tau_i}\in w\left(\langle v_{i-1},
v_i\rangle\right)$ \
\\ Hence, there exists a vector ${\bs b}=\langle b_1,
b_2, \ldots, b_k\rangle$ with $b_i=0$ or $1$ such that:
\begin{eqnarray}
\prod_{i=1}^k M_i^{b_i} N_i^{1-b_i}=\prod_{i=1}^k
M_i^{1-b_i}N_i^{b_i} \\
w\left(\mathcal{C}\right)=\sum_{i=1}^kr_i=\sum_{i=1}^k \left(2 b_i -
1\right) \ln\left(\frac{\sigma_i}{\tau_i}\right)
\end{eqnarray}
 Let $\bs{\alpha}$ be a strong equilibrium point of $\CE$ that is not a
z-point. Then, by Definition~\ref{Strong Equilibrium Point}, for
$i=1$ to $k$,
$\sigma_iM_i\left(\bs{\alpha}\right)-\tau_iN_i\left(\bs{\alpha}\right)=0$\
\\$\Rightarrow\sigma_iM_i\left(\bs{\alpha}\right)=\tau_iN_i\left(\bs{\alpha}\right)$ for
$i=1$ to $k$\
\\$\Rightarrow\left(\sigma_iM_i\left(\bs{\alpha}\right)\right)^{b_i}=\left(\tau_iN_i\left(\bs{\alpha}\right)\right)^{b_i}$
and
$\left(\tau_iN_i\left(\bs{\alpha}\right)\right)^{1-b_i}=\left(\sigma_iM_i\left(\bs{\alpha}\right)\right)^{1-b_i}$
for $i=1$ to $k$\
\\$\Rightarrow\left(\sigma_iM_i\left(\bs{\alpha}\right)\right)^{b_i}\left(\tau_iN_i\left(\bs{\alpha}\right)\right)^{1-b_i}=\left(\sigma_iM_i\left(\bs{\alpha}\right)\right)^{1-b_i}\left(\tau_iN_i\left(\bs{\alpha}\right)\right)^{b_i}$
for $i=1$ to $k$\
\\$\Rightarrow\prod_{i=1}^k\left(\sigma_iM_i\left(\bs{\alpha}\right)\right)^{b_i}\left(\tau_iN_i\left(\bs{\alpha}\right)\right)^{1-b_i}=\prod_{i=1}^k\left(\sigma_iM_i\left(\bs{\alpha}\right)\right)^{1-b_i}\left(\tau_iN_i\left(\bs{\alpha}\right)\right)^{b_i}$\
\\$\Rightarrow\prod_{i=1}^k{\sigma_i}^{b_i}{\tau_i}^{1-b_i}=\prod_{i=1}^k{\sigma_i}^{1-b_i}{\tau_i}^{b_i}$
[From Equation~(1) and since $\bs{\alpha}$ is not a z-point]\
\\$\Rightarrow\prod_{i=1}^k\frac{{\sigma_i}^{b_i}{\tau_i}^{1-b_i}}{{\sigma_i}^{1-b_i}{\tau_i}^{b_i}}=1$\
\\$\Rightarrow \sum_{i=1}^k \left(2 b_i -
1\right)\ln\left(\frac{\sigma_i}{\tau_i}\right)=0$ [Taking
logarithm]\
\\$\Rightarrow w\left(\mathcal{C}\right)=0$ [From Equation~(2)]\
\\Hence, $\CE$ has no energy cycle.
\
\\$\boldsymbol{\left(3\right)\Rightarrow\left(4\right)}:$ \ \\Let
$\CE=\{\sigma_1M_1-\tau_1N_1, \sigma_2M_2 - \tau_2N_2,
\ldots,\sigma_mM_m - \tau_mN_m\}$ and for all $j=1$ to $m$,
$M_j\prec N_j$ and $\sigma_j, \tau_j>0$. Let $\Gamma=\Gamma_\CE$. We
shall prove that if the linear equation
$\Gamma\bs{\alpha}=\langle\ln\left(\sigma_1/\tau_1\right), \ldots,
\ln\left(\sigma_m/\tau_m\right)\rangle^T$ has no solution in $\MR^n$
then $\CE$ has an energy cycle. For $j=1$ to $m$, let $\Gamma_j$ be
the $j^{th}$ row of $\Gamma$. If the system of linear equations
$\Gamma\bs{\alpha}=\langle\ln\left(\sigma_1/\tau_1\right), \ldots,
\ln\left(\sigma_m/\tau_m\right)\rangle^T$ has no solution in $\MR^n$
then, from linear algebra \cite[p.~164,~Theorem]{narayan03textbook} and the fact that $\Gamma$ is
a matrix of integers, it follows that there exists $l$, there exist
(not necessarily distinct) integers $j_1, j_2, \ldots, j_l\in\{1,2,
\ldots, m\}$, there exist $a_1, a_2, \ldots, a_l\in\{+1,-1\}$ such
that:
\begin{eqnarray}
a_1 \Gamma_{j_1} + a_2 \Gamma_{j_2} + \cdots + a_l \Gamma_{j_l}={\bs
0}
\\a_1 \ln\left(\sigma_{j_1}/\tau_{j_1}\right) + a_2
\ln\left(\sigma_{j_2}/\tau_{j_2}\right) + \cdots + a_l
\ln\left(\sigma_{j_l}/\tau_{j_l}\right)\neq0
\end{eqnarray}

Consider the sequence $\mathcal{C}$ of $l+1$ vertices in the
event-graph defined recursively by
\[
v_0=\prod_{i=1,a_i=+1}^lM_{j_i}\prod_{i=1,a_i=-1}^lN_{j_i}
\] and for $i=1$ to $l$,
\[v_i=\frac{v_{i-1}N_{j_i}^{a_i}}{M_{j_i}^{a_i}}\]
\\Observe that by (3),
\[
\prod_{i=1}^{l}\left(\frac{N_{j_i}}{M_{j_i}}\right)^{a_i}=1
\]
Hence,
\[v_0=\prod_{i=1, a_i=+1}^l
M_{j_i}^{a_i}\prod_{i=1, a_i=-1}^lN_{j_i}^{-a_i} = \prod_{i=1,
a_i=+1}^lN_{j_i}^{a_i}\prod_{i=1, a_i=-1}^lM_{j_i}^{-a_i}=v_l\]
Hence, $\mathcal{C}$ is a cycle. Further, for $i=1$ to $l$,\
\\$a_i \ln\frac{\sigma_{j_i}}{\tau_{j_i}}\in w\left(\langle v_{i-1},
v_i\rangle\right)$\
\\From Equation (4), \[w\left(\mathcal{C}\right)=a_1 \ln\left(\sigma_{j_1}/\tau_{j_1}\right) + a_2
\ln\left(\sigma_{j_2}/\tau_{j_2}\right) + \cdots + a_l
\ln\left(\sigma_{j_l}/\tau_{j_l}\right)\neq0\] Hence, $\mathcal{C}$
is an energy cycle.
\end{proof}

Horn and Jackson~\cite{horn72general} and Feinberg~\cite{feinberg95existence} have proved that chemical reaction networks with appropriate properties admit Lyapunov
functions. While finite, natural event-systems are closely related
to the chemical reaction networks considered by Horn and Jackson and
by Feinberg, they are not identical. Consequently, we will prove the
existence of Lyapunov functions for finite, natural event-systems~(Theorem~\ref{LyapunovExists}).

The Lyapunov function is analogous in form and properties to
``Entropy of the Universe'' in thermodynamics. The Lyapunov function
composed with an event-process is monotonic with respect to time, providing an analogy to the second law of thermodynamics.

\begin{definition}
\label{def:Lyapunov} Let $\CE$ be a finite, natural event-system of
dimension $n$ with positive strong $\CE$-equilibrium point $\bs{c}=\langle c_1,c_2,\ldots, c_n\rangle$. Then
$g_{\CE,\bs{c}}:\MR^n_{>0}\rightarrow\MR$ is given by
\[
g_{\CE,\bs{c}}\left(x_1, x_2, \ldots, x_n\right)=\sum_{i=1}^n
\left(x_i\left(\ln\left(x_i\right)-1-\ln\left(c_i\right)\right)+c_i\right)
\]
\end{definition}

The function $g_{\CE,\bs{c}}$ will turn out to be the desired Lyapunov function.

Note that if $\CE_1$ and $\CE_2$ are two finite natural event-systems of the same dimension and if $\bs{c}$ is a positive strong
$\CE_1$-equilibrium point as well as a positive strong
$\CE_2$-equilibrium point, then the functions $g_{\CE_1,\bs{c}}$ and
$g_{\CE_2,\bs{c}}$ are identical.

\begin{lemma}$\label{OrbitalDerivative}$
Let $\CE=\{\sigma_1M_1 - \tau_1N_1, \sigma_2M_2 - \tau_2N_2, \ldots,
\sigma_mM_m - \tau_mN_m\}$ be a finite, natural event-system of
dimension $n$ with positive strong $\CE$-equilibrium point $\bs{c}$, such
that for all $j=1$ to $m$, $\sigma_j, \tau_j>0$ and $M_j\prec N_j$.
Then for all $\bs{x}\in\MR^n_{>0}$, \[\nabla g_{\CE,
c}\left(\bs{x}\right)\cdot {\bs
P}_\CE\left(\bs{x}\right)=\sum_{j=1}^m\left(\sigma_jM_j\left(\bs{x}\right)-\tau_jN_j\left(\bs{x}\right)
\right)\ln\left(\frac{\tau_jN_j\left(\bs{x}\right)}{\sigma_jM_j\left(\bs{x}\right)}\right)\]
\end{lemma}
\begin{proof}
Let $g=g_{\CE,c}$. Let $\bs{x}=\langle x_1, x_2, \ldots, x_n\rangle\in
\MR^n_{>0}$. Let ${\bs P}={\bs P}_\CE$.
\begin{eqnarray*}
\nabla g\left(\bs{x}\right)\cdot {\bs P}\left(\bs{x}\right) &=&
\sum_{i=1}^n\left(\frac{\partial g}{\partial
x_i}\left(\bs{x}\right)\cdot P_i\left(\bs{x}\right)\right)\\
&=&\sum_{i=1}^n\ln\left(\frac{x_i}{c_i}\right)\left(\sum_{j=1}^m
\gamma_{j,i} \left(\sigma_j M_j\left(\bs{x}\right) - \tau_j
N_j\left(\bs{x}\right)\right)\right)\\
&=&\sum_{j=1}^m\left(\sigma_jM_j\left(\bs{x}\right)-\tau_jN_j\left(\bs{x}\right)
\right)\sum_{i=1}^n\ln\left(\left(\frac{x_i}{c_i}\right)^{\gamma_{j,i}}\right)\\
&=&\sum_{j=1}^m\left(\sigma_jM_j\left(\bs{x}\right)-\tau_jN_j\left(\bs{x}\right)
\right)\ln\left(\prod_{i=1}^n\left(\frac{x_i}{c_i}\right)^{\gamma_{j,i}}\right)\\
&=&\sum_{j=1}^m\left(\sigma_jM_j\left(\bs{x}\right)-\tau_jN_j\left(\bs{x}\right)
\right)\ln\left(\frac{\tau_jN_j\left(\bs{x}\right)}{\sigma_jM_j\left(\bs{x}\right)}\right)
\end{eqnarray*}
The last equality follows from the definition of $\Gamma_\CE$ and
the fact that $\bs{c}$ is a strong-equilibrium point.
\end{proof}

\begin{lemma}\label{ineq} For all $x\in\MR_{>0}$,
$(1-x)\ln\left(x\right)\leq 0$ with equality iff $x=1$.
\end{lemma}
\begin{proof}
If $0<x<1$ then $1-x>0$ and $\ln(x)<0$. If $x>1$ then $1-x<0$ and
$\ln(x)>0$. In either case, the product is strictly negative. If
$x=1$ then $(1-x)\ln\left(x\right)=0$
\end{proof}

\begin{theorem}
\label{SecondLaw} Let $\CE$ be a finite, natural event-system of
dimension $n$ with positive strong $\CE$-equilibrium point ${\bs
c}$. Then for all $\bs{x}\in\MR^n_{> 0}$, $\nabla g_{\CE, {\bs
c}}\left(\bs{x}\right)\cdot {\bs P}_\CE\left(\bs{x}\right)\leq 0$ with
equality iff $\bs{x}$ is a strong $\CE$-equilibrium point.
\end{theorem}
\begin{proof}
Let $\CE=\{\sigma_1M_1 - \tau_1N_1, \sigma_2M_2 - \tau_2N_2, \ldots,
\sigma_mM_m - \tau_mN_m\}$ be a finite, natural event-system of
dimension $n$ with positive strong $\CE$-equilibrium point $\bs{c}$, such
that for all $j=1$ to $m$, $\sigma_j, \tau_j>0$ and $M_j\prec N_j$.
Let ${\bs P}={\bs P}_\CE$ and let $g=g_{\CE,c}$. By Lemma~
\ref{OrbitalDerivative}, for all $\bs{x}\in\MR^n_{> 0}$,
\[
\nabla g\left(\bs{x}\right)\cdot {\bs P}\left(\bs{x}\right) =
\sum_{j=1}^m\left(\sigma_jM_j\left(\bs{x}\right)-\tau_jN_j\left(\bs{x}\right)
\right)\ln\left(\frac{\tau_jN_j\left(\bs{x}\right)}{\sigma_jM_j\left(\bs{x}\right)}\right)
\]
From Lemma~\ref{ineq} and the observation that for $\jtom$,
$M_j\left(\bs{x}\right),N_j\left(\bs{x}\right)>0$ when $\bs{x}\in\MR_{>0}^n$ and by
assumption $\sigma_j,\tau_j>0$, we have,
\[
\nabla g\left(\bs{x}\right)\cdot {\bs P}\left(\bs{x}\right)\leq 0\] with
equality iff for all $\jtom$,
$\sigma_jM_j\left(\bs{x}\right)=\tau_jN_j\left(\bs{x}\right)$. This occurs iff
$\bs{x}$ is a strong $\CE$-equilibrium point.
\end{proof}

Recall that a function $g$ is a Lyapunov function at a point $\bs{p}$ for
a vector field $\bs{v}$ iff $g$ is smooth, positive definite at $\bs{p}$ and
$L_{\bs{v}}g$ is negative semi-definite at $\bs{p}$~\cite[p.\hspace{3pt}131]{irwin80smooth}. For a finite natural event-system $\CE$, ${\bs P}_\CE$ induces a vector
field on $\MR^n$. We will show that, if $\bs{c}$ is a positive strong
$\CE$-equilibrium point, then $g_{\CE,\bs{c}}$ is a Lyapunov function at
$\bs{c}$ for the vector field induced by ${\bs P}_\CE$.

\begin{theorem}[Existence of Lyapunov Function]
\label{LyapunovExists} Let $\CE$ be a finite, natural event-system
of dimension $n$ with positive strong $\CE$-equilibrium point $\bs{c}$.
Then $g_{\CE,\bs{c}}$ is a Lyapunov function for the vector field induced
by ${\bs P}_{\CE}$ at $\bs{c}$.
\end{theorem}
\begin{proof}
Let $g=g_{\CE,\bs{c}}$. For $\iton$:
\[\frac{\partial g}{\partial x_i}=\ln\left(\frac{x_i}{c_i}\right)\] which are
all in $\cal{C}^{\infty}$ as functions on $\MR^n_{>0}$, hence $g$ is in
$\cal{C}^{\infty}$.
\[
\frac{\partial g}{\partial
x_i}\left(\bs{c}\right)=\ln\left(\frac{c_i}{c_i}\right)=0
\] establishes that $\nabla g\left(\bs{c}\right)={\bs 0}$.
For $\iton$, for $k=1,2, \ldots, n$:
\[
\frac{\partial^2g}{\partial x_k \partial
x_i}=\frac{\delta_{i,k}}{x_i}
\]
where $\delta_{i,k}$ is the Kronecker delta function. Hence, for all
$\bs{x}\in\MR^n_{>0}$, the Hessian of $g$ at $\bs x$ is positive definite.
Therefore, $g$ is strictly convex over $\MR^n_{>0}$. Further, $g\left(\bs{c}\right)=0$ and $\nabla g\left(\bs{c}\right)={\bs 0}$ and $g$ is
strictly convex together imply that $g$ is positive definite at $\bs{c}$.
To establish $g$ as a Lyapunov function, it remains to show that the
directional derivative $L_{\bs{P}}g$ of $g$ in the direction of the vector
field induced by ${\bs P}={\bs P}_{\CE}$ is negative semi-definite
at $\bs{c}$. This follows from Theorem~\ref{SecondLaw} since for all
$\bs{x}\in\MR^n_{>0}$, $L_{\bs{P}}g\left(\bs{x}\right)= \nabla g\left(\bs{x}\right) \cdot
{\bs P}\left(\bs{x}\right)\leq 0$.
\end{proof}

Henceforth, the function $g_{\CE,\bs{c}}$ will be called the Lyapunov
function of $\CE$ at $\bs{c}$. The next theorem shows that finite, natural event-systems satisfy a form of ``detailed balance.''

\begin{theorem}\label{alleqstrong}
If $\CE$ is a natural, finite event-system of dimension $n$ then all
positive $\CE$-equilibrium points are strong $\CE$-equilibrium
points.
\end{theorem}
\begin{proof}
Let ${\bs P}={\bs P}_\CE$. Let $\bs{c}\in\MR^n_{>0}$ be a positive strong
$\CE$-equilibrium point. Let $\bs{x}$ be a positive $\CE$-equilibrium
point. That is, ${\bs P}(\bs{x})={\bs 0}$. Hence, $\nabla g_{\CE,
\bs{c}}\left(\bs{x}\right) \cdot {\bs P}_\CE\left(\bs{x}\right)=0$. By Theorem~
\ref{SecondLaw}, $\bs{x}$ is a strong $\CE$-equilibrium point.
\end{proof}

The following lemma was proved by Feinberg~\cite[Proposition~B.1]{feinberg95existence}.

\begin{lemma}\label{convex}
Let $n>0$ be an integer. Let $U$ be a linear subspace of $\MR^n$,
and let ${\bs a}=\langle a_1,a_2, \ldots, a_n\rangle$ and ${\bs b}$
be elements of $\MR^n_{>0}$. There is a unique element
$\bs{\mu}=\langle\mu_1,\mu_2, \cdots, \mu_n\rangle\in U^\perp$ such that
$\langle a_1\mathrm{e}^{\mu_1},a_2\mathrm{e}^{\mu_2}, \ldots
,a_n\mathrm{e}^{\mu_n}\rangle -{\bs b}$ is an element of $U$.
\end{lemma}

The next theorem follows from one proved by Horn and
Jackson~\cite[Lemma~4B]{horn72general}. Our proof is derived
from Feinberg's~\cite[Proposition~5.1]{feinberg95existence}.

\begin{theorem}\label{HornJackson}
Let $\CE$ be a finite, natural event-system of dimension $n$. Let
$H$ be a positive conservation class of $\CE$. Then $H$ contains
exactly one positive strong $\CE$-equilibrium point.
\end{theorem}
\begin{proof}
Let $\Gamma=\Gamma_\CE$. Let ${\bs c^*}=\langle c^*_1,c^*_2, \ldots,
c^*_n\rangle$ be a positive strong $\CE$-equilibrium point. Let
${\bs p}\in H\cap \MR^n_{>0}$. For all ${\bs c}\in\MR^n_{>0}$,\
\\ \\
$(1)$ ${\bs c}$ is a strong $\CE$-equilibrium point \
\\$\Leftrightarrow\Gamma\langle
\ln(c_1),\ln(c_2), \ldots, \ln(c_n)\rangle^T=\Gamma\langle
\ln(c^*_1),\ln(c^*_2),\cdots,\ln(c^*_n)\rangle^T$.
(Lemma~\ref{StrongEqCorr})\
\\$\Leftrightarrow\Gamma\left\langle
\ln\left(\frac{c_1}{c^*_1}\right),\ln\left(\frac{c_2}{c^*_2}\right),
\ldots, \ln\left(\frac{c_n}{c^*_n}\right)\right\rangle^T=0$\
\\$\Leftrightarrow$ There exists $\mu=\langle\mu_1,\mu_2, \ldots,
\mu_n\rangle\in\ker\Gamma\cap\MR^n$ such that\
\\$\left\langle\ln\left(\frac{c_1}{c^*_1}\right),\ln\left(\frac{c_2}{c^*_2}\right),
\ldots, \ln\left(\frac{c_n}{c^*_n}\right)\right\rangle^T={\bs \mu}$.\
\\$\Leftrightarrow$ There exists $\bs{\mu}=\langle\mu_1,\mu_2, \ldots,
\mu_n\rangle\in\ker\Gamma\cap\MR^n$ such that
$c_i=c^*_i\mathrm{e}^{\mu_i}$ for $\iton$.\
\\ \\
$(2)$ ${\bs c}\in H\cap\MR^n\Leftrightarrow{\bs c}-{\bs
p}\in\left(\ker\Gamma\right)^{\perp}\cap\MR^n$. (By
Definition~\ref{Conservation Class})\
\\ \\
From (1) and (2), ${\bs c}$ is a positive strong $\CE$-equilibrium
point in $H$ iff there exists $\bs{\mu}\in\ker\Gamma\cap\MR^n$ such that
${\bs c}=\langle c^*_1\mathrm{e}^{\mu_1}, c^*_2\mathrm{e}^{\mu_2},
\ldots, c^*_n\mathrm{e}^{\mu_n}\rangle$ and\ \\$\langle
c^*_1\mathrm{e}^{\mu_1}, c^*_2\mathrm{e}^{\mu_2}, \ldots,
c^*_n\mathrm{e}^{\mu_n}\rangle-{\bs p} \in
\left(\ker\Gamma\right)^{\perp}\cap\MR^n$. Applying
Lemma~\ref{convex} with ${\bs a}={\bs c^*},{\bs b}={\bs p}$ and
$U=(\ker\Gamma)^{\perp}\cap\MR^n$, it follows that there exists a
unique $\bs{\mu}$ of the desired form. Hence, there exists a unique
positive strong $\CE$-equilibrium point in $H$ given by ${\bs
c}=\langle c^*_1\mathrm{e}^{\mu_1}, c^*_2\mathrm{e}^{\mu_2}, \ldots,
c^*_n\mathrm{e}^{\mu_n}\rangle$.
\end{proof}

To prove the main theorem of this section~(Theorem~\ref{AsymptoticStability}), we will first establish several technical lemmas.

Lemma~\ref{reduce} shows that an event that remains zero at all times along a process can be ignored.

\begin{lemma}\label{reduce}
Let $\CE$ be a finite event-system of dimension $n$, let
$\Omega\subseteq \MC$ be non-empty, open and simply-connected, and
let $\boldsymbol{f}=\langle f_1, f_2,\ldots, f_n\rangle$ be an
$\CE$-process on $\Omega$. Then either for all $t \in \Omega$,
$\bs{f}(t)$ is a strong $\CE$-equilibrium point or there exist a
finite event-system $\hat{\CE}$ of dimension $\hat{n}\leq n$, an
$\hat{\CE}$-process $\boldsymbol{\hat{f}}=\langle \hat{f_1},
\hat{f_2}, \ldots, \hat{f}_{\hat{n}} \rangle$ on $\Omega$, and a
permutation $\pi$ on $\{1,2,\ldots,n\}$ such that:
\begin{enumerate}
\item\label{reduce.a} If $\CE$ is physical then $\hat{\CE}$ is physical.
\item\label{reduce.b} If $\CE$ is natural then $\hat{\CE}$ is natural.
\item\label{reduce.b2} If $\bs{c}=\langle c_1, c_2, \ldots, c_n \rangle$ is a positive strong $\CE$-equilibrium point, then $\bs{\hat{c}}=\langle c_{\pi^{-1}(1)}, c_{\pi^{-1}(2)}, \ldots, c_{\pi^{-1}(\hat{n})} \rangle$ is a positive strong $\hat{\CE}$-equilibrium point.
\item\label{reduce.c} For all $e\in\hat{\CE}$, there exists $t\in\Omega$ such that
$e(\boldsymbol{\hat{f}}(t))\neq0$.
\item\label{reduce.d} If $\hat{\CE}$ is natural, $I\subseteq\Omega\cap\MR_{\geq 0}$ is connected, $0\in I$ and
$\boldsymbol{f}(0)$ is a non-negative point then for all $t\in
I\cap\MR_{>0}$, $\boldsymbol{\hat{f}}(t)$ is a positive point.
\item\label{reduce.e} For $\iton$, if $\pi(i)\leq\hat{n}$ then for all $t\in\Omega$, $f_i(t)=\hat{f}_{\pi(i)}(t)$.
\item\label{reduce.f} For $\iton$, if $\pi(i)>\hat{n}$ then for all $t_1,t_2\in\Omega$, $f_i(t_1)=f_i(t_2)$.
\end{enumerate}
\end{lemma}

\begin{proof}
Let $m=|\CE|$. Let
$\CE_1=\{e\in\CE\mid$~there~exists~$t\in\Omega,e(\boldsymbol{f}(t))\neq0\}$.
If $\CE_1 = \varnothing$ then for all $t \in \Omega,
e(\bs{f}(t))=0$, so $\bs{f}(t)$ is a strong $\CE$-equilibrium point
and the Lemma holds. Assume $\CE_1 \neq \varnothing$ and let
$\hat{m}=|\CE_1|$. For $j=1,2,\ldots,m$, let
$\sigma_j,\tau_j\in\MR_{>0}$ and
$M_j=\prod_{i=1}^{n}X_i^{a_{j,i}},N_j=\prod_{i=1}^{n}X_i^{b_{j,i}}\in
\mathbb{M}_\infty$ be such that $M_j\prec N_j$ and
$\{\sigma_1M_1-\tau_1N_1,\sigma_2M_2-\tau_2N_2,\ldots,\sigma_{\hat{m}}M_{\hat{m}}-\tau_{\hat{m}}N_{\hat{m}}\}=\CE_1$
and
$\{\sigma_1M_1-\tau_1N_1,\sigma_2M_2-\tau_2N_2,\ldots,\sigma_mM_m-\tau_mN_m\}=\CE$.

Let
$C=\{i\mid$~there~exists~$j\leq\hat{m}$~such~that~either~$a_{j,i}\neq0$~or~$b_{j,i}\neq0\}$.
Let $\hat{n}=|C|$. Let $\pi$ be a permutation on $\{1,2,\ldots,n\}$
such that $\pi(C)=\{1,2,\ldots,\hat{n}\}$.

For $j=1,2,\ldots,\hat{m}$, let
$e_{\pi,j}=\sigma_j\prod_{i=1}^{\hat{n}}X_i^{a_{j,{\pi^{-1}(i)}}}-\tau_j\prod_{i=1}^{\hat{n}}X_i^{b_{j,{\pi^{-1}(i)}}}$.
Let $\hat{\CE}=\{e_{\pi,1},e_{\pi,2},\ldots,e_{\pi,\hat{m}}\}$.

It follows that $\hat{\CE}$ is a finite event-system of dimension
$\hat{n}\leq n$. For $i=1,2,\ldots,\hat{n}$, let
$\hat{f}_i=f_{\pi^{-1}(i)}$. Let $\boldsymbol{\hat{f}}=\langle
\hat{f}_1, \hat{f}_2,\ldots,\hat{f}_{\hat{n}}\rangle$.

Let $(\gamma_{j,i})_{m\times n}=\Gamma_\CE$. Let
$(\hat{\gamma}_{j,i})_{\hat{m}\times\hat{n}}=\Gamma_{\hat{\CE}}$. It
follows that for $j=1,2,\ldots,\hat{m}$, for $i=1,2,\ldots,\hat{n}$,
\begin{align}\label{reduce.1}
\hat{\gamma}_{j,i}=b_{j,\pi^{-1}(i)}-a_{j,\pi^{-1}(i)}=\gamma_{j,\pi^{-1}(i)}.
\end{align}

We claim that $\boldsymbol{\hat{f}}$ is an $\hat{\CE}$-process on
$\Omega$. To see this, for $k=1,2,\ldots,\hat{n}$, for all
$t\in\Omega:$
\begin{align*}
\hat{f}_k'(t) &= f_{\pi^{-1}(k)}'(t) \text{\qquad[Definition of $\hat{f}_k$.]}\\
&= \left[\left(\sum_{j=1}^m\gamma_{j,\pi^{-1}(k)}\left(\sigma_j\prod_{i=1}^nX_i^{a_{j,i}}-\tau_j\prod_{i=1}^nX_i^{b_{j,i}}\right)\right)\circ\boldsymbol{f}\right]\left(t\right)\qquad[\boldsymbol{f}\text{ is an $\CE$-process on }\Omega.]\\
&= \left[\left(\sum_{j=1}^{\hat{m}}\gamma_{j,\pi^{-1}(k)}\left(\sigma_j\prod_{i=1}^nX_i^{a_{j,i}}-\tau_j\prod_{i=1}^nX_i^{b_{j,i}}\right)\right)\circ\boldsymbol{f}\right]\left(t\right)\qquad\text{[Definition of $\CE_1$.]}\\
&=\left[\left(\sum_{j=1}^{\hat{m}}\gamma_{j,\pi^{-1}(k)}\left(\sigma_j\prod_{i\in
C}X_i^{a_{j,i}}-\tau_j\prod_{i\in
C}X_i^{b_{j,i}}\right)\right)\circ\boldsymbol{f}\right]\left(t\right) \qquad[j\leq\hat{m},i\notin C\Rightarrow a_{j,i}=b_{j,i}=0.]\\
\end{align*}
\begin{align*}
&= \left[\left(\sum_{j=1}^{\hat{m}}\gamma_{j,\pi^{-1}(k)}\left(\sigma_j\prod_{i=1}^{\hat{n}}X_{\pi^{-1}(i)}^{a_{j,{\pi^{-1}(i)}}}-\tau_j\prod_{i=1}^{\hat{n}}X_{\pi^{-1}(i)}^{b_{j,{\pi^{-1}(i)}}}\right)\right)\circ\boldsymbol{f}\right]\left(t\right)\qquad[\pi(C)=\{1,2,\ldots,\hat{n}\}.]\\
&= \sum_{j=1}^{\hat{m}}\gamma_{j,\pi^{-1}(k)}\left(\sigma_j\prod_{i=1}^{\hat{n}}(f_{\pi^{-1}(i)}(t))^{a_{j,{\pi^{-1}(i)}}}-\tau_j\prod_{i=1}^{\hat{n}}(f_{\pi^{-1}(i)}(t))^{b_{j,{\pi^{-1}(i)}}}\right)\qquad\text{[By composition.]}\\
&= \sum_{j=1}^{\hat{m}}\hat{\gamma}_{j,k}\left(\sigma_j\prod_{i=1}^{\hat{n}}(f_{\pi^{-1}(i)}(t))^{a_{j,{\pi^{-1}(i)}}}-\tau_j\prod_{i=1}^{\hat{n}}(f_{\pi^{-1}(i)}(t))^{b_{j,{\pi^{-1}(i)}}}\right)\qquad\text{[From~(\ref{reduce.1}).]}\\
&= \sum_{j=1}^{\hat{m}}\hat{\gamma}_{j,k}\left(\sigma_j\prod_{i=1}^{\hat{n}}(\hat{f}_i(t))^{a_{j,{\pi^{-1}(i)}}}-\tau_j\prod_{i=1}^{\hat{n}}(\hat{f}_i(t))^{b_{j,{\pi^{-1}(i)}}}\right)    \qquad\text{[Definition of $\hat{f}_i$.]}\\
&=\left[\left(\sum_{j=1}^{\hat{m}}\hat{\gamma}_{j,k}e_{\pi,j}\right)\circ\boldsymbol{\hat{f}}\right](t)
\qquad\text{[Definition of }e_{\pi,j}.]
\end{align*}
This establishes the claim.

With $\hat{\CE},\hat{n},\boldsymbol{\hat{f}}$ and $\pi$ as
described,
we will now establish (1) through (6).\\\
\\(1) Follows from the definition of $\hat{\CE}$.\\\
\\(2) Follows from 3.
\\(3) Suppose $\CE$ is natural. Hence, there
exists a positive strong $\CE$-equilibrium point $\langle
c_1,c_2,\ldots,c_n\rangle$.  For $j=1,2,\ldots,\hat{m}:$
\begin{align*}
e_{\pi,j}(c_{\pi^{-1}(1)},c_{\pi^{-1}(2)},\ldots,c_{\pi^{-1}(\hat{n})})&=
\sigma_j\displaystyle\prod_{i=1}^{\hat{n}}c_{\pi^{-1}(i)}^{a_{j,{\pi^{-1}(i)}}}-\tau_j\prod_{i=1}^{\hat{n}}c_{\pi^{-1}(i)}^{b_{j,{\pi^{-1}(i)}}}\\
&=\sigma_j\prod_{i\in C}c_i^{a_{j,i}}-\tau_j\prod_{i\in
C}c_i^{b_{j,i}}\qquad[j\leq\hat{m},i\notin C\Rightarrow a_{j,i}=b_{j,i}=0.]\\
&=e_j(c_1,c_2,\ldots,c_n)\\
&=0
\end{align*}
Hence, $\bs{\hat{c}}$ is a positive strong $\hat{\CE}$-equilibrium point.\\\
\\(4) Suppose $j\leq\hat{m}$. Then for all $t\in\Omega:$
\begin{align*}
e_{\pi,j}(\boldsymbol{\hat{f}}(t))&=\sigma_j\prod_{i=1}^{\hat{n}}(\hat{f}_i(t))^{a_{j,{\pi^{-1}(i)}}}-\tau_j\prod_{i=1}^{\hat{n}}(\hat{f}_i(t))^{b_{j,{\pi^{-1}(i)}}}\\
&=\sigma_j\prod_{i=1}^{\hat{n}}(f_{\pi^{-1}(i)}(t))^{a_{j,{\pi^{-1}(i)}}}-\tau_j\prod_{i=1}^{\hat{n}}(f_{\pi^{-1}(i)}(t))^{b_{j,{\pi^{-1}(i)}}}\\
&=\sigma_j\prod_{i\in C}(f_i(t))^{a_{j,i}}-\tau_j\prod_{i\in
C}(f_i(t))^{b_{j,i}}\\
&=\sigma_j\prod_{i=1}^n(f_i(t))^{a_{j,i}}-\tau_j\prod_{i=1}^n(f_i(t))^{b_{j,i}}\qquad[j\leq\hat{m},i\notin C\Rightarrow a_{j,i}=b_{j,i}=0.]\\
&=\left(\left(\sigma_j\prod_{i=1}^nX_i^{a_{j,i}}-\tau_j\prod_{i=1}^nX_i^{b_{j,i}}\right)\circ\boldsymbol{f}\right)(t)\\
&=e_j(\boldsymbol{f}(t))
\end{align*}
Since $j\leq\hat{m}$, therefore $e_j\in\CE_1$ and there exists
$t\in\Omega$ such that $e_j(\boldsymbol{f}(t))\neq0$. Hence, for all
$e_{\pi,j}\in\hat{\CE}$, there exists $t\in\Omega$ such that
$e_{\pi,j}(\boldsymbol{\hat{f}}(t))\neq0$.\\\
\\(5) Suppose
$\hat{\CE}$ is natural, $I\subseteq\Omega\cap\MR_{\geq 0}$ is
connected, $0\in I$ and $\boldsymbol{f}(0)$ is a non-negative point.
It follows that $\boldsymbol{\hat{f}}(0)$ is a non-negative point
and, from Theorem~\ref{thm:stays_positive}, for all $t\in I$,
$\boldsymbol{\hat{f}}(t)$ is a non-negative point. Suppose, for the
sake of contradiction, that there exist $i_0\leq\hat{n}$ and $t_0\in
I\cap\MR_{>0}$ such that $\hat{f}_{i_0}(t_0)=0$. From
Theorem~\ref{thm:stays_positive} again, $\hat{f}_{i_0}(0)=0$ and for
all $t\in I:\hat{f}_{i_0}(t)=0$. Since $I$ is an interval and
$0,t_0\in I$, $I$ contains an accumulation point. Hence, since
$\hat{f}_{i_0}$ is analytic on $\Omega$ and $\Omega$ is connected,
for all $t\in\Omega:$
\begin{align}
\hat{f}_{i_0}(t)=0.\label{8}
\end{align}
It follows that for all $t\in\Omega:$
\begin{align}
0=\hat{f}'_{i_0}(t)=\sum_{j=1}^{\hat{m}}\hat{\gamma}_{j,i_0}e_{\pi,j}(\boldsymbol{\hat{f}}(t)).
\label{5}
\end{align}
\
\\We claim that for $j=1,2,\ldots,\hat{m}$, for all
$t\in\Omega:\hat{\gamma}_{j,i_0}e_{\pi,j}(\boldsymbol{\hat{f}}(t))\geq0$.\\\
\\Case 1: Suppose $\hat{\gamma}_{j,i_0}=0$. Then
$\hat{\gamma}_{j,i_0}e_{\pi,j}(\boldsymbol{\hat{f}}(t))=0\geq0$.\\\
\\Case 2: Suppose $\hat{\gamma}_{j,i_0}>0$. Then $b_{j,\pi^{-1}(i_0)}>0$.
Hence,
\begin{align*}
e_{\pi,j}(\boldsymbol{\hat{f}}(t))&=\sigma_j\prod_{i=1}^{\hat{n}}(\hat{f}_i(t))^{a_{j,{\pi^{-1}(i)}}}-\tau_j\prod_{i=1}^{\hat{n}}(\hat{f}_i(t))^{b_{j,{\pi^{-1}(i)}}}\\
&=\sigma_j\prod_{i=1}^{\hat{n}}(\hat{f}_i(t))^{a_{j,{\pi^{-1}(i)}}}\qquad[\text{Since $b_{j,\pi^{-1}(i_0)}>0$ and from~\ref{8}, }\hat{f}_{i_0}(t)=0.]\\
&\geq0\qquad\text{[$\boldsymbol{\hat{f}}(t)$ is a non-negative
point, by Theorem~\ref{thm:stays_positive}]}
\end{align*}
Hence,
$\hat{\gamma}_{j,i_0}e_{\pi,j}(\boldsymbol{\hat{f}}(t))\geq0$.\\\
\\Case 3: Suppose $\hat{\gamma}_{j,i_0}<0$. Then $a_{j,\pi^{-1}(i_0)}>0$.
Hence,
\begin{align*}
e_{\pi,j}(\boldsymbol{\hat{f}}(t))&=\sigma_j\prod_{i=1}^{\hat{n}}(\hat{f}_i(t))^{a_{j,{\pi^{-1}(i)}}}-\tau_j\prod_{i=1}^{\hat{n}}(\hat{f}_i(t))^{b_{j,{\pi^{-1}(i)}}}\\
&=-\tau_j\prod_{i=1}^{\hat{n}}(\hat{f}_i(t))^{b_{j,{\pi^{-1}(i)}}}\qquad[\text{Since $a_{j,\pi^{-1}(i_0)}>0$ and from~\ref{8}, }\hat{f}_{i_0}(t)=0.]\\
&\leq0\qquad\text{[$\boldsymbol{\hat{f}}(t)$ is a non-negative
point, by Theorem~\ref{thm:stays_positive}]}
\end{align*}
Hence,
$\hat{\gamma}_{j,i_0}e_{\pi,j}(\boldsymbol{\hat{f}}(t))\geq0$. This
completes the proof of the claim.

From \ref{5} and the claim, it now follows that for
$j=1,2,\ldots,\hat{m}$, for all $t\in\Omega:$
\begin{align}
\hat{\gamma}_{j,i_0}e_{\pi,j}(\boldsymbol{\hat{f}}(t))=0\label{9}
\end{align}
Since $i_0\leq\hat{n}$, there exists $j_0\leq\hat{m}$ such that
either $a_{j_0,i_0}\neq0$ or $b_{j_0,i_0}\neq0$. If
$\hat{\gamma}_{j_0,i_0}\neq0$ then, from \ref{9},
$e_{\pi,j_0}(\boldsymbol{\hat{f}}(t))=0$. If
$\hat{\gamma}_{j_0,i_0}=0$ then, since
$\hat{\gamma}_{j_0,i_0}=b_{j_0,i_0}-a_{j_0,i_0}$, it follows that
$a_{j_0,i_0}\neq0$ and $b_{j_0,i_0}\neq0$. Hence, $X_{i_0}$ divides
$e_{\pi,j_0}$. From \ref{8}, it follows that
$e_{\pi,j_0}(\boldsymbol{\hat{f}}(t))=0$. Hence, irrespective of the
value of $\hat{\gamma}_{j_0,i_0}$, for all
$t\in\Omega:e_{\pi,j_0}(\boldsymbol{\hat{f}}(t))=0$. Since
$e_{\pi,j_0}$ is an element of $\hat{\CE}$, this leads to a
contradiction with Lemma~\ref{reduce}.\ref{reduce.c}. Hence, for all
$i\leq\hat{n}$, for all $t\in I\cap\MR_{>0}:$ $\hat{f}_i(t)>0$.\\\
\\(6) Follows from the definition of $\boldsymbol{\hat{f}}$.\\\
\\(7) For $\iton$, if $\pi(i)>\hat{n}$
then $i\notin C$. That is, for $\jtom:
\gamma_{j,i}=b_{j,i}-a_{j,i}=0-0=0.$ Hence, for all
$t\in\Omega:f_i'(t)=\sum_{j=1}^m\gamma_{j,i}e_j(\boldsymbol{f}(t))=0$.
Hence, since $f_i$ is analytic on $\Omega$, and $\Omega$ is
simply-connected, for all $t_1,t_2\in\Omega: f_i(t_1)=f_i(t_2)$.
\end{proof}

We have described, for finite, natural event-systems, Lyapunov functions on the positive orthant. We next extend the definition of these Lyapunov functions to admit values at non-negative points.

\begin{definition}\label{extended lyap}
Let $\CE$ be a finite, natural event-system of dimension $n$ with
positive strong $\CE$-equilibrium point $\boldsymbol{c} = \langle
c_1, c_2, \ldots, c_n \rangle$. For all $v\in \MR_{>0}$, let
$\overline{g_v} : \MR_{\geq 0} \rightarrow \MR$ be such that for all
$x \in \MR_{\geq 0}$
\begin{equation}\label{gv}
  \overline{g_v}(x)=
  \begin{cases}
    x (\ln(x) - 1 - \ln(v)) + v,    &\text{if $x > 0$;} \\
    v,                              &\text{otherwise.}
  \end{cases}
\end{equation}
Then the \emph{extended lyapunov function}
$\overline{g_{\CE,\boldsymbol{c}}} : \MR^n_{\geq 0} \rightarrow \MR$
is
\begin{equation}\label{2}
  \overline{g_{\CE,\boldsymbol{c}}}(x_1, x_2, \ldots, x_n) = \sum_{i=1}^n \overline{g_{c_i}}(x_i)
\end{equation}
\end{definition}

The next lemma lists some properties of extended Lyapunov functions.

\begin{lemma}\label{ex lyap properties}
Let $\CE$ be a finite, natural event-system of dimension $n$ with
positive strong $\CE$-equilibrium point $\boldsymbol{c} = \langle
c_1, c_2, \ldots, c_n \rangle$. Then:
\begin{enumerate}
\item\label{ex lyap cont} $\overline{g_{\CE,\boldsymbol{c}}}$ is continuous on
$\MR^n_{\geq 0}$.
\item\label{ex lyap positive} For all $x_1,x_2,\ldots,x_n\in \MR_{\geq 0}$, $\overline{g_{\CE,\boldsymbol{c}}}(x_1,x_2,\ldots,x_n) \geq 0$ with equality iff $\langle x_1,x_2,\ldots,x_n\rangle=\boldsymbol{c}$.
\item\label{xlnx}For all $r\in\MR_{\geq 0}$, the set
$\{\boldsymbol{x}\in\MR^n_{\geq0}\mid
\overline{g_{\CE,\boldsymbol{c}}}(\boldsymbol{x})\leq r\}$ is
bounded.
\item\label{gmon}If $\Omega \subseteq \MC$ is open, simply connected and such that $0 \in \Omega$, $\boldsymbol{f} = \langle f_1, f_2, \ldots, f_n \rangle$ is an $\CE$-process on $\Omega$ such that $\boldsymbol{f}(0)$ is a non-negative point, and $I \subseteq \MR_{\geq 0} \cap \Omega$ is an interval such that $0 \in I$ then $(\overline{g_{\CE,\boldsymbol{c}}} \circ \boldsymbol{f})$ is monotonically non-increasing on $I$.
\end{enumerate}
\end{lemma}
\begin{proof}
For $\iton$, let $\overline{g_{c_i}}(x)$ be as defined in
Equation~\ref{gv}.\
\\
\\1. For $\iton,\overline{g_{c_i}}$ is continuous on $\MR_{>0}$ and
$\lim_{x\to 0^+}\overline{g_{c_i}}(x) = c_i =
\overline{g_{c_i}}(0)$, so $\overline{g_{c_i}}$ is continuous on
$\MR_{\geq 0}$. Since $\overline{g_{\CE,\boldsymbol{c}}}$ is the
finite sum of continuous functions on $\MR_{\geq0}$,
$\overline{g_{\CE,\boldsymbol{c}}}$ is continuous on $\MR^n_{\geq
0}$. \
\\\\2. Let $j\in\{1,2,\ldots,n\}$. Let $\overline{g}=\overline{g_{c_j}}$.
For all $x\in\MR_{>0}$,
$\overline{g}\,'(x)=\ln\left(\displaystyle\frac{x}{c_j}\right)$. If
$0<x<c_j$ then, by substitution, $\overline{g}\,'(x)<0$. Similarly, if
$x>c_j$ then $\overline{g}\,'(x)>0$. Hence, $\overline{g}$ is
monotonically decreasing in $(0,c_j)$ and monotonically increasing
in $(c_j,\infty)$. From continuity of $\overline{g}$ in $\MR_{\geq
0}$, it follows that
\begin{equation}\label{A}
\text{For all }x\in\MR_{\geq0},\overline{g}(x)\geq
\overline{g}(c_j)=0\text{ with equality iff }x=c_j.
\end{equation}
From Equations~(\ref{2}) and (\ref{A}), the claim follows.\
\\
\\3. Observe that $\lim_{x\to +\infty}\overline{g}(x)=+\infty.$
It follows that:
\begin{equation}\label{B}
\text{For all }r\in\MR_{\geq 0},\text{ the set }\{x\in\MR_{\geq
0}\mid \overline{g}(x)\leq r\}\text{ is bounded.}
\end{equation}
If $x_1,x_2,\ldots,x_n\in\MR_{\geq0}$ are such that
$\overline{g_{\CE,\boldsymbol{c}}}(x_1,x_2,\ldots,x_n)\leq r$, it
follows from Equations~(\ref{2}) and (\ref{A}) that for
$\iton:\overline{g_{c_i}}(x_i)\leq r$. The claim now follows from
Equation~(\ref{B}).\
\\
\\4. Let $\Omega \subseteq \MC$ be open, simply connected, and
such that $0 \in \Omega$; let $\bs{f} = \langle f_1, f_2, \ldots,
f_n \rangle$ be an $\CE$-process on $\Omega$ such that $\bs{f}(0)$
is a non-negative point; and let $I \subseteq \MR_{\geq 0} \cap
\Omega$ be an interval such that $0 \in I$. By Lemma~\ref{reduce}
there exists $\hat{n}$, $\hat{\CE}$, $\bs{\hat{f}}$, and $\pi$
satisfying
\ref{reduce}.\ref{reduce.a}\nobreakdash--\ref{reduce}.\ref{reduce.f}.
Let $\bs{\hat{c}} = \langle \hat{c}_1, \hat{c}_2, \ldots,
\hat{c}_{\hat{n}} \rangle = \langle c_{\pi^{-1}(1)},
c_{\pi^{-1}(2)}, \ldots, c_{\pi^{-1}(\hat{n})} \rangle$. By
Lemma~\ref{reduce}.\ref{reduce.b}, $\bs{\hat{c}}$ is a positive
strong equilibrium point of $\hat{\CE}$.  Then for all $t \in I$,
\begin{align*}
\left(\overline{g_{\CE,\bs{c}}}\circ\bs{f}\right)(t)
 &= \sum_{i=1}^n \overline{g_{c_i}}\left(f_i(t)\right) &&\text{[Equation (\ref{2}).]} \\
 &= \sum_{i:\pi(i)\leq\hat{n}} \overline{g_{c_i}}\left(f_i(t)\right) +
    \sum_{i:\pi(i)>\hat{n}} \overline{g_{c_i}}\left(f_i(t)\right) \\
 &= \sum_{i=1}^{\hat{n}} \overline{g_{c_{\pi^{-1}(i)}}}\left(f_{\pi^{-1}(i)}(t)\right) +
    \sum_{i:\pi(i)>\hat{n}} \overline{g_{c_i}}\left(f_i(t)\right) \\
 &= \sum_{i=1}^{\hat{n}} \overline{g_{\hat{c}_i}}\left(\hat{f}_i(t)\right) +
    \sum_{i:\pi(i)>\hat{n}} \overline{g_{c_i}}\left(f_i(t)\right) &&\text{[Definition of $\bs{\hat{c}}$ and Lemma~\ref{reduce}.\ref{reduce.e}.]}\\
 &= \left(\overline{g_{\hat{\CE},\bs{\hat{c}}}}\circ\bs{\hat{f}}\right)(t) +
    \sum_{i:\pi(i)>\hat{n}} \overline{g_{c_i}}\left(f_i(t)\right) &&\text{[Equation (\ref{2}).]} \\
 &= \left(\overline{g_{\hat{\CE},\bs{\hat{c}}}}\circ\bs{\hat{f}}\right)(t) + \text{constant} &&\text{[Lemma \ref{reduce}.\ref{reduce.f}.]}
\end{align*}

By Definition~\ref{extended lyap}, for all
$\bs{x}\in\MR_{>0}^{\hat{n}}$,
$\overline{g_{\hat{\CE},\bs{\hat{c}}}}(\bs{x})=g_{\hat{\CE},\bs{\hat{c}}}(\bs{x})$.
By Lemma~\ref{reduce}.\ref{reduce.d}, for all $t \in I \cap
\MR_{>0}$, $\bs{\hat{f}}(t) \in \MR_{>0}^{\hat{n}}$. So for all $t
\in I \cap \MR_{>0}$,
$\left(\overline{g_{\hat{\CE},\bs{\hat{c}}}}\circ\bs{\hat{f}}\right)(t)
 = \left(g_{\hat{\CE},\bs{\hat{c}}}\circ\bs{\hat{f}}\right)(t)$. Then, for all $t \in I \cap \MR_{>0}$,
\begin{align*}
\left(\overline{g_{\hat{\CE},\bs{\hat{c}}}}\circ\bs{\hat{f}}\right)'(t)
 &= \left(g_{\hat{\CE},\bs{\hat{c}}}\circ\bs{\hat{f}}\right)'(t) \\
 &= \nabla g_{\hat{\CE},\bs{\hat{c}}}\left(\bs{\hat{f}}(t)\right)\cdot\bs{\hat{f}}\,'(t) &&\text{[Chain rule.]} \\
 &= \nabla g_{\hat{\CE},\bs{\hat{c}}}\left(\bs{\hat{f}}(t)\right)\cdot\bs{P}_{\hat{\CE}}\left(\bs{\hat{f}}(t)\right)
 &&\text{[Definition \ref{Process}.]} \\
 &\leq 0
 && \text{[Theorem~\ref{SecondLaw}.]}
\end{align*}
Therefore
$\left(\overline{g_{\hat{\CE},\bs{\hat{c}}}}\circ\bs{\hat{f}}\right)$
is non-increasing on $I \cap \MR_{>0}$.

By Definition~\ref{Process}, $\bs{\hat{f}}$ is continuous on $I$; by
Theorem~\ref{thm:stays_positive}, $\bs{\hat{f}}(I) \subseteq
\MR^{\hat{n}}_{\geq0}$; and by Lemma~\ref{ex lyap
properties}.\ref{ex lyap cont},
$\overline{g_{\hat{\CE},\bs{\hat{c}}}}$ is continuous on
$\MR^{\hat{n}}_{\geq0}$; so
$\left(\overline{g_{\hat{\CE},\bs{\hat{c}}}}\circ\bs{\hat{f}}\right)$
is continuous on $I$. Therefore
$\left(\overline{g_{\hat{\CE},\bs{\hat{c}}}}\circ\bs{\hat{f}}\right)$
is non-increasing on $I$.  Thus
$(\overline{g_{\CE,\bs{c}}}\circ\bs{f})$ is a constant plus a
monotonically non-increasing function on $I$, so
$(\overline{g_{\CE,\bs{c}}}\circ\bs{f})$ is monotonically
non-increasing on $I$.
\end{proof}

The next lemma makes use of properties of the extended Lyapunov function to show that $\CE$-processes starting at non-negative points are uniformly bounded in forward real time.

\begin{lemma}\label{bounded}
Let $\CE$ be a finite, natural event-system of dimension $n$. Let
$\bs{\alpha}\in\MR^n_{\geq 0}$. There exists $k\in\MR_{\geq0}$ such
that for all $\Omega\subseteq \MC$ open and simply connected and
such that $0\in\Omega$, for all $\CE$-processes
$\boldsymbol{f}=\langle f_1, f_2, \ldots, f_n \rangle$ on $\Omega$
such that $\bs{f}(0)=\bs{\alpha}$, for all intervals $I\subseteq
\Omega\cap\MR_{\geq 0}$ such that $0 \in I$, for all $t\in I$, for
$i=1,2,\ldots,n$: $f_i(t)\in\MR$ and $0\leq f_i(t)<k$.
\end{lemma}
\begin{proof}
Since $\CE$ is natural, let $\boldsymbol{c}\in\MR^n_{>0}$ be a
positive strong $\CE$-equilibrium point. Let
$\overline{g}=\overline{g_{\CE,\boldsymbol{c}}}$.

Let $\ell=\overline{g}(\boldsymbol{\alpha})$. Let
$S=\{\boldsymbol{x}\in\MR^n_{\geq 0}\mid
\overline{g}(\boldsymbol{x})\leq \ell\}$. By Lemma~\ref{ex lyap
properties}.\ref{xlnx}, $S$ is bounded. Hence, let $k$ be such that
for all $\bs{x}\in S:|\bs{x}|_{\infty} < k$.

Let $\Omega \subseteq \MC$ be open, simply connected, and such that
$0 \in \Omega$; let $\bs{f} = \langle f_1, f_2, \ldots, f_n \rangle$
be an $\CE$-process on $\Omega$ such that $\bs{f}(0)=\bs{\alpha}$;
and let $I \subseteq \MR_{\geq 0} \cap \Omega$ be an interval such
that $0 \in I$.

From Theorem~\ref{thm:stays_positive}, for all $t\in I$, for $\iton:
f_i(t)\in\MR$ and $f_i(t)\geq 0$.

Consider the function:
\[\overline{g}\circ{\boldsymbol{f}}|_{I}:I\rightarrow\MR\]
From Lemma~\ref{ex lyap properties}.\ref{gmon}, for all $t\in I$,
$\overline{g}\circ{\boldsymbol{f}}|_I$ is monotonically
non-increasing on $I$. That is, for all $t\in I$,
\begin{equation}\label{1}
\overline{g}({\boldsymbol{f}}(t))\leq \ell
\end{equation}

It follows from Equation~\ref{1} and the definition of $S$ that
$\boldsymbol{f}(I)\subseteq S$. By the definition of $k$, it follows
that for all $t\in I$, for $\iton$, $f_i(t)<k$.
\end{proof}

The next lemma shows that, because $\CE$-processes starting at non-negative points are uniformly bounded in real time, they can be continued forever along forward real time.

\begin{lemma}[Existence and uniqueness of $\CE$-process.]
\label{lem:ex&uni} Let $\CE$ be a finite, natural event-system of
dimension $n$. Let $\bs\alpha\in\MR^n_{\geq 0}$. There exist a
simply-connected open set $\Omega\subseteq \MC$, an $\CE$-process
${\boldsymbol{f}}=\langle f_1,f_2,\ldots,f_n\rangle$ on $\Omega$ and
$k\in\MR_{\geq 0}$ such that:
\begin{enumerate}
\item \label{item1}$\MR_{\geq 0}\subseteq \Omega$.
\item \label{item2}${\boldsymbol{f}}(0)=\bs{\alpha}$.
\item \label{item3}For all $t\in\MR_{\geq 0}$, for $\iton:f_i(t)\in\MR$ and $0\leq f_i(t)< k$.
\item \label{item4}For all simply-connected open sets $\widetilde{\Omega}\subseteq\MC$, for all
$\CE$-processes $\bs{\tilde{f}}$ on $\widetilde{\Omega}$, for all
intervals $I\subseteq \widetilde{\Omega}\cap\MR_{\geq 0}$, if $0\in I$
and $\bs{\tilde{f}}(0)=\alpha$, then for all $t\in I$,
${\boldsymbol{f}}(t)=\bs{\tilde{f}}(t)$.
\end{enumerate}
\end{lemma}

\begin{proof}
Claim: There exists $k\in\MR_{\geq 0}$ such that for all intervals
$I\subseteq\MR_{\geq 0}$ with $0\in I$, for all real-$\CE$-processes
$\boldsymbol{\tilde{h}}=\langle \tilde{h}_1, \tilde{h}_2, \ldots,
\tilde{h}_n\rangle$ on $I$ with $\bs{\tilde{h}}(0)=\bs{\alpha}$, for
all $t\in I$, for $\iton$: $0\leq \tilde{h}_i(t) \leq k$.

To see this, let $I\subseteq\MR_{\geq 0}$ be an interval such that
$0\in I$. Let $\boldsymbol{\tilde{h}}=\langle \tilde{h}_1,
\tilde{h}_2, \ldots, \tilde{h}_n\rangle$ be a real-$\CE$-process on
$I$ such that $\bs{\tilde{h}}(0)=\bs{\alpha}$.

From Lemma~\ref{lem:solution}, there exist an open, simply-connected
$\widetilde{\Omega}\subseteq\MC$ and an $\CE$-process
$\boldsymbol{\tilde{f}}=\langle
\tilde{f}_1,\tilde{f}_2,\ldots,\tilde{f}_n\rangle$ on
$\widetilde{\Omega}$ such that:
\begin{enumerate}
\item $I\subset\widetilde{\Omega}$
\item For all $t\in I: \boldsymbol{\tilde{f}}(t)=\boldsymbol{\tilde{h}}(t)$.
\end{enumerate}

From Lemma~\ref{bounded}, there exists $k\in\MR_{\geq0}$ such that
for all $t\in I$, for $i=1,2,\ldots,n$: $\tilde{f}_i(t)\in\MR$ and
$0\leq \tilde{f}_i(t)<k$. That is, for all $t\in I$, for $\iton:
0\leq \tilde{h}_i(t) < k$. This proves the claim.

Therefore, by~\cite[p.~397,~Corollary]{hirsch04textbook}, there exists
$k\in\MR_{\geq0}$, there is a real-$\CE$-process $\bs{h}=\langle
h_1,h_2,\ldots,h_n\rangle$ on $\MR_{\geq 0}$ such that
$\bs{h}(0)=\bs{\alpha}$ and for all $t\in \MR_{\geq 0}$,for $\iton:
0\leq h_i(t) < k$.. By Lemma~\ref{lem:solution}, there exist an
open, simply-connected $\Omega\subseteq\MC$ and an $\CE$-process
$\bs{f}$ on $\Omega$ such that $\MR_{\geq 0}\subseteq\Omega$ and for
all $t\in\MR_{\geq 0}$, $\bs{f}(t)=\bs{h}(t)$. Therefore, for all
$t\in\MR_{\geq 0}$, for $\iton: f_i(t)\in\MR$ and $0\leq f_i(t)< k$.
Hence, Parts (\ref{item1},\ref{item2},\ref{item3}) are established.
Part(\ref{item4}) follows from Lemma~\ref{lem:uni_real_global}.
\end{proof}

The next lemma shows that the $\omega$-limit points of $\CE$-processes that start at non-negative points satisfy detailed balance.

\begin{lemma}\label{limits are se}
Let $\CE$ be a finite, natural event-system of dimension $n$, let
$\Omega \subseteq \MC$ be open and simply-connected, let $\bs{f}$ be
an $\CE$-process on $\Omega$, and let $\bs{q} \in \MC^n$. If
$\MR_{\geq 0} \subseteq \Omega$ and $\bs{f}(0)$ is a non-negative
point and $\bs{q}$ is an $\omega$-limit point of $\bs{f}$, then
$\bs{q} \in \MR^n_{\geq 0}$ and is a strong $\CE$-equilibrium point.
\end{lemma}

\begin{proof}
Suppose $\MR_{\geq 0} \subseteq \Omega$, $\bs{f}(0)$ is a
non-negative point, $S$ is the set of $\omega$-limit points of
$\bs{f}$, and $\bs{q} \in S$. By Lemma~\ref{invariant} $\bs{q} \in
\MR^n_{\geq 0}$. By Lemma~\ref{lem:ex&uni} there exists an open,
simply-connected $\Omega_{\bs{q}} \subseteq \MC$ such that
$\MR_{\geq 0} \subseteq \Omega_{\bs{q}}$ and an $\CE$-process
$\bs{h} = \langle h_1, h_2, \ldots, h_n \rangle$ on
$\Omega_{\bs{q}}$ such that $\bs{h}\left(0\right) = \bs{q}$.

Let $\bs{c}$ be a positive strong $\CE$-equilibrium point. By
Lemma~\ref{ex lyap properties}.\ref{ex lyap positive},
$\overline{g_{\CE,\bs{c}}}\left(\bs{f}\left(t\right)\right)$ is
bounded below and, by Lemma~\ref{ex lyap properties}.\ref{gmon}, is
monotonically non-increasing on $\MR_{\geq 0}$. Therefore
$\lim_{t\rightarrow\infty}
\overline{g_{\CE,\bs{c}}}\left(\bs{f}\left(t\right)\right)$ exists.
Since $\overline{g_{\CE,\bs{c}}}$ is continuous, for all
$\bs{\alpha} \in S$,
$\overline{g_{\CE,\bs{c}}}\left(\bs{\alpha}\right) =
\lim_{t\rightarrow\infty}
\overline{g_{\CE,\bs{c}}}\left(\bs{f}\left(t\right)\right)$. By
Lemma~\ref{invariant}, for all $t \in \MR_{\geq 0}$, $\bs{h}(t) \in
S$. Hence,
$\overline{g_{\CE,\bs{c}}}\left(\bs{h}\left(t\right)\right)$ is
constant on $\MR_{\geq 0}$.

By Lemma~\ref{reduce} either $\bs{q}$ is a strong $\CE$-equilibrium
or there exists a finite event-system $\hat{\CE}$ of dimension
$\hat{n} \leq n$, an $\hat{\CE}$-process $\bs{\hat{h}} = \langle
\hat{h}_1, \hat{h}_2, \ldots, \hat{h}_{\hat{n}} \rangle$ on
$\Omega_{\bs{q}}$, and a permutation $\pi$ on $\{1, 2, \ldots, n\}$
satisfying 1--7 of Lemma~\ref{reduce}.

Assume $\bs{q}$ is not a strong $\CE$-equilibrium point. By
Lemma~\ref{reduce}.\ref{reduce.e}, for $i = 1, 2, \ldots, \hat{n}$,
for all $t \in \Omega_{\bs{q}}$, $\hat{h}_i\left(t\right) =
h_{\pi^{-1}\left(i\right)}\left(t\right)$. Let $\bs{\hat{c}} =
\langle \hat{c}_1, \hat{c}_2, \ldots, \hat{c}_{\hat{n}} \rangle =
\langle c_{\pi^{-1}\left(1\right)}, c_{\pi^{-1}\left(2\right)},
\ldots, c_{\pi^{-1}\left(\hat{n}\right)} \rangle$. By
Lemma~\ref{reduce}.\ref{reduce.b2}, $\bs{\hat{c}}$ is an
$\hat{\CE}$-strong equilibrium point.

For all $v \in \MR_{> 0}$, let $\overline{g_v}$ be as defined in
Equation~\ref{gv} in Definition~\ref{extended lyap}. Then for all $t
\in \MR_{\geq 0}$,
\begin{align*}
\overline{g_{\CE, \bs{c}}}\left(\bs{h}\left(t\right)\right) - \overline{g_{\hat{\CE}, \bs{\hat{c}}}}\left(\bs{\hat{h}}\left(t\right)\right) & = \sum_{i=1}^n \overline{g_{c_i}}\left(h_i\left(t\right)\right) - \sum_{j=1}^{\hat{n}} \overline{g_{\hat{c}_j}}\left(\hat{h}_j\left(t\right)\right) \\
& = \sum_{i=1}^n \overline{g_{c_i}}\left(h_i\left(t\right)\right) - \sum_{j=1}^{\hat{n}} \overline{g_{c_{\pi^{-1}\left(j\right)}}}\left(h_{\pi^{-1}\left(j\right)}\left(t\right)\right) \\
& = \sum_{i=1}^n \overline{g_{c_{\pi^{-1}\left(i\right)}}}\left(h_{\pi^{-1}\left(i\right)}\left(t\right)\right) - \sum_{j=1}^{\hat{n}} \overline{g_{c_{\pi^{-1}\left(j\right)}}}\left(h_{\pi^{-1}\left(j\right)}\left(t\right)\right) \\
& = \sum_{i=\hat{n}+1}^n
\overline{g_{c_{\pi^{-1}\left(i\right)}}}\left(h_{\pi^{-1}\left(i\right)}\left(t\right)\right)
\end{align*}

But, by Lemma~\ref{reduce}.\ref{reduce.f}, if $\pi\left(i\right) >
\hat{n}$ then $h_i\left(t\right)$ is constant. Hence,
$\overline{g_{c_{\pi^{-1}\left(i\right)}}}\left(h_{\pi^{-1}\left(i\right)}\left(t\right)\right)$
is constant for $i = \hat{n}+1,\hat{n}+2,\ldots, n$, so
$\overline{g_{\CE, \bs{c}}}\left(\bs{h}\left(t\right)\right) -
\overline{g_{\hat{\CE},
\bs{\hat{c}}}}\left(\bs{\hat{h}}\left(t\right)\right)$ is constant.
Since $\overline{g_{\CE,\bs{c}}}\left(\bs{h}\left(t\right)\right)$
and $\overline{g_{\CE, \bs{c}}}\left(\bs{h}\left(t\right)\right) -
\overline{g_{\hat{\CE},
\bs{\hat{c}}}}\left(\bs{\hat{h}}\left(t\right)\right)$ are both
constant, $\overline{g_{\hat{\CE},
\bs{\hat{c}}}}\left(\bs{\hat{h}}\left(t\right)\right)$ must be
constant. By Lemma~\ref{reduce}.\ref{reduce.d}, for all $t \in
\MR_{>0}$, $\hat{\bs{h}}\left(t\right)$ is a positive point, so by
Definitions~\ref{def:Lyapunov} and~\ref{extended lyap},
$\overline{g_{\hat{\CE},
\bs{\hat{c}}}}\left(\bs{\hat{h}}\left(t\right)\right) =
g_{\hat{\CE}, \bs{\hat{c}}}\left(\bs{\hat{h}}\left(t\right)\right)$.
Since $g_{\hat{\CE},
\bs{\hat{c}}}\left(\bs{\hat{h}}\left(t\right)\right)$ is constant,
$\frac{d}{dt} g_{\hat{\CE},
\bs{\hat{c}}}\left(\bs{\hat{h}}\left(t\right)\right) = \nabla
g_{\hat{\CE},
\bs{\hat{c}}}\left(\bs{\hat{h}}\left(t\right)\right)\cdot
\bs{P}_\CE\left(\bs{\hat{h}}\left(t\right)\right) = 0$. Then by
Theorem~\ref{SecondLaw} and continuity $\bs{\hat{h}}\left(0\right)$
must be a strong $\hat{\CE}$-equilibrium point, so for all $e \in
\hat{\CE}$, for all $t \in \Omega_{\bs{q}}$, $e(\bs{\hat{h}}(t)) =
0$, which contradicts Lemma~\ref{reduce}.\ref{reduce.c}. Therefore
$\bs{q}$ is a strong $\CE$-equilibrium point.
\end{proof}

The next theorem consolidates our results concerning natural event-systems.
It also establishes that positive strong equilibrium points are locally attractive relative to their conservation classes. Together with the existence of a Lyapunov function, this implies that positive strong equilibrium points are asymptotically stable relative to their conservation classes~\cite[Theorem~5.57]{irwin80smooth}.

\begin{theorem}\label{AsymptoticStability}
Let $\CE$ be a finite, natural event-system of dimension $n$. Let
$H$ be a positive conservation class of $\CE$. Then:
\begin{enumerate}
\item For all $\bs{x}\in H\cap\MR^n_{\geq0}$, there
exist $k\in\MR_{\geq 0}$, an open, simply-connected
$\Omega\subseteq\MC$ and an $\CE$-process ${\boldsymbol{f}}=\langle
f_1,f_2,\ldots,f_n\rangle$ on $\Omega$ such that:
\label{AsymptoticStability_XinR}
\begin{enumerate}
\item $\MR_{\geq 0}\subseteq\Omega$. \label{AsymptoticStability_RinO}
\item $\bs{f}(0)=\bs{x}$.
\item For all $t\in\MR_{\geq 0}$, $\bs{f}(t)\in H\cap\MR^n_{\geq
0}$. \label{AsymptoticStability_finHR}
\item For all $t\in\MR_{\geq 0}$, for $\iton$, $0\leq f_i(t)\leq k$.
\item For all open, simply-connected $\widetilde{\Omega}\subseteq\MC$, for all $\CE$-processes
$\bs{\tilde{f}}$ on $\widetilde{\Omega}$, if $0\in\widetilde{\Omega}$ and
$\bs{\tilde{f}}(0)=\bs{x}$ then for all intervals
$I\subseteq\widetilde{\Omega}\cap\MR_{\geq 0}$ such that $0\in I$, for
all $t\in I: \bs{f}(t)=\bs{\tilde{f}}(t)$.
\end{enumerate}
\item\label{Thm:AsymptoticStability.2} There exists $\bs{c}\in H$ such that:
\begin{enumerate}
\item\label{Thm:AsymptoticStability.2a} $\bs{c}$ is a positive strong $\CE$-equilibrium point.
\item For all $\bs{d}\in H$, if $\bs{d}$ is a positive strong
$\CE$-equilibrium point, then $\bs{d}=\bs{c}$.
\item\label{Thm:AsymptoticStability.2c} There exists $U\subseteq H\cap \MR^n_{> 0}$ such that
\begin{enumerate}
\item $U$ is open in $H\cap\MR^n_{>0}$.
\item $\bs{c}\in U$.
\item For all $\bs{x}\in U$, there exist an open,
simply-connected $\Omega\subseteq\MC$ and an $\CE$-process $\bs{f}$
on $\Omega$ such that
\begin{enumerate}
\item $\MR_{\geq 0}\subseteq\Omega$.
\item $\bs{f}(0)=\bs{x}$.
\item $\bs{f}(t)\rightarrow \bs{c}$ as $t\rightarrow\infty$ along the
positive real line. (i.e. for all $\varepsilon\in\MR_{>0}$, there
exists $t_0\in\MR_{>0}$ such that for all
$t\in\MR_{>t_0}:||\bs{f}(t)-\bs{c}||_2 < \varepsilon$.)
\end{enumerate}
\end{enumerate}
\end{enumerate}
\end{enumerate}
\end{theorem}

\begin{proof}\
\\1. Follows from Lemma~\ref{lem:ex&uni} and
Theorem~\ref{thm:consclass}.\
\\2.(a) and 2.(b) follow from Theorem~\ref{HornJackson}.\
\\2.(c) Let $\bs{c}\in H$ be a positive strong-$\CE$-equilibrium point as in Theorem~\ref{AsymptoticStability}.\ref{Thm:AsymptoticStability.2a}.
Let $g=g_{\CE,\bs{c}}$. Let $T=H\cap\MR^n_{>0}$. For all $\bs{x}\in
H\cap\MR^n$, for all $r\in\MR_{>0}$, let
\begin{align*}
B_r(\bs{x})=\big\{\bs{y}\in H\cap\MR^n\mid \|\bs{x}-\bs{y}\|_2<r\big\}\\
S_r(\bs{x})=\big\{\bs{y}\in H\cap\MR^n\mid \|\bs{x}-\bs{y}\|_2=r\big\}\\
\overline{B_r(\bs{x})}=\big\{\bs{y}\in H\cap\MR^n\mid
\|\bs{x}-\bs{y}\|_2\leq r\big\}
\end{align*}

Since $\MR^n_{>0}$ is open in $\MR^n$, it follows that $T$ is open
in $H\cap\MR^n$. Therefore, there exists $\delta\in\MR_{>0}$ such
that $B_{2\delta}(\bs{c})\subseteq T$. Let $\delta\in\MR_{>0}$ be
such that $B_{2\delta}(\bs{c})\subseteq T$. It follows that
$\overline{B_\delta(\bs{c})}\subseteq T$.

Since $g$ is continuous and $S_\delta(\bs{c})$ is compact, let
$\bs{x_0}\in S_\delta(\bs{c})$ be such that
$g(\bs{x_0})=\inf_{\bs{x}\in S_\delta(\bs{c})} g(\bs{x})$. Let
$U=B_\delta(\bs{c})\cap\{\bs{x}\in T\mid g(\bs{x})<g(\bs{x_0})\}$.
It follows that $U$ is open in $T$. Since $\bs{x_0}\neq\bs{c}$, and
by Lemma~\ref{ex lyap properties}.\ref{ex lyap positive},
$g(\bs{x_0})=\overline{g_{\CE,\bs{c}}}(\bs{x_0})>0=g(\bs{c})$.
Hence, $\bs{c}\in U$.

Let $\bs{x}\in U$. From Lemma~\ref{lem:ex&uni}, there exist an open,
simply-connected $\Omega\subset\MC$ and an $\CE$-process $\bs{f}$ on
$\Omega$ such that $\MR_{\geq 0}\subseteq\Omega$ and
$\bs{f}(0)=\bs{x}$.

We claim that for all $t\in\MR_{\geq 0}$, $\bs{f}(t)\in
B_\delta(\bs{c})$. Suppose not. Then there exists $t_0\in\MR_{\geq
0}$ such that $\bs{f}(t_0)\in S_\delta(\bs{c})$. From the definition
of $\bs{x_0}$, $g(\bs{x_0})\leq g(\bs{f}(t_0))$. Since
$\bs{f}(0)=\bs{x}\in U$, $g(\bs{f}(0))<g(\bs{x_0})$. Hence,
$g(\bs{f}(0))<g(\bs{f}(t_0))$, contradicting Lemma~\ref{ex lyap
properties}.\ref{gmon}.

To see that $\bs{f}(t)\rightarrow \bs{c}$ as $t\rightarrow\infty$
along the positive real line, suppose not. Then there exists
$\varepsilon\in\MR_{>0}$ such that $\varepsilon<\delta$ and there exists
an increasing sequence of real numbers $\{t_i\in\MR_{>0}\}_{
i\in\mathbb{Z}_{>0}}$ such that $t_i\to\infty$ as $i\to\infty$ and
for all $i$, $\bs{f}(t_i)\in \overline{B_\delta(\bs{c})}\setminus
B_\varepsilon(\bs{c})$. Since $\overline{B_\delta(\bs{c})}\setminus
B_\varepsilon(\bs{c})$ is compact, there exists a convergent
subsequence. By Definition~\ref{omega limit}, the limit of this
subsequence is an $\omega$-limit point $\bs{q}$ of $\bs{f}$ such
that $\bs{q}\in\overline{B_\delta(\bs{c})}\setminus
B_\varepsilon(\bs{c})$. From Lemma~\ref{limits are se}, $\bs{q}$ is a
strong-$\CE$-equilibrium point. Since
$\bs{q}\in\overline{B_\delta(\bs{c})}$, $\bs{q}\in T$. From
Theorem~\ref{HornJackson}, $\bs{q}=\bs{c}$. Hence, $\bs{c}\notin
B_\varepsilon(\bs{c})$, a contradiction.
\end{proof}

We have established that positive strong equilibrium points are asymptotically stable relative to their conservation classes. A stronger result would be that if an $\CE$-process starts at a positive point then it asymptotically tends to the positive strong equilibrium point in its conservation class. Such a result is related to the widely-held notion that, for systems of chemical reactions, concentrations approach equilibrium. We have been unable to prove this result. We will now state it as an open problem. This problem has a long history. It appears to have been first suggested in \cite[Lemma~4C]{horn72general}, where it was accompanied by
an incorrect proof. The proof was retracted in \cite{horn74dynamics}.

\begin{open}\label{MainTheorem}
Let $\CE$ be a finite, natural event-system of dimension $n$. Let
$H$ be a positive conservation class of $\CE$. Then
\begin{enumerate}
\item For all $\bs{x}\in H\cap\MR^n_{\geq0}$, there
exist $k\in\MR_{\geq 0}$, an open, simply-connected
$\Omega\subseteq\MC$ and an $\CE$-process ${\boldsymbol{f}}=\langle
f_1,f_2,\ldots,f_n\rangle$ on $\Omega$ such that:
\begin{enumerate}
\item $\MR_{\geq 0}\subseteq\Omega$.
\item $\bs{f}(0)=\bs{x}$.
\item For all $t\in\MR_{\geq 0}$, $\bs{f}(t)\in H\cap\MR^n_{\geq
0}$.
\item For all $t\in\MR_{\geq 0}$, for $\iton$, $0\leq f_i(t)<k$.
\item For all open, simply-connected $\widetilde{\Omega}\subseteq\MC$, for all $\CE$-processes
$\bs{\tilde{f}}$ on $\widetilde{\Omega}$, if $0\in\widetilde{\Omega}$ and
$\bs{\tilde{f}}(0)=\bs{x}$ then for all intervals
$I\subseteq\widetilde{\Omega}\cap\MR_{\geq 0}$, if $0\in I$ then for all
$t\in I: \bs{f}(t)=\bs{\tilde{f}}(t)$.
\end{enumerate}
\item
There exists $\bs{c}\in H$ such that:
\begin{enumerate}
\item $\bs{c}$ is a positive strong $\CE$-equilibrium point.
\item For all $\bs{d}\in H$, if $\bs{d}$ is a positive strong
$\CE$-equilibrium point, then $\bs{d}=\bs{c}$.
\item For all $\bs{x}\in H\cap \MR^n_{>0}$, there exist an open,
simply-connected $\Omega\subseteq\MC$ and an $\CE$-process $\bs{f}$
on $\Omega$ such that:
\begin{enumerate}
\item $\MR_{\geq 0}\subseteq\Omega$.
\item $\bs{f}(0)=\bs{x}$.
\item $\bs{f}(t)\rightarrow \bs{c}$ as $t\rightarrow\infty$ along the
positive real line. (i.e. for all $\varepsilon\in\MR_{>0}$, there
exists $t_0\in\MR_{>0}$ such that for all
$t\in\MR_{>t_0}:||\bs{f}(t)-\bs{c}||_2 < \varepsilon$.)
\end{enumerate}
\end{enumerate}
\end{enumerate}
\end{open}

In light of Theorem~\ref{AsymptoticStability}, Open~Problem~\ref{MainTheorem} is equivalent to the following statement.

\begin{open}\label{MainTheorem2}
Let $\CE$ be a finite, natural event-system of dimension $n$. Let
$\bs{x} \in \MR^n_{>0}$. Then there exists an open, simply-connected
$\Omega \subseteq \MC$, an $\CE$-process $\bs{f}$ on $\Omega$ and a
positive strong $\CE$-equilibrium point $\bs{c}$ such that:
\begin{enumerate}
\item $\MR_{\geq 0}\subseteq\Omega$.
\item $\bs{f}(0)=\bs{x}$.
\item $\bs{f}(t)\rightarrow \bs{c}$ as $t\rightarrow\infty$ along the
positive real line. (i.e. for all $\varepsilon\in\MR_{>0}$, there
exists $t_0\in\MR_{>0}$ such that for all
$t\in\MR_{>t_0}:||\bs{f}(t)-\bs{c}||_2 < \varepsilon$.)
\end{enumerate}
\end{open}

\section{Finite Natural Atomic Event-systems}\label{sec:atomic}
In this section, we settle Open~\ref{MainTheorem} in the affirmative for the case of finite, natural, atomic event-systems. The atomic hypothesis appears to be a natural assumption to make concerning systems of chemical reactions. Therefore, our result may be considered a validation of the notion in chemistry that concentrations tend to equilibrium. We will prove the following theorem:

\begin{theorem}\label{AtomicMainTheorem}
Let $\CE$ be a finite, natural, atomic event-system of dimension $n$. Let $\bs{\alpha} \in \MR^n_{>0}$. Then there exists an open, simply-connected $\Omega \subseteq \MC$, an $\CE$-process $\bs{f}$ on $\Omega$, and a positive strong $\CE$-equilibrium point $\bs{c}$ such that:
\begin{enumerate}
\item $\MR_{\geq 0} \subseteq \Omega$,
\item $\bs{f}(0) = \bs{\alpha}$, and
\item $\bs{f}(t) \rightarrow \bs{c}$ as $t \rightarrow \infty$ along the positive real line (i.e. for all $\varepsilon \in \MR_{>0}$, there exists $t_0 \in \MR_{>0}$ such that for all $t \in \MR_{>t_0} : \| \bs{f}(t) - \bs{c} \|_2 < \varepsilon$).
\end{enumerate}
\end{theorem}

It follows from Theorem~\ref{AsymptoticStability} that the point $\bs{c}$ depends only on the conservation class of $\bs{\alpha}$ and not on $\bs{\alpha}$ itself. That is, two $\CE$-processes starting at positive points in the same conservation class asymptotically converge to the same $\bs{c}$.

Implicit in the atomic hypothesis is the idea that atoms are neither created nor destroyed, but rather are conserved by chemical reactions. Our proof uses a formal analog of this idea.
Recall from Definition~\ref{atomic} that if $\CE$ is atomic then $C_\CE(M)$ contains a unique monomial from $\mathbb{M}_{A_\CE}$.

\begin{definition}
Let $\CE$ be a finite, natural, atomic event-system of dimension $n$. The {\em atomic decomposition map} $\bs{D}_\CE : \mathbb{M}_{\{X_1, X_2, \ldots, X_n\}} \rightarrow \MZ^n_{\geq 0}$ is the function $M \mapsto \langle b_1, b_2, \ldots, b_n \rangle$ such that $X_1^{b_1} X_2^{b_2} \cdots X_n^{b_n} \in C_\CE(M) \cap \mathbb{M}_{A_\CE}$.
\end{definition}

The next lemma lists some properties of the atomic decomposition map. Note that though the event-graph $G_\CE$ is directed, if $M$ and $N$ are monomials and there exists a path in $G_\CE$ from $M$ to $N$ then there also exists a path in $G_\CE$ from $N$ to $M$. Informally, this is because all events are ``reversible.''

\begin{lemma}\label{atomicdecomposition}
Let $\CE$ be a finite, natural, atomic event-system of dimension
$n$ and let $M,N\in \mathbb{M}_{\{X_1,X_2, \ldots, X_n\}}$. Then:
\begin{enumerate}
\item $\bs{D}_{\CE}(M)=\bs{D}_{\CE}(N)$ if and only if $C_\CE(M) = C_\CE(N)$.
\item $\bs{D}_\CE(MN)=\bs{D}_\CE(M)+\bs{D}_\CE(N)$.
\end{enumerate}
\end{lemma}
\begin{proof}
Let $\bs{D}=\bs{D}_\CE$.

(1)
$\bs{D}(M) = \bs{D}(N) = \langle b_1, b_2, \ldots, b_n \rangle$ if and only if $X_1^{b_1} X_2^{b_2} \cdots X_n^{b_n} \in C_\CE(M)$ and $X_1^{b_1} X_2^{b_2} \cdots X_n^{b_n} \in C_\CE(N)$. Then $C_\CE(M) = C_\CE(N)$.

(2)
Let $\bs{D}(M) = \langle b_1, b_2, \ldots, b_n \rangle$ and $\bs{D}(N) = \langle c_1, c_2, \ldots, c_n \rangle$. Then, in $G_\CE$ there is a
path from $M$ to $X_1^{b_1} X_2^{b_2} \cdots X_n^{b_n} \in \mathbb{M}_{A_\CE}$ and a path from $N$ to $X_1^{c_1} X_2^{c_2} \cdots X_n^{c_n} \in \mathbb{M}_{A_\CE}$. It follows that there is a path from $MN$ to $X_1^{b_1+c_1} X_2^{b_2+c_2} \cdots X_n^{b_n+c_n} \in \mathbb{M}_{A_\CE}$.
Hence $\bs{D}(MN) = \langle b_1 + c_1, b_2 + c_2, \ldots, b_n + c_n \rangle = \bs{D}(M) + \bs{D}(N)$.
\end{proof}

\begin{definition}
Let $\CE$ be a finite, natural, atomic event-system of dimension
$n$. For all $i \in \{1, 2, \ldots, n\}$, for all $M \in \mathbb{M}_{\{X_1, X_2, \ldots, X_n\}}$, $D_{\CE,i}(M)$ is the $i^\text{th}$ component of $\bs{D}_\CE(M)$.
\end{definition}

\begin{definition}
Let $\CE$ be a finite, natural, atomic event-system of dimension $n$. For all $i \in \{1, 2, \ldots, n\}$ the function $\kappa_{\CE,i}:\MC^n\rightarrow\MC$ is given by \[\langle z_1,z_2, \ldots, z_n\rangle\longmapsto\sum_{j=1}^{n}D_{\CE,i}(X_j)z_j.\]
\end{definition}

\begin{lemma}\label{Conservation of Atoms}
Let $\CE$ be a finite, natural, atomic event-system of dimension
$n$. Then for all $i \in \{1, 2, \ldots, n\}$, the function $\kappa_{\CE,i}$ is a
conservation law of $\CE$.
\end{lemma}
\begin{proof}
Let $m = | \CE |$, and for $\jtom$, let $\sigma_j, \tau_j \in \MR_{>0}$ and $M_j, N_j \in \mathbb{M}_\infty$ with $M_j\prec N_j$ be such that $\CE=\{\sigma_1M_1-\tau_1N_1, \ldots, \sigma_mM_m-\tau_mN_m\}$. For $\iton$, let $a_{j,i}, b_{j,i} \in \MZ_{>0}$ be such that $M_j = X_1^{a_{j,1}} X_2^{a_{j,2}} \cdots X_n^{a_{j,n}}$ and $N_j = X_1^{b_{j,1}} X_2^{b_{j,2}} \cdots X_n^{b_{j,n}}$. Let $(\gamma_{j,i})_{m \times n} = \Gamma_\CE$.

Then for $\jtom$:

\begin{align*}
\quad & \sigma_j M_j - \tau_j N_j\in\CE \\
\Rightarrow \quad & M_j\in C_\CE(N_j) & & \text{[Definition \ref{EventGraph}]} \\
\Rightarrow \quad & \bs{D}_\CE(M_j)=\bs{D}_\CE(N_j) & & \text{[Lemma \ref{atomicdecomposition}]} \\
\Rightarrow \quad & \sum_{i=1}^na_{j,i}\bs{D}_\CE(X_i)=\sum_{i=1}^nb_{j,i}\bs{D}_\CE(X_i) & & \text{[Lemma \ref{atomicdecomposition}]} \\
\Rightarrow \quad & \sum_{i=1}^n (b_{j,i}-a_{j,i})\bs{D}_\CE(X_i) = \bs{0} & & \\
\Rightarrow \quad & \sum_{i=1}^n\gamma_{j,i}\bs{D}_\CE(X_i)= \bs{0} & & \text{[Definition \ref{gamma}]}
\end{align*}
It follows that for all $j \in \{1,2,\ldots,m\}$, for all $k\in\{1,2,\ldots,n\}$,
\[\sum_{i=1}^n \gamma_{j,i} D_{\CE,k}(X_i) = 0 \]

Therefore, for all $k \in \{1,2,\ldots,n\}$, $\Gamma_\CE \cdot \langle D_{\CE,k}(X_1), D_{\CE,k}(X_2), \ldots, D_{\CE,k}(X_n) \rangle^T = \bs{0}$. Since the vector $\langle D_{\CE,k}(X_1),D_{\CE,k}(X_2), \ldots, D_{\CE,k}(X_n)\rangle^T$ is in the kernel of $\Gamma_\CE$, by Theorem~\ref{linear conservation law}, $\kappa_{\CE,k}$ is a conservation law of $\CE$.
\end{proof}

\begin{lemma}\label{eventgraphequality}
Let $\CE$ be a finite, natural event-system of dimension $n$. Let $M,N \in \mathbb{M}_\infty$ and let $\bs{q} \in \MC^n$. If $M \in C_\CE (N)$ and $\bs{q}$ is a strong $\CE$-equilibrium point and $M(\bs{q})=0$, then $N(\bs{q}) = 0$.
\end{lemma}

\begin{proof}
Let $\langle v_0, v_1 \rangle$ be an edge in $G_\CE$. Then there exist $e \in \CE$ and $\sigma, \tau \in \MR_{>0}$ and $T,U,V \in \mathbb{M}_\infty$ such that $e = \sigma U - \tau V$ and $v_0 = TU$ and $v_1 = TV$.

Assume $v_0(\bs{q}) = 0$. Then either $T(\bs{q}) = 0$ or $U(\bs{q}) = 0$. If $T(\bs{q}) = 0$ then $v_1(\bs{q}) = 0$. If $U(\bs{q}) = 0$ and $\bs{q}$ is a strong $\CE$-equilibrium point, then $e(\bs{q}) = \sigma U(\bs{q}) - \tau V(\bs{q}) = 0$, so $V(\bs{q}) = 0$. Therefore $v_1(\bs{q}) = 0$. The lemma follows by induction.
\end{proof}

We are now ready to prove Theorem~\ref{AtomicMainTheorem}.

\begin{proof}[Proof of Theorem~\ref{AtomicMainTheorem}]
Since $\bs{\alpha}$ is a positive point, it is in some positive conservation class $H$. By Theorem~\ref{AsymptoticStability}:
\begin{enumerate}
\item There exists exactly one positive strong $\CE$-equilibrium point $\bs{c} \in H$.
\item There exist an open and simply-connected $\Omega \subseteq \MC$ and an $\CE$-process $\bs{f}$ on $\Omega$ such that $\MR_{\geq 0} \subset \Omega$ and $\bs{f}(0) = \bs{\alpha}$.
\item For all $t \in \MR_{\geq 0}$, $\bs{f}(t) \in H \cap \MR^n_{\geq 0}$.
\item There exists $k\in \MR_{\geq 0}$ such that for $\iton$, for all $t\in \MR_{\geq 0}$, $f_i(t)\in \MR$ and $0 \leq f_i(t) \leq k$. \end{enumerate}

Let $\{t_j\}_{j \in \MZ_{>0}}$ be an infinite sequence of non-negative reals such that $t_j \rightarrow \infty$ as $j \rightarrow \infty$. Then $\{\bs{f}(t_j)\}_{j\in\MZ_{>0}}$ is an infinite sequence contained in a compact subset of $\MR^n$, so it must have a convergent subsequence. Let $\bs{q} = \langle q_1, q_2, \ldots, q_n \rangle \in \MC^n$ be the limit point of a convergent subsequence of $\{\bs{f}(t_j)\}_{j\in\MZ_{>0}}$. $H$ and $\MR^n_{\geq 0}$ are both closed in $\MC^n$, so $\bs{q} \in H \cap \MR^n_{\geq 0}$. Since $\CE$ is natural and $\bs{q}$ is an $\omega$-limit of $\bs{f}$, $\bs{q}$ must be a strong $\CE$-equilibrium point by Lemma~\ref{limits are se}.

Assume, for the sake of contradiction, that $\bs{q} \notin \MR^n_{>0}$. Let $i \in \{1, 2, \ldots, n\}$ be such that $q_i=0$. Let $N \in C_\CE(X_i) \cap \mathbb{M}_{A_\CE}$. Since $\CE$ is atomic, a unique such $N$ exists.  It follows from the definition of event graph that $X_i \in C_{\CE}(N)$. By Lemma~\ref{eventgraphequality}, $N(\bs{q}) = X_i(\bs{q}) = q_i = 0$. It follows that $N \neq 1$. Hence, there exists $X_a \in A_\CE$ such that $X_a$ divides $N$ and $X_a(\bs{q}) = 0$.

For all $j\in\{1,2,\ldots,n\}$ such that $D_{\CE,a}(X_j) \neq 0$, let $M_j \in C_\CE(X_j) \cap \mathbb{M}_{A_\CE}$. Then $X_a$ divides $M_j$, so $M_j(\bs{q}) = 0$. Again by Lemma~\ref{eventgraphequality}, $X_j(\bs{q}) = M_j(\bs{q}) = 0$, so $q_j = 0$. It follows that for all $j \in \{1,2,\ldots,n\}$ either $D_{\CE,a}(X_j) = 0$ or $q_j = 0$ so
\[ \kappa_{\CE,a}(\bs{q}) = \sum_{j=1}^n D_{\CE,a}(X_j) q_j = 0. \]
Since $\kappa_{\CE,a}$ is a conservation law of $\CE$ by Lemma~\ref{Conservation of Atoms} and $\bs{q}$ is an $\omega$-limit point of $\bs{f}$, it follows that
\begin{equation}\label{kappa zero}
\kappa_{\CE,a}(\bs{\alpha}) = 0.
\end{equation}

For all $j$, $D_{\CE,a}(X_j)$ is nonnegative, and $\bs{\alpha}$ is a positive point, so for all $j \in \{1, 2, \ldots, n\}$, $D_{\CE,a}(X_j) \alpha_j \geq 0$. But $D_{\CE,a}(X_a)=1$ and $\alpha_a > 0$ so $\kappa_{\CE,a}(\bs{\alpha}) > 0$, contradicting equation~(\ref{kappa zero}). Therefore $\bs{q} \in \MR^n_{>0}$. Since $\bs{c}$ is the unique positive strong $\CE$-equilibrium point in $H$, $\bs{c} = \bs{q}$.

Let $U \subseteq H \cap \MR^n_{>0}$ be the open set stated to exist in Theorem~\ref{AsymptoticStability}.\ref{Thm:AsymptoticStability.2c}. Since $\bs{c}$ is an $\omega$-limit point of $\bs{f}$, there exists $t_0 \in \MR_{>0}$ such that $\bs{f}(t_0) \in U$. Again by Theorem~\ref{AsymptoticStability}, there exist $\widetilde{\Omega} \subseteq \MC$ and an $\CE$-process $\bs{\tilde{f}}$ on $\widetilde{\Omega}$ such that $\MR_{\geq 0} \subseteq \widetilde{\Omega}$ and $\bs{\tilde{f}}(0) = \bs{f}(t_0)$ and $\bs{\tilde{f}}(t) \rightarrow \bs{c}$ as $t \rightarrow \infty$. By Lemma~\ref{time invariance}, for all $t \in \MR_{\geq 0}$, $\bs{f}(t+t_0) = \bs{\tilde{f}}(t)$. Therefore, $\bs{f}(t) \rightarrow \bs{c}$ as $t \rightarrow \infty$.
\end{proof} 
\section{Conclusion}\label{conclusion}
We have endeavored to place the kinetic theory of chemical reactions on a firm mathematical foundation and to make the law of mass action available for purely mathematical consideration.

With regard to chemistry, we have proven that many of the expectations acquired through empirical study are warranted. In particular:
\begin{enumerate}
\item For finite event-systems, the stoichiometric coefficients determine conservation laws that processes must obey~(Theorem~\ref{thm:consclass}). In fact, we can show (manuscript in preparation):
    \begin{enumerate}
    \item For finite, physical event-systems, the stoichiometric coefficients determine all linear conservation laws;
    \item For finite, natural event-systems, the stoichiometric coefficients determine \emph{all} conservation laws.
    \end{enumerate}
\item For finite, physical event-systems, a process begun with positive (non-negative) concentrations will retain positive (non-negative) concentrations through forward real time where it is defined~(Theorem~\ref{thm:stays_positive}). For finite, natural event-systems, a process begun with positive (non-negative) concentrations will retain positive (non-negative) concentrations through \emph{all} forward real time~(Theorem~\ref{AsymptoticStability}) --- that is, it will be defined through all forward real time.
\item Finite, natural event-systems must obey the ``second law of thermodynamics''~(Theorem~\ref{LyapunovExists}). In addition, the flow of energy is very restrictive --- finite, natural event-systems can contain no energy cycles~(Theorem~\ref{NatL}).
\item For finite, natural event-systems, every positive conservation class contains exactly one positive equilibrium point. This point is a strong equilibrium point and is asymptotically stable relative to its conservation class~(Theorem~\ref{AsymptoticStability}).
\end{enumerate}

Unfortunately, we, like our predecessors, are unable to settle the problem of whether a process begun with positive concentrations must approach equilibrium. We consider this the fundamental open problem in the field~(Open~Problem~\ref{MainTheorem}). For finite, natural event-systems that obey a mathematical analogue of the atomic hypothesis, we settle Open~Problem~\ref{MainTheorem} in the affirmative~(Theorem~\ref{AtomicMainTheorem}). In particular, we show that for finite, natural, atomic event-systems, every positive conservation class contains exactly one non-negative equilibrium point. This point is a positive strong equilibrium point and is globally stable relative to the intersection of its conservation class with the positive orthant.

In terms of expanding the mathematical aspects of our theory, there are several potentially fruitful avenues including:
\begin{enumerate}
\item \textbf{Complex-analytic aspects of event-systems.} While we exploit some of the complex-analytic properties of processes in this paper, we believe that a deeper investigation along these lines is warranted. For example, if we do not restrict the domain of a process to be simply-connected, then each component of a process becomes a complete analytic function in the sense of Weierstrass.
\item \textbf{Infinite event-systems.} Issues of convergence arise when considering infinite event-systems. To obtain a satisfactory theory, some constraints may be necessary. For example, a bound on the maximum degree of events may be worth considering. It may also be possible to generalize the notion of an atomic event-system to the infinite-dimensional case in such a way that each atom has an associated conservation law. One might then restrict initial concentrations to those for which each conservation law has a finite value. Additional constraints are likely to be needed as well.
\item \textbf{Algebraic-geometric aspects of event-systems.} Every finite event-system that generates a prime ideal has a corresponding affine toric variety~(as defined in~\cite[p.\hspace{2pt}15]{eisenbud1994binomial}). The closed points of this variety are the strong equilibria of the event-system. Further, every affine toric variety is isomorphic to an affine toric variety whose ideal is generated by a finite event system. One could generalize event-systems to allow irreversible reactions. In that case, it appears that the prime ideals generated by such event-systems are exactly the ideals corresponding to affine toric varieties.

    We can show (proof not provided) that finite, natural, atomic event-systems generate prime ideals. We are working towards settling Open~Problem~\ref{MainTheorem} in the affirmative for every finite, natural event-system that generates a prime ideal.
\end{enumerate} 
\section{Acknowledgements}\label{sec:ack}
This work benefitted from discussions with many people, named here in alphabetical order: Yuliy Baryshnikov, Yuriy Brun, Qi  Cheng, Ed  Coffman, Ashish Goel, Jack Hale, Lila
Kari, David Kempe, Eric Klavins, John Reif, Paul
Rothemund, Robert Sacker, Rolfe Schmidt, Bilal Shaw, David Soloveichik, Hal Wasserman, Erik
Winfree. 
\bibliographystyle{amsplain}
\bibliography{EventSystems}

\providecommand{\bysame}{\leavevmode\hbox to3em{\hrulefill}\thinspace}
\providecommand{\MR}{\relax\ifhmode\unskip\space\fi MR }
\providecommand{\MRhref}[2]{%
  \href{http://www.ams.org/mathscinet-getitem?mr=#1}{#2}
}
\providecommand{\href}[2]{#2}
\begin{thebibliography}{10}

\bibitem{adleman99toward}
Leonard Adleman, \emph{Toward a mathematical theory of self-assembly}, Tech.
  Report 00-722, University of Southern California, October 1999, Department of
  Computer Science.

\bibitem{AH79}
Lars Ahlfors, \emph{Complex analysis}, International Series in Pure and Applied
  Mathematics, McGraw-Hill, 1979.

\bibitem{bernstein99nonnegative}
D.~S. Bernstein and Santosh~P. Bhat, \emph{Nonnegativity, reducibility, and
  semistability of mass action kinetics}, IEEE Conference on Decision and
  Control, IEEE Publications, December 1999, pp.~2206--2211.

\bibitem{eisenbud1994binomial}
David Eisenbud and Bernd Sturmfels, \emph{Binomial ideals}, Duke Math. J.
  \textbf{84} (1996), no.~1, 1--45.

\bibitem{feinberg95existence}
Martin Feinberg, \emph{The existence and uniqueness of steady states for a
  class of chemical reaction networks}, Arch. Rational Mech. Anal. \textbf{132}
  (1995), 311--370.

\bibitem{gatermann02family}
Karin Gatermann and Birkett Huber, \emph{A family of sparse polynomial systems
  arising in chemical reaction systems}, Journal of Symbolic Computation
  \textbf{33} (2002), no.~3, 275--305.

\bibitem{Waage}
Cato~M. Guldberg and Peter Waage, \emph{Studies concerning affinity}, Journal
  of chemical education \textbf{63} (1986), 1044.

\bibitem{horn74dynamics}
Friedrich J.~M. Horn, \emph{The dynamics of open reaction systems},
  Mathematical aspects of chemical and biochemical problems and quantum
  chemistry (New York), Proc. SIAM-AMS Sympos. Appl. Math., vol. VIII, 1974.

\bibitem{horn72general}
Friedrich J.~M. Horn and Roy Jackson, \emph{General mass action kinetics},
  Arch. Rational Mech. Anal. \textbf{49} (1972), 81--116.

\bibitem{irwin80smooth}
M.~C. Irwin, \emph{Smooth dynamical systems}, Academic Press, 1980.

\bibitem{edm87}
K.~Ito (ed.), \emph{Encyclopedic dictionary of mathematics}, MIT Press,
  Cambridge, Massachusetts, 1987.

\bibitem{hirsch04textbook}
Robert L.~Devaney Morris W.~Hirsch, Stephen~Smale, \emph{Differential
  equations, dynamical systems, and an introduction to chaos}, 2nd ed.,
  Elsevier Academic Press, Amsterdam, 2004.

\bibitem{narayan03textbook}
Shanti Narayan and P.~K. Mittal, \emph{A textbook of matrices}, 10th ed., S.
  Chand and Company Ltd., New Delhi, 2003.

\bibitem{sontag01structure}
Eduardo~D. Sontag, \emph{Structure and stability of certain chemical networks
  and applications to the kinetic proofreading model of {T}-cell receptor
  signal transduction}, IEEE Trans. Automatic Control \textbf{46} (2001),
  1028--1047.

\end{thebibliography}
\end{document}